%% file: draft.tex
\documentclass{elsarticle}
\usepackage[margin=2.5cm]{geometry}
\usepackage{standalone}
\usepackage{tikz}
\usepackage{amssymb,amsmath,amsfonts,amsthm,mathrsfs,mathtools}
\usepackage{multirow}
\usepackage{tabularx}
\usepackage{diagbox}
\usepackage{indentfirst}
\usepackage{float}
\usepackage{array}
\usepackage{booktabs}
\usepackage{siunitx}
\usepackage{colortbl,xcolor}

\usepackage{graphicx,subcaption}
\usepackage{lineno,hyperref}
\usepackage{array,calc}
\usepackage{pgfplots}
\usepackage{pgfplotstable}
\usepackage{varwidth}
\usepackage{enumitem}
\usepackage{placeins}

\usepackage[linesnumbered,algoruled,boxed,lined]{algorithm2e}
\usepackage{pgfplots}
\usepackage{pgfplotstable}
\usepackage{mystyle}
\usepackage{newtxtext,newtxmath}
\usepackage{blkarray}

\usepackage{arydshln}
\usepackage{subcaption} 

\DeclareMathOperator*{\A}{ \mathlarger{\mathlarger{\mathlarger{\boldsymbol{\mathsf{A}}}}} }

\DeclareMathOperator{\Ima}{Im}

\pgfplotsset{compat=1.14}
\modulolinenumbers[5]

\journal{Journal of \LaTeX\ Templates}

\biboptions{numbers,sort&compress}

\begin{document}

\begin{frontmatter}

	\title{Isogeometric B\'ezier dual mortaring: The biharmonic problem}

	\author[byu_address]{Di Miao\corref{mycorrespondingauthor}}
	\ead{miaodi1987@gmail.com}

	\author[byu_address,coreform_address]{Michael A. Scott}

	\author[coreform_address]{Michael J. Borden}

	\author[coreform_address]{Derek C. Thomas}

	\author[byu_address]{Zhihui Zou}

	\cortext[mycorrespondingauthor]{Corresponding author}

	\address[byu_address]{Department of Civil and Environmental Engineering, Brigham Young University, 368 CB, Provo, UT 84602, USA}
	\address[coreform_address]{Coreform LLC, P.O. Box 970336, Orem, UT 84097, USA}

	\begin{abstract}
		In this paper we develop an isogeometric B\'ezier dual mortar method for the biharmonic problem on multi-patch domains. The well-posedness of the discrete biharmonic problem requires a discretization with $C^1$ continuous basis functions. Hence, two Lagrange multipliers are required to apply both $C^0$ and $C^1$ continuity constraints on each intersection. The dual mortar method utilizes dual basis functions to discretize the Lagrange multiplier spaces. In order to preserve the sparsity of the coupled problem, we develop a dual mortar suitable $C^1$ constraint and utilize the \Bezier dual basis to discretize the Lagrange multiplier spaces. The \Bezier dual basis functions are constructed through \Bezier projection and possess the same support size as the corresponding B-spline basis functions. We prove that this approach leads to a well-posed discrete problem and specify requirements to achieve optimal convergence. Although the \Bezier dual basis is sub-optimal due to the lack of polynomial reproduction, our formulation successfully postpones the domination of the consistency error for practical problems. We verify the theoretical results and demonstrate the performance of the proposed formulation through several benchmark problems. 
	\end{abstract}


\end{frontmatter}

\section{Introduction}
Isogeometric analysis, introduced by Hughes et al.~\cite{HUGHES20054135}, leverages computer aided design (CAD) representations directly in finite element analysis. It has been shown that this approach can alleviate the model preparation burden of going from a CAD design to an analysis model and improve overall solution accuracy and robustness~\cite{bazilevs2006isogeometric, da2011some, da2014mathematical}. Additionally, the higher-order smoothness inherent in CAD basis functions make it possible to solve higher-order partial differential equations, e.g. the biharmonic equation~\cite{kapl_isogeometric_2015, kapl_isogeometric_2017}, the Kirchhoff-Love shell problem~\cite{kiendl2009isogeometric, kiendl2010bending, kiendl2015isogeometric} and the Cahn-Hilliard equation~\cite{gomez2008isogeometric, borden2014higher} directly without resorting to complex mixed discretization schemes.\par

\begin{figure}[ht]
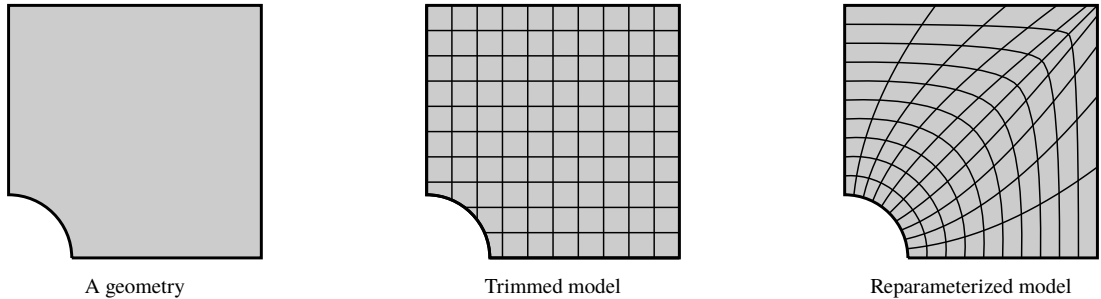

	\captionsetup[subfigure]{labelformat=empty, font = footnotesize}
	\centering
	\begin{subfigure}[b]{0.33\textwidth}
		\centering
		\includestandalone[scale=.7]{geometry}
		\caption{A geometry}
	\end{subfigure}
	\begin{subfigure}[b]{0.33\textwidth}
		\centering
		\includestandalone[scale=.7]{trimmed_geometry}
		\caption{Trimmed model}
	\end{subfigure}
	\begin{subfigure}[b]{0.33\textwidth}
		\centering
		\includestandalone[scale=.7]{reparameterized_geometry}
		\caption{Reparameterized model}
	\end{subfigure}
	\caption{A geometry and two modelling strategies: trimming and reparameterization.}
	\label{fig:geometries}
\end{figure}

CAD models are often built from collections of non-uniform rational B-splines (NURBS). Adjacent NURBS patches often have inconsistent knot layouts, different parameterizations, and may not even be physically connected. Additionally, trimming curves~\cite{kim2009isogeometric, schmidt2012isogeometric} are often employed to further simplify the design process and to extend the range of objects that can be modeled by NURBS at the expense of further complicating the underlying parameterization of the object. While usually not an issue from a design perspective, these inconsistencies in the NURBS patch layout, including trimming, must be accommodated in the isogeometric model to achieve accurate simulation results. As shown in Figure~\ref{fig:geometries}, two primary approaches are often employed. First, the exact trimmed CAD model, shown in Figure~\ref{fig:geometries} in the middle, is used directly in the simulation~\cite{schmidt2012isogeometric}. To accomplish this requires additional algorithms for handling cut cells and the weak imposition of boundary conditions and may result in reduced solution accuracy and robustness. Second, the CAD model is reparameterized~\cite{xu2014high}, as shown in Figure~\ref{fig:geometries} on the right, into a watertight spline representation like multi-patch NURBS, subdivision surfaces~\cite{peters2008subdivision}, or T-splines~\cite{sederberg_t-splines_2003} which can then be used as a basis for analysis directly. The reparameterization process often results in more accurate and robust simulation results but is only semi-automatic using prevailing approaches. In both cases, existing techniques are primarily surface-based due to the predominance of surface-based CAD descriptions.

In the present work, we assume that some form of reparameterization has been performed (see Figure~\ref{fig:geometries} on the right) on a CAD model to either remove some or all of the trimming curves and/or to restructure the underlying patch layout to improve the parameterization (e.g., reduce degree, distortion, complexity, etc.). However, we relax the requirement that adjacent patches share a consistent parameterization along shared edges and instead introduce a dual mortaring along the shared interfaces. Relaxing this requirement can simplify the reparameterization process leading to more robust approaches~\cite{xu2011parameterization}. This dual mortaring can be built into the simulation technology directly or can be used to build a weakly continuous basis which can then be used for either design or analysis. The present work leverages B\'ezier dual mortaring along each patch interface. In particular, in this work, the B\'ezier dual mortaring approach, introduced in~\cite{zou2018isogeometric}, is extended to biharmonic problems, which require the weak satisfaction of $C^1$ continuity. To preserve the sparsity of the condensed linear system, we propose a dual mortar suitable $C^1$ constraint and the corresponding Lagrange multiplier. Several different treatments of extraordinary points are considered. Numerical benchmarks illustrate the accuracy and robustness of the proposed method for both $2^\text{nd}$ and $4^\text{th}$ order problems.\par

The remainder of the paper is organized as follows: In Section~\ref{sec:prior}, we discuss prior work. In Section~\ref{sec:dual}, we give a brief review of dual basis functions. The dual mortar method is illustrated in Section~\ref{sec:dual-mortar-method}. The dual mortar suitable $C^1$ constraint and corresponding dual mortar formulation is given in Section~\ref{sec:dual-mortar-formulation}. In Section~\ref{sec:error_analysis}, we study the approximation error of the proposed method. Numerical results for several $2^\text{nd}$ order and $4^\text{th}$ order PDEs are given in Section~\ref{sec:numerical_examples}, followed by conclusions in Section~\ref{sec:conclusion}.

\section{Prior work}
\label{sec:prior}

Researchers in both the design and analysis communities have made significant progress in handling multi-patch NURBS models and the connections between adjacent patches. In the design community, some of the earliest efforts resulted in the so-called subdvision surfaces~\cite{peters2008subdivision,catmull1978recursively,10000144165} which allow for the construction of smooth spline bases over unstructured meshes. However, from an analysis perspective, subdivision surfaces are non-polynomial and require expensive numerical integration schemes~\cite{nguyen2014comparative} to achieve accuracy. T-splines~\cite{sederberg_t-splines_2003}, also introduced in the design community, overcome some of the limitations of subdivision surfaces, and introduce the possibility of local refinement. One of the most difficult aspects of building smooth CAD models from multi-patch NURBS objects is achieving smoothness and approximability in the neighborhood of extraordinary points~\cite{toshniwal2017smooth}. In a surface, an extraordinary point is a vertex in the mesh which has more or less than four adjacent edges. The concept of geometric smoothness has been used extensively in the construction of smooth surfaces over extraordinary points~\cite{farin2002handbook} and some of these concepts have been utilized in the context of isogeometric analysis~\cite{scott2013isogeometric, peters1992joining, peters1994constructing}. Loop~\cite{pan2015isogeometric} and Catmull-Clark~\cite{wei2015truncated, wei2016extended, burkhart2010iso} subdivision has been utilized in isogeometric analysis to generate smooth surfaces and solids. T-spline-based isogeometric analysis was first proposed in~\cite{bazilevs_isogeometric_2010} and the approximation properties of analysis-suitable T-splines was studied in~\cite{da2011isogeometric}. T-splines have been widely applied in different areas including fracture~\cite{verhoosel2011isogeometric}, boundary element analysis~\cite{scott2013isogeometric}, fluid-structure interaction~\cite{bazilevs2006isogeometric1} and shells~\cite{benson2010generalized}. In~\cite{kapl_isogeometric_2015}, geometric continuity was used to construct $C^1$ smooth functions on two-patch planar domains and multi-patch planar domains~\cite{kapl_isogeometric_2017}. A possible issue with this approach is so-called $C^1$ locking~\cite{collin_analysis-suitable_2016}. In~\cite{chan2018isogeometric}, a local degree elevation approach is proposed to overcome this form of locking and is applied to the construction of geometrically continuous functions on smooth surfaces.\par




From the analysis perspective, the pointwise satisfaction of continuity constraints between adjacent patches is often unnecessarily rigorous. A reasonable approximation can be achieved even if these constraints are applied in a variational setting. The Lagrange multiplier method is a general framework which can be used to apply constraints to variational problems. In the context of isogeometric analysis, various types of Lagrange multiplier approaches have been applied to problems in solids~\cite{hesch_isogeometric_2012, seitz2016isogeometric} and fluids~\cite{bazilevs2012isogeometric}. While general in applicability, the solvability and optimality of the Lagrange multiplier method is significantly influenced by the \textit{inf-sup} condition~\cite{babuvska1973finite,boffi_mixed_2013}. In the context of domain coupling, to satisfy the \textit{inf-sup} condition, there are special considerations when building the Lagrange multiplier space to ensure stability~\cite{barbosa1991finite}. This has been further studied in~\cite{bernardi_basics_2005, bernardi_domain_1993, belgacem_mortar_1998} for finite element analysis and in~\cite{brivadis_isogeometric_2015} for isogeometric analysis. Whereas the Lagrange multiplier method leads to a saddle point problem, the mortar method, first introduced by Bernardi~\cite{bernardi_domain_1993}, considers a constrained solution space and gives rise to a positive definite variational problem. Wohlmuth~\cite{wohlmuth2000mortar} used dual basis functions to discretize the Lagrange multiplier spaces, which further simplifies the mortar formulation. A dual mortar method for isogeometric analysis was first developed by Seitz et al.~\cite{seitz_isogeometric_2016}. Other variational solutions include the Nitsche method~\cite{riviere2008discontinuous, nguyen_nitsches_2014, guo_nitsches_2015} and the penalty method~\cite{apostolatos2015domain, kiendl2010bending}. The stability parameter in the Nitsche method needs to be approximated by eigenvalue problems associated with element intersections, which increases the computational cost. The penalty method results in a variational problem that is inconsistent with the original problem. Hence, the convergence to the exact solution is not guaranteed. \par


\section{Dual basis functions}\label{sec:dual}

In this section, we give a brief introduction to the concept of global and \Bezier dual basis functions. \Bezier dual basis functions will be used in Section~\ref{sec:dual-mortar-method} to facilitate the solution of domain coupling problems in the dual mortar method. A dual basis is defined as a set of basis functions $\{\hat{N}_i\}_{i=1}^{n}$, which are dual to a corresponding set of primal basis functions $\{{N}_i\}_{i=1}^{n}$ in the sense that
\begin{equation}\label{eq:def:dual}
	\langle\hat{N}_i,N_j\rangle_\Omega:=\int_\Omega\hat{N}_iN_jd\Omega=\delta_{ij}, \quad\forall{}i,j\in\left[1,2,\dots,n\right],
\end{equation}
where $\delta_{ij}$ is the Kronecker delta.

\subsection{Global dual basis}
The global dual basis functions  $\{\hat{N}_i^G\}_{i=1}^{n}$ for a given set of primal basis function $\{{N}_i\}_{i=1}^{n}$ can be computed as
\begin{equation}
	\hat{N}_i^G=\sum_j G^{-1}_{ij}N_j,\label{eq:global_dual}
\end{equation}
where  $G^{-1}_{ij}$ are the components of the inverse of the Gramian matrix $\mathbf{G}$ with components $G_{ij}=\langle{N}_i,N_j\rangle_\Omega$.

\begin{figure}[ht]
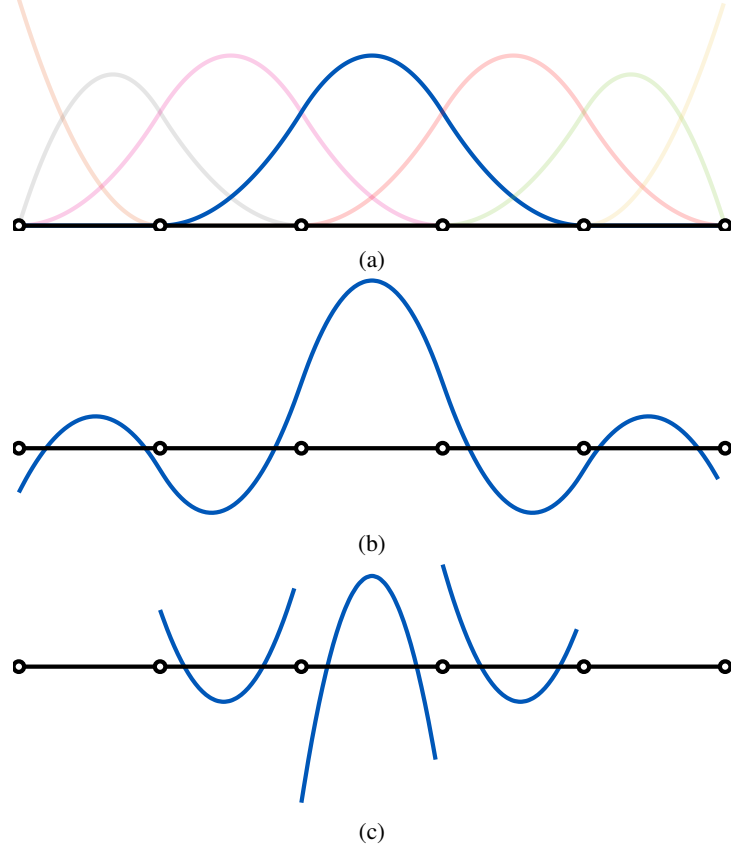

	\center
	\begin{subfigure}{\linewidth}
		\center
		\includestandalone[scale=1]{basisfunctions}
		\caption{}\label{fig:b-spline-func}
	\end{subfigure}
	\begin{subfigure}{\linewidth}
		\center
		\includestandalone[scale=1]{dual_global_one}
		\caption{}\label{fig:global_dual}
	\end{subfigure}
	\begin{subfigure}{\linewidth}
		\center
		\includestandalone[scale=1]{bezier_dual_one}
		\caption{}\label{fig:bezier_dual}
	\end{subfigure}
	\caption{A comparison of a B-spline basis function (a) with the corresponding global dual basis function (b) and the \Bezier dual basis function (c).}
	\label{fig:bezier_extraction_illustration}
\end{figure}

In this work, we choose B-spline functions as the primal basis. One important property of B-spline functions is that they have compact support. This leads to sparse linear systems when these functions are used to define the trial and weighting function spaces in a Galerkin method. The global dual basis functions, however, do not have compact support and will result in dense linear systems when used as the weighting function space in a Galerkin method. Dual basis supports are shown in Figure~\ref{fig:bezier_extraction_illustration} where we have highlighted one B-spline function in Figure~\ref{fig:b-spline-func} and then shown the corresponding global dual basis function in Figure~\ref{fig:global_dual}.

\subsection{\Bezier dual basis}

To maintain the sparsity of linear systems we will use \Bezier dual basis functions, which are computed locally and have compact support. These functions are computed using the \Bezier projection operator introduced in~\cite{thomas_bezier_2015}. The \Bezier dual basis has been shown to be effective in reducing volumetric and shear locking~\cite{MIAO2018273}, alleviating membrane locking in Kirchhoff-Love shells~\cite{greco2018reconstructed}, and as a dual mortaring strategy for elasticity problems~\cite{zou2018isogeometric}.

The construction of \Bezier dual basis functions leverages \Bezier extraction~\cite{borden2011isogeometric} and can be performed in several simple steps. For a given \Bezier element, $\Omega_e$, the element extraction operator $\mathbf{C}^e$ is computed. The element extraction operator maps a set of Bernstein polynomials $\{B_i\}_{i=1}^m$ defined over a \Bezier element, where $m$ depends on the polynomial degree of the \Bezier element in each parametric direction, into the set of global B-spline basis functions that have support over that element. The element reconstruction operator, $\mathbf{R}^e = \mathbf{C}^{e-1}$, and the Gramian matrix, $\mathbf{G}^{e}$, of the Bernstein polynomials defined over the element are then computed. The element extraction operator for the dual basis is then simply
\begin{equation}
	\hat{\mathbf{D}}^e=\text{diag}(\omega^e)\mathbf{R}^e\left[\mathbf{G}^{e}\right]^{-1}
\end{equation}
where $\text{diag}(\omega^e)$ is a diagonal matrix that contains the \Bezier projection weights
\begin{equation}
	\omega_i^e=\frac{\int_{\Omega^e}N_id\Omega}{\int_{\Omega}N_id\Omega}.
\end{equation}
The restriction of a \Bezier dual basis functions $\hat{N}_i^B$ to $\Omega_e$ is then computed as
\begin{equation}
	\hat{N}^B_i\vert_{\Omega_e}=\sum_{j=1}^m\hat{D}^e_{ij}B_j.
\end{equation}
From this local definition of the dual basis over an element we have
\begin{equation}
	\int_{\Omega_e}\hat{N}^B_iN_jd\Omega=\omega^e_i\delta_{ij},
\end{equation}
and
\begin{equation}
	\A_e{} \int_{\Omega_e}\hat{N}^B_iN_jd\Omega=\delta_{ij},
\end{equation}
where $\A$ is the standard assembly operator \cite{Hug00} . The \Bezier dual basis of the B-spline basis function highlighted in Figure~\ref{fig:b-spline-func} is shown in Figure~\ref{fig:bezier_dual}. Note that the \Bezier dual basis function has the same compact support as the primal B-spline basis function.

\begin{remark}
	The \Bezier dual basis functions define a quasi-interpolation operator $\mathcal{T}(f)=\sum_i \langle\hat{N}^B_i,f\rangle N_i$, which possesses the following properties:
	\begin{itemize}
		\item Optimal approximation: for $p^\text{th}$ order spline basis function and $f\in C^\infty$, the approximation error is given by~\cite{thomas_bezier_2015}
		      \begin{equation}
			      \|\mathcal{T}(f)-f\|_{L^2}\leq Ch^{p+1}\|f\|_{H^{p+1}}.
		      \end{equation}
		\item Boundary interpolation: for two sets of $p^\text{th}$ order spline basis functions $\{{N^s_i}\}_{i=1}^{n_s}$ and $\{{N^m_i}\}_{i=1}^{n_m}$ defined on $\left[{0,L}\right]$, if the first and last elements of $s$ are subsets of the first and last elements of $m$, then
		      \begin{equation}
			      \mathcal{T}^s(f^m)(0)=f^m(0)\;\text{ and }\;\mathcal{T}^s(f^m)(L)=f^m(L),\quad\forall f^m\in\spn\{{N^m_i}\}_{i=1}^{n_m}.\label{eq:boundary_interpolation}
		      \end{equation}
	\end{itemize}
	The second property is critical for the coercivity of the biharmonic problem on multi-patch domains.
\end{remark}

\section{The dual mortar method}
\label{sec:dual-mortar-method}
We introduce the dual mortar method in the context of an abstract formulation for a constrained problem: find $u\in\mathcal{X}$ and $\lambda\in\mathcal{M}$ such that
\begin{equation}\label{eq:LM-form}
	\left\{\begin{alignedat}{2}
		a(v,u)+b(v,\lambda)&=l(v)\quad &&\forall{}v\in\mathcal{X},\\
		b(\mu,u)&=0\quad &&\forall{}\mu\in\mathcal{M},
	\end{alignedat}\right.
\end{equation}
where $a(\cdot,\cdot)$ is a bilinear form representing a potential energy, $l(\cdot)$ is a linear form representing the external load and $b(\cdot,\cdot)$ is a bilinear form representing a set of constraints on the solution $u$. In Section~\ref{sec:dual-mortar-formulation}, $b(\cdot,\cdot)$ will represent the continuity constraints across patch boundaries.

If we introduce a pair of discrete function spaces $\mathcal{X}^h \subset \mathcal{X}$ and $\mathcal{M}^h \subset \mathcal{M}$ we can represent the weak form~\eqref{eq:LM-form} as the matrix problem
\begin{equation}\label{eq:disc-LM-form}
	\mathbf{K}^\text{LM}\mathbf{U}^{\text{LM}}=\begin{bmatrix}
		\mathbf{K} & \mathbf{B}^T \\
		\mathbf{B} & \mathbf{0}
	\end{bmatrix}\mathbf{U}^{\text{LM}}
	=
	\begin{bmatrix}
		\mathbf{F} \\
		\mathbf{0}
	\end{bmatrix},
\end{equation}
where $\mathbf{K}$ is the discretized stiffness matrix, $\mathbf{F}$ is the discretized external force vector, $\mathbf{B}$ is the discretized constraints matrix and $\mathbf{U}^{\text{LM}}$ is a vector containing the control values of the displacement field $u^h \in \mathcal{X}^h$ and Lagrange multiplier field $\lambda^h \in \mathcal{M}^h$. The stiffness matrix $\mathbf{K}^\text{LM}$ for the discrete problem~\eqref{eq:disc-LM-form} always contains both positive and negative eigenvalues, for which iterative methods are known to be less efficient than for symmetric positive definite systems. More importantly, the constraint matrix $\mathbf{B}$ might be row-wise linearly dependent if the Lagrange multiplier space $\mathcal{M}^h$ is not discretized with caution. This will lead to a non-invertible linear system. The mortar method resolves these issues by introducing a constrained function space
\begin{equation}
	\mathcal{K}\coloneq\{u\in\mathcal{X}\, | \, b(\lambda,u)=0, \quad \forall{}\lambda\in\mathcal{M}\}.
\end{equation}
The saddle point problem~\eqref{eq:LM-form} can now be transformed into a minimization problem: find $u\in\mathcal{K}$ such that
\begin{equation}
	a(v, u)=l(v),\quad \forall{}v\in\mathcal{K}.
\end{equation}
Given $\mathbf{N}$, the vector containing the basis functions of $\mathcal{X}^h$, the vector containing the basis functions of $\mathcal{K}^h$ is given by
\begin{equation}
	\mathbf{N}^k=\mathbf{C}^T\mathbf{N},\label{eq:basis-null-space}
\end{equation}
where the matrix $\mathbf{C}$ is the vector basis of the null space of the constraint matrix $\mathbf{B}$. If the Lagrange multiplier space is discretized by a set of dual basis functions, the constraint matrix $\mathbf{B}$ can be written as~\cite{gilbert1987computing}
\begin{equation}\label{eq:constraint-form}
	\mathbf{B}=\begin{bmatrix}
		\mathbf{B}_1 & \mathbf{B_2}
	\end{bmatrix},
\end{equation}
where $\mathbf{B}_1$ is an identity matrix, and the bandwidth of $\mathbf{B}_2$ depends on the support size of dual basis functions. For a constraint matrix $\mathbf{B}$ constructed using \Bezier dual basis functions, $\mathbf{B}_2$ is a sparse matrix with limited bandwidth, while the global dual basis functions leads to a dense $\mathbf{B}_2$.\par

For a $\mathbf{B}$ in the form~\eqref{eq:constraint-form} with $\mathbf{B}_1 = \mathbf{I}$, the vector basis of its null space can be obtained from
\begin{equation}
	\mathbf{C}=\begin{bmatrix}
		-\mathbf{B}_2 \\
		\mathbf{I}
	\end{bmatrix}.
	\label{eq:null-space}
\end{equation}
The mortar linear system can now be written as
\begin{equation}
	\mathbf{K}^{\text{mortar}}\mathbf{U}^{\text{mortar}}=\mathbf{C}^T\mathbf{K}\mathbf{C}\mathbf{U}^{\text{mortar}}=\mathbf{C}^T\mathbf{F}.\label{eq:mortar-form-discretized}
\end{equation}
The relation between the mortar displacement nodal value vector $\mathbf{U}^{\text{mortar}}$ and $\mathbf{U}^{\text{LM}}$ is given by
\begin{equation}
	\mathbf{U}^{\text{LM}}=\mathbf{C}\mathbf{U}^{\text{mortar}}.
\end{equation}
With a sparse $\mathbf{C}$ obtained from the \Bezier dual basis, the stiffness matrix of the mortar formulation $\mathbf{K}^{\text{mortar}}$ will remain sparse, resulting in an efficient linear system.

\section{A dual mortar formulation for the multi-patch biharmonic problem}
\label{sec:dual-mortar-formulation}

In this section, we present a formulation for the biharmonic problem over multi-patch tensor product domains. Because the biharmonic problem requires trial and test functions that are in $H^2$, we will use the dual mortar method to add constraints between patch boundaries to weakly enforce $C^1$ continuity. We begin by introducing concepts from domain decomposition.

\subsection{Domain decomposition}\label{sec:domain_decompostion}

\begin{figure}[ht]
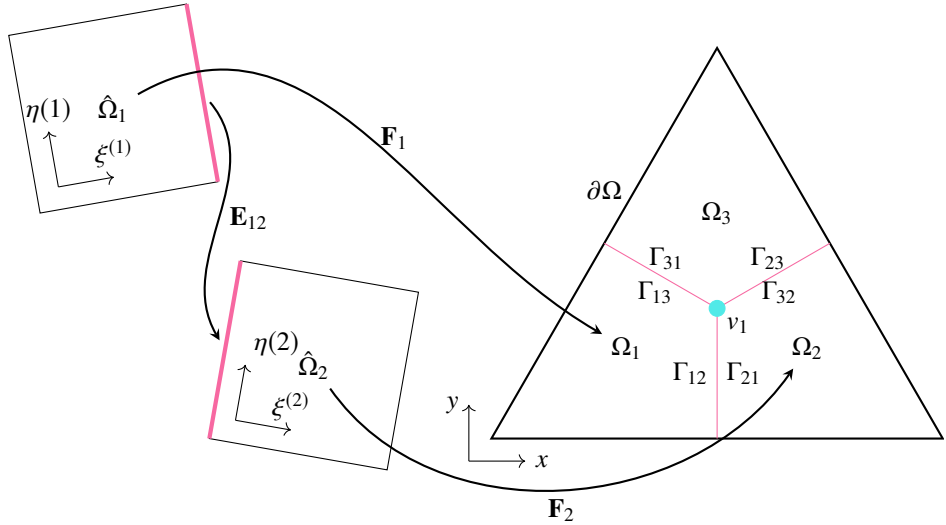

	\center
	\includestandalone[scale=1]{patch_partition}
	\caption{An example of a domain decomposition of a triangular domain. The patches are defined on different parametric domains and are connected via geometric mappings.}
	\label{fig:partition}
\end{figure}

Let $\Omega$ be a bounded open domain in $\mathbb{R}^2$ with its boundary denoted by $\partial\Omega$. We assume that $\Omega$ can be subdivided into $K$ non-overlapping patches $\Omega_k$ for $1\leq{}k\leq{}K$, i.e.
\begin{equation}
	\bar{\Omega} = \bigcup_{k=1}^{K} \bar{\Omega}_{k}  \quad \textrm{and} \quad {\Omega}_{k}\bigcap{\Omega}_{l}=\emptyset, \quad \forall k\neq{l}
\end{equation}
where $\bar{\Omega}_k$ is the closure of $\Omega_k$. For simplicity, we only consider the case where the intersection of two patches is either empty, a single vertex, or the entire edge. We denote the common interface of two neighboring patches as $\Gamma_{kl} = \partial{\Omega}_k\bigcap\partial{\Omega}_l$ so that $\Gamma_{kl}=\emptyset$ if $\Omega_k$ is not a neighbor of $\Omega_l$. We also define the skeleton $\mathbf{S}=\bigcup_{k,l\in{K}, k<{l}}\Gamma_{kl}$ as the union of all interfaces in $\Omega$. The set $\mathbf{V}$ denotes the set of all vertices in $\Omega$. A representative example of a multi-patch geometry is shown in Figure~\ref{fig:partition}.

For each patch, there exists a bijective geometric mapping from the parametric domain $\hat{\Omega}_k$ to the physical domain $\Omega_k$, which is defined as
\begin{equation}
	\mathbf{F}_{k}\left(\xi_{k},\eta_{k}\right)\colon{\hat{\Omega}_k}\mapsto{{\Omega}_k}\in\mathbb{R}^2,
\end{equation}
where $\left(\xi_{k},\eta_{k}\right)$ are the coordinates of the parametric domain. For simplicity and without loss of generality, we assume the parametric domain is ${\hat{\Omega}_k}=\left[{0,1}\right]\times\left[{0,1}\right]$ for all patches.

We can use the mappings $\mathbf{F}_k$ to create connections between neighboring patches. Due to the fact that $\mathbf{F}_k$ is a bijection, there exists an inverse mapping denoted by $\mathbf{F}^{-1}_k$. We can construct a bijective transformation on the intersection $\Gamma_{kl}$ as
\begin{equation}
	\mathbf{E}_{kl}=\mathbf{F}_{l}^{-1}\circ\mathbf{F}_{k},
\end{equation}
which maps a parametric point on $\partial{\hat{\Omega}}_k\bigcap \hat{\Gamma}_{kl}$ to a physical point on the intersection $\Gamma_{kl}$ and then to a parametric point on $\partial{\hat{\Omega}}_l\bigcap \hat{\Gamma}_{kl}$. With the mapping $\mathbf{E}_{kl}$ in hand, we are now ready to formulate the biharmonic problem over multi-patch domains.

\subsection{The biharmonic problem}

The strong form governing equation and boundary conditions of the homogeneous biharmonic problem defined over the domain $\Omega$ is given by the following:

\begin{equation}\label{eq:strong-form}
	\left\{\begin{alignedat}{2}
		\Delta^2u&=f\quad &&\text{on }\Omega,\\
		u=\frac{\partial{u}}{\partial{\mathbf{n}}}&=0\quad &&\text{on }\partial\Omega,
	\end{alignedat}\right.
\end{equation}
with $f\in{}L^2(\Omega)$ and $u\in{}C^4(\Omega)$. This problem can be restated in the weak form as: find $u\in{}H^2_0(\Omega)$ such that
\begin{equation}\label{eq:weak-form}
	a_b(u,v)=l(v),\quad\forall{}v\in{}H^2_0(\Omega),
\end{equation}
with
\begin{equation}
	\begin{split}
		a_b(u,v)&=\int_{\Omega}\nabla u \nabla v d \Omega,\\
		l(v) &= \int_{\Omega} fv d \Omega,
	\end{split}
\end{equation}
where $H^2_0(\Omega)$ is the Sobolev space containing all functions in the space $H^2(\Omega)$ that also satisfy the homogeneous Dirichlet boundary conditions in~\eqref{eq:strong-form}.

For a partitioned domain $\Omega$, as defined in Section~\ref{sec:domain_decompostion}, the construction of a finite dimensional subspace of $H^2_0(\Omega)$ is a nontrivial task because there is no guarantee that the discretization of neighboring domains is smooth enough across shared boundaries to satisfy the $H^2(\Omega)$ requirement. In order to handle multi-patch geometries, we will recast the biharmonic problem in terms of the following Lagrange multiplier formulation: find $u\in{\mathcal{X}_b}$, $\lambda_0\in{\mathcal{M}_0}$ and $\lambda_1\in{\mathcal{M}_1}$ such that:
\begin{equation}
	\left\{\begin{alignedat}{2}
		a_b(u,v)+b_0(\lambda_0,v)+b_1(\lambda_1,v)&=l(v)\quad&&\forall v\in{\mathcal{X}_b},\\
		b_0(\mu_0,u)&=0 \quad&&\forall \mu_0\in{\mathcal{M}_0},\\
		b_1(\mu_1,u)&=0 \quad&&\forall \mu_1\in{\mathcal{M}_1},
	\end{alignedat}\right.
\end{equation}
where
\begin{equation}
	\mathcal{X}_b\coloneq\{v\in{}L^2(\Omega)\,\vert\,{}v\vert_{\Omega_k}\in{}H^2(\Omega_k),\,{}1\leq{}k\leq{}K\text{  and }v\vert_{\partial{}\Omega}=\frac{\partial{}v}{\partial{\mathbf{n}}}\vert_{\partial{}\Omega}=0\},
\end{equation}
is an unconstrained broken Sobolev space endowed with the norm $\|u\|_{H^2_*(\Omega)}=\left(\sum_{k=1}^K \|u\|_{H^2(\Omega_k)}\right)^{1/2}$, $\mathcal{M}_0$ and $\mathcal{M}_1$ are Lagrange multiplier spaces, and $b_0$ and $b_1$ impose the required constraints on $u$ to satisfy the $H^2$ requirement. We will define $b_0$ and $b_1$ in the following section.

\begin{remark}
	Strictly speaking, the restriction of $u\in{}H^2(\Omega)$ to the boundary $\partial\Omega$ is ill-defined. To rigorously define the value of $u$ and its normal derivative on $\partial\Omega$, we need the help of the trace operator $T$. A standard trace theorem states~\cite{grisvard2011elliptic}: Given $\Omega$ with a boundary $\partial \Omega$ of class $C^{k,1}$ (i.e., $k$ times continuously differentiable and its $k^\text{th}$ order derivatives are Lipschitz continuous). Assume that $l\leq k$. Then the mapping
	\begin{equation}
		u \rightarrow \{Tu, T\frac{\partial u}{\partial \mathbf{n}}, \dots, T\frac{\partial^l u}{\partial^l \mathbf{n}}\},
	\end{equation}
	is a bounded linear operator from $H^{k+1}(\Omega)$ onto $\prod_{j=0}^lH^{k-j+\frac{1}{2}}(\partial \Omega)$. And if $u\in{}H^{k+1}(\Omega)\cap{}C^\infty(\Omega)$, we have the following relation
	\begin{equation}
		\{u, \frac{\partial u}{\partial \mathbf{n}}, \dots, \frac{\partial^l u}{\partial^l \mathbf{n}}\} =
		\{Tu, T\frac{\partial u}{\partial \mathbf{n}}, \dots, T\frac{\partial^l u}{\partial^l \mathbf{n}}\}
	\end{equation}
	Hence, for $u\in{}H^2(\Omega)$, we have $Tu\in H^{\frac{3}{2}}(\partial \Omega)$ and $T\frac{\partial u}{\partial \mathbf{n}}\in H^{\frac{1}{2}}(\partial \Omega)$. If $u\circ \mathbf{F} \in H^2({\hat{\Omega}})$ and each entries of $\nabla \mathbf{F}^{-1}$ are in $L^{\infty}(\Omega)$, then~\cite{ciarlet1972interpolation}
	\begin{equation}
		\begin{bmatrix}
			\frac{\partial u\circ \mathbf{F}}{\partial \xi} \circ \mathbf{F}^{-1} \\
			\frac{\partial u\circ \mathbf{F}}{\partial \eta} \circ \mathbf{F}^{-1}
		\end{bmatrix}
		\in H^{1}(\Omega)^2.
	\end{equation}
	As a result, $T\left(\frac{\partial u\circ \mathbf{F}}{\partial \xi} \circ \mathbf{F}^{-1}\right)\in H^{\frac{1}{2}}(\partial \Omega)$ and $T\left(\frac{\partial u\circ \mathbf{F}}{\partial \eta} \circ \mathbf{F}^{-1}\right)\in H^{\frac{1}{2}}(\partial \Omega)$.\par
	This result guarantees $H^2_0(\Omega)$ and the inter-patch constraints that will be discussed in the next section are all well-defined. In order to avoid cumbersome notation, the trace operator $T(\cdot)$ is suppressed and we will use $\{ \frac{\partial u}{\partial \xi}, \frac{\partial u}{\partial \eta} \}$ to refer to $\{ \frac{\partial u\circ \mathbf{F}}{\partial \xi} \circ \mathbf{F}^{-1},\frac{\partial u\circ \mathbf{F}}{\partial \eta} \circ \mathbf{F}^{-1} \}$.
\end{remark}

\subsection{Dual-compatible $C^1$ constraints}

We now develop a set of constraints to impose $C^1$ continuity across patch boundaries under the dual mortar framework. Since $C^1(\Omega) \subset H^2(\Omega)$, we will be able to use these constraints to solve the multi-patch biharmonic problem.\par

\begin{figure}[ht]
	\center
	\includegraphics[scale=.4]{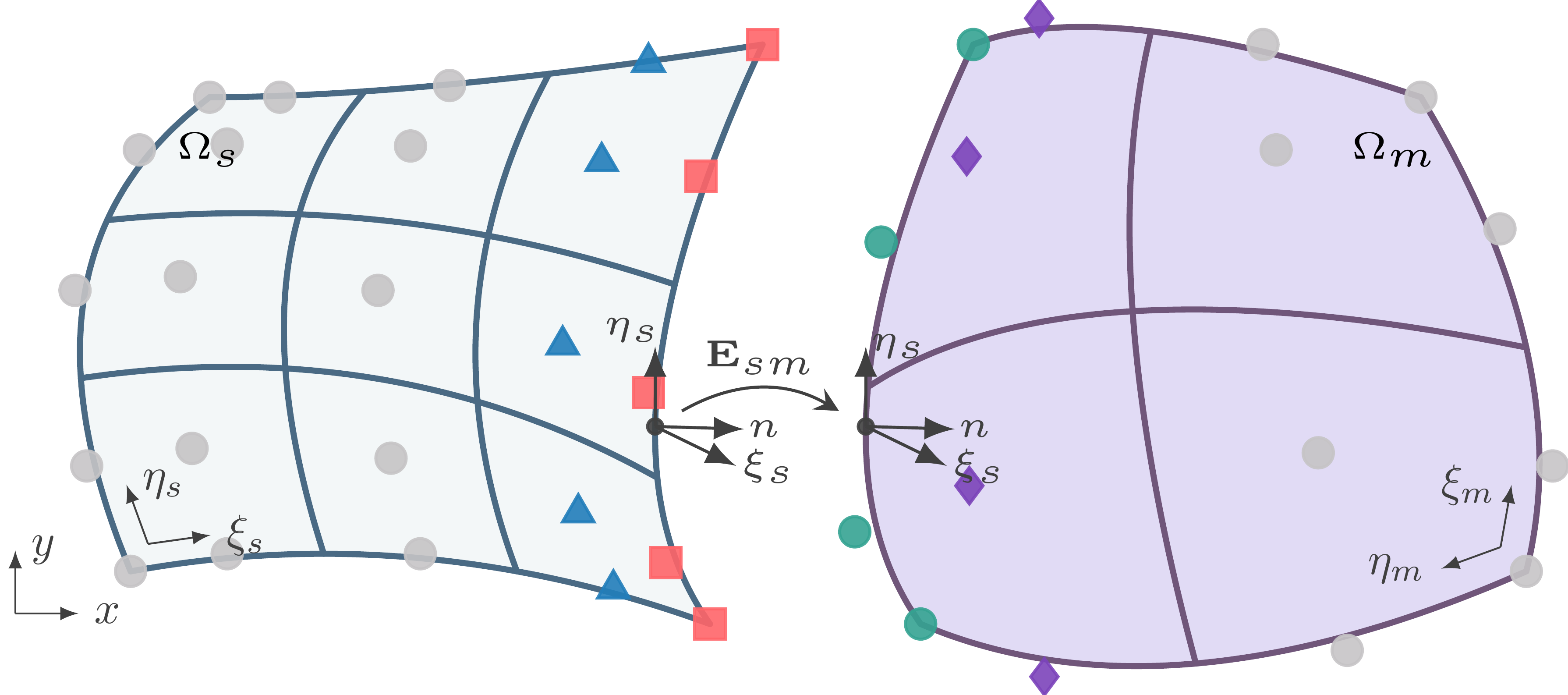}
	\caption{A two-patch planar domain $\Omega$ consisting of two patches $\Omega_s$ and $\Omega_m$ that are defined by the two mappings $\mathbf{F}_s$ and $\mathbf{F}_m$.}
	\label{fig:two-patch-mesh}
\end{figure}

To illustrate our method, we consider the construction of $C^1$ constraints for the two-patch domain shown in Figure~\ref{fig:two-patch-mesh}. We call $\Omega_s$ the slave domain and $\Omega_m$ the master domain. For a function $u \in C^1( \Omega_s \cup \Omega_m )$ with
\begin{equation}
	u=
	\begin{cases}
		u_s\quad \text{in }\Omega_s \\
		u_m\quad \text{in }\Omega_m,
	\end{cases}
\end{equation}
and $u_s\in{}C^1(\Omega_s)$, $u_m\in{}C^1(\Omega_m)$ the following two constraints are required across the intersection $\Gamma_{sm}$:
\begin{subequations}
	\begin{align}
		[u]_{\Gamma_{sm}}                                                   & =0,\label{eq:biharmonic_c0}                                                        \\
		\left[\frac{\partial{u}}{\partial{\mathbf{n}}}\right]_{\Gamma_{sm}} & =0,\quad\text{with }\mathbf{n}=\mathbf{n}_s=-\mathbf{n}_m,\label{eq:biharmonic_c1}
	\end{align}
\end{subequations}
where $\mathbf{n}_k$ is the outward normal direction of $\partial{\Omega_k}$ and
\begin{equation}
	[\cdot]_{\Gamma_{sm}}\coloneq\cdot\vert_{\Omega_s}-\cdot\vert_{\Omega_m}
\end{equation}
is the jump operator. The continuity constraint~\eqref{eq:biharmonic_c0} can naturally be incorporated into the framework of the dual mortar formulation. The smoothness constraint~\eqref{eq:biharmonic_c1}, however, can not be directly imposed. First, the existence of a dual basis for $\frac{\partial{N_i}}{\partial{\mathbf{n}}}\vert_{\Gamma_{sm}}$ is doubtful. Even if these dual basis functions do exist, since they are biorthogonal to the normal derivative of the basis functions, their formulation will depend on the parameterization of $\Gamma_{sm}$ and the geometric information of $\Omega_s$. This complex geometric dependence would destroy the simplicity of the dual basis formulation. To overcome this issue, we instead propose a smoothness constraint involving parametric derivatives only.

\begin{lemma}
	Given two differentiable bijective geometric mappings $\mathbf{F}_{s}\colon\hat{\Omega}_s\rightarrow{\Omega}_s$ and $\mathbf{F}_{m}\colon\hat{\Omega}_m\rightarrow{\Omega}_m$, a $C^0$-continuous function $u$ is $C^1$-continuous in the physical domain if and only if
	\begin{equation}
		\left[\frac{\partial{u}}{\partial{\mathbf{\xi}_s}}\right]_{\Gamma_{sm}}=0\quad\text{ and }\quad\left[\frac{\partial{u}}{\partial{\mathbf{\eta}_s}}\right]_{\Gamma_{sm}}=0. \label{eq:biharmonic_c1_mod}
	\end{equation}

	\begin{proof}
		It suffices to consider two neighboring patches as shown in Figure~\ref{fig:two-patch-mesh}. In this configuration, if $u$ is a $ C^0$-continuous function then $[\frac{\partial{u}}{\partial{\mathbf{\eta}_s}}]_{\Gamma_{sm}}=0$. If $u$ is also $C^1$-continuous, we have
		\begin{equation}
			0=\left[\frac{\partial{u}}{\partial{\mathbf{n}}}\right]_{\Gamma_{sm}}=\left[\frac{\partial{u}}{\partial{\mathbf{\xi}_s}}\right]_{\Gamma_{sm}}\frac{\partial{}\xi_s}{\partial{\mathbf{n}}}+\left[\frac{\partial{u}}{\partial{\mathbf{\eta}_s}}\right]_{\Gamma_{sm}}\frac{\partial{}\eta_s}{\partial{\mathbf{n}}}\Longrightarrow\left[\frac{\partial{u}}{\partial{\mathbf{\xi}_s}}\right]_{\Gamma_{sm}}\frac{\partial{}\xi_s}{\partial{\mathbf{n}}}=0
		\end{equation}
		The fact that $\mathbf{F}_{s}$ is bijective and $\frac{\partial{}\eta_s}{\partial{\mathbf{n}}} = 0$ indicates $\frac{\partial{}\xi_s}{\partial{\mathbf{n}}}\neq{}0$. Hence, $\left[\frac{\partial{u}}{\partial{\mathbf{\xi}_s}}\right]_{\Gamma_{sm}}=0$. On the other hand,
		\begin{equation}
			\begin{cases}
				\left[\frac{\partial{u}}{\partial{\mathbf{\xi}_s}}\right]_{\Gamma_{sm}}=0 \\
				\left[\frac{\partial{u}}{\partial{\mathbf{\eta}_s}}\right]_{\Gamma_{sm}}=0
			\end{cases}
			\Longrightarrow
			\left[\frac{\partial{u}}{\partial{\mathbf{n}}}\right]_{\Gamma_{sm}}=\left[\frac{\partial{u}}{\partial{\mathbf{\xi}_s}}\right]_{\Gamma_{sm}}\frac{\partial{}\xi_s}{\partial{\mathbf{n}}}+\left[\frac{\partial{u}}{\partial{\mathbf{\eta}_s}}\right]_{\Gamma_{sm}}\frac{\partial{}\eta_s}{\partial{\mathbf{n}}}=0
		\end{equation}
		This concludes the proof.
	\end{proof}
\end{lemma}

Hence, the constraints in~\eqref{eq:biharmonic_c1_mod} are equivalent to constraint~\eqref{eq:biharmonic_c1}. On an intersection that is parallel to the $\eta_s$ direction in the parametric domain, the constraint $\left[\frac{\partial{u}}{\partial{\mathbf{\xi}_s}}\right]_{\Gamma_{sm}}=0$ is utilized; on an intersection that is parallel to the $\xi_s$ direction in the parametric domain, the constraint $\left[\frac{\partial{u}}{\partial{\mathbf{\eta}_s}}\right]_{\Gamma_{sm}}=0$ is utilized.

\begin{remark}
	In order to demonstrate the advantages of the constraints in~\eqref{eq:biharmonic_c1_mod}, we consider the following intergral:
	\begin{equation}
		\begin{split}
			\int_{\Gamma_{sm}}\frac{\partial{}N_a(\xi_s,\eta_s)}{\partial{}\xi_s}\hat{N}_{j}(\eta_s)d\Gamma&=\int_{\Gamma_{sm}}\frac{\partial{}N_{n_{\xi_s}-1}(1)N_{i}(\eta_s)}{\partial{}\xi_s}\hat{N}_{j}(\eta_s)d\Gamma\\
			&=\frac{\partial{}N_{n_{\xi_s}-1}(1)}{\partial{}\xi_s}\int_{\Gamma_{sm}}N_{i}(\eta_s)\hat{N}_{j}(\eta_s)d\Gamma
		\end{split}\quad i,j\in\left\{1,2,\dots, n_{\eta_s}\right\}
		\label{eq:advan_constraint}
	\end{equation}
	where $n_{\xi_s}$ and $n_{\eta_s}$ are the number of nodes in the $\xi_s$ and $\eta_s$ directions of the slave patch, respectively, and the index $a = n_{\xi_s}i-1$. This integral is one term that is involved in the discretization of the constraint~\eqref{eq:biharmonic_c1} and is constructed by a Lagrange multiplier basis function $\hat{N}_j(\eta_s)$ and an activated basis function of the slave patch $N_a(\xi_s,\eta_s)$ that is one column away from the intersection (denoted by the blue triangles in Figure~\ref{fig:two-patch-mesh}). Due to the tensor product structure of multivariate spline basis functions, the derivative in one direction ($\xi_s$ for this case) will not influence the contributions coming from other directions. Hence, the dual basis function of an activated basis function in the constraint $\left[\frac{\partial{u}}{\partial{\mathbf{\xi}_s}}\right]_{\Gamma_{sm}}=0$ can be constructed by the dual basis function of its $\eta_s$ component divided by $\frac{\partial{}N_{n_{\xi_s}-1}(1)}{\partial{}\xi_s}$.
\end{remark}
The only issue now is how to evaluate the derivative of $u_m$ w.r.t. $\xi_s$ or $\eta_s$ directions. This can be done by considering the following chain rule
\begin{align}
	\begin{bmatrix}
		\tfrac{\partial{u_m}}{\partial{\xi_s}} \\
		\tfrac{\partial{u_m}}{\partial{\eta_s}}
	\end{bmatrix}
	=
	\begin{bmatrix}
		\frac{\partial\xi_m}{\partial\xi_s}  & \frac{\partial\xi_m}{\partial\eta_s}  \\
		\frac{\partial\eta_m}{\partial\xi_s} & \frac{\partial\eta_m}{\partial\eta_s}
	\end{bmatrix}^T
	\cdot
	\begin{bmatrix}
		\tfrac{\partial{u_m}}{\partial{\xi_m}} \\
		\tfrac{\partial{u_m}}{\partial{\eta_m}}
	\end{bmatrix}
	=
	\nabla\mathbf{E}_{sm}^T
	\cdot
	\begin{bmatrix}
		\tfrac{\partial{u_m}}{\partial{\xi_m}} \\
		\tfrac{\partial{u_m}}{\partial{\eta_m}}
	\end{bmatrix}.
\end{align}
The Jacobian of the composition mapping $\mathbf{E}_{sm}$ can be written as
\begin{equation}
	\nabla\mathbf{E}_{sm}=\nabla(\mathbf{F}_{m}^{-1}\circ\mathbf{F}_{s})=\nabla(\mathbf{F}_{m}^{-1})\cdot{}\nabla\mathbf{F}_{s}=(\nabla\mathbf{F}_{m})^{-1}\cdot{}\nabla\mathbf{F}_{s}.
\end{equation}\par

\subsection{The dual mortar formulation}
The Lagrange multiplier formulation for the multi-patch biharmonic problem can be defined as: find $u\in{\mathcal{X}_b}$, $\lambda_0\in{\mathcal{M}_0}$ and $\lambda_1\in{\mathcal{M}_1}$ such that
\begin{equation}
	\left\{\begin{alignedat}{2}
		a_b(u,v)+b_0(\lambda_0,v)+b_1(\lambda_1,v)&=l(v)\quad&&\forall v\in{\mathcal{X}_b},\\
		b_0(\mu_0,u)&=0 \quad&&\forall \mu_0\in{\mathcal{M}_0},\\
		b_1(\mu_1,u)&=0 \quad&&\forall \mu_1\in{\mathcal{M}_1},
	\end{alignedat}\right.\label{eq:biharmonic_mixed}
\end{equation}
with
\begin{subequations}
	\begin{align}
		b_0(\lambda_0,v) & =\sum_{\Gamma\in\mathbf{S}}\int_\Gamma[u]_{\Gamma}\lambda_0d\Gamma,\label{eq:operator-b0}                                                                                                                                                                                                                          \\
		b_1(\lambda_1,v) & =\sum_{\Gamma\in\mathbf{S}}\left(\int_\Gamma\left[\frac{\partial{u}}{\partial{\xi_s}}\right]_{\Gamma}\lambda_1d\Gamma\text{ if }\Gamma\parallel\eta_s\text{ or }\int_\Gamma\left[\frac{\partial{u}}{\partial{\eta_s}}\right]_{\Gamma}\lambda_1d\Gamma\text{ if }\Gamma\parallel\xi_s\right).\label{eq:operator-b1}
	\end{align}
\end{subequations}
The constrained function space required by the dual mortar formulation of the multi-patch biharmonic problem can then be defined as
\begin{equation}
	\mathcal{K}_b:=\left\{u\in{}\mathcal{X}_b\,\vert\,b_0(\mu_0,u)=0 \text{ and }b_1(\mu_1,u)=0\,\forall(\mu_0,\mu_1)\in{\mathcal{M}_0\times{}\mathcal{M}_1}\right\}.\label{eq:constrained_space}
\end{equation}

\subsection{Discretization}

For each intersection, the two adjacent patches are classified as either slave $\Omega_s$ or master $\Omega_m$. One patch can, at the same time, be a master for one intersection and a slave for another intersection. To approximate the solution of the variational problem, we use B-spline basis functions $\{N^s_i\}_{i\in{I_s}}$ and $\{N^m_i\}_{i\in{I_m}}$ to discretize $\Omega_s$ and $\Omega_m$, respectively. Note that, on the intersection, we select the side with finer trace mesh as the slave patch and denote the other side as the master patch. This selection strategy can minimize the error from variational crimes~\cite{strang1973analysis,brenner_mathematical_2007}.
An appropriate indexing is chosen so that there is no overlap between the index sets $I_s$ and $I_m$ (i.e., given $n_s$ basis functions in $\Omega_s$, we can assume the starting index in the index set $I_m$ is $n_s+1$). The discretized geometrical mappings are represented by
\begin{align}
	\mathbf{F}_s & =\sum_{i\in{I_s}}\mathbf{P}_i^sN_i^s, \\
	\mathbf{F}_m & =\sum_{i\in{I_m}}\mathbf{P}_i^mN_i^m,
\end{align}
where the control points $\mathbf{P}_i^s,\mathbf{P}_i^m\in\mathbb{R}^2$. The discrete space $\mathcal{X}_b^h \subset \mathcal{X}_b$ contains the discretized test and weighting functions. In other words,
\begin{equation}
	u^h=\sum_{i\in{I_s\cup{}I_m}}U_iN_i,\quad v^h=\sum_{i\in{I_s\cup{}I_m}}V_iN_i
\end{equation}
with
\begin{align}
	N_i=
	\begin{cases}
		N_i^s \quad & i\in{I_s}, \\
		N_i^m \quad & i\in{I_m}.
	\end{cases}
\end{align}
The discrete Lagrange multiplier spaces $\mathcal{M}_0^h \subset \mathcal{M}_0$ and $\mathcal{M}_1^h \subset \mathcal{M}_1^h$ are created using the dual basis. Depending on the orientation of the intersection, we have that
\begin{itemize}
	\item for the intersection $\xi_s=0$,
	      \begin{equation}
		      \begin{split}
			      \lambda_0^h&=\sum_{i=1}^{n_{\eta_s}}\Lambda^0_i\hat{N}^s_i(\eta_s),\quad  \mu_0^h=\sum_{i=1}^{n_{\eta_s}}\delta\Lambda^0_i\hat{N}^s_i(\eta_s)\\
			      \lambda_1^h&=\sum_{i=1}^{n_{\eta_s}}\Lambda^1_i\frac{\hat{N}^s_i(\eta_s)}{c}, \quad \mu_1^h=\sum_{i=1}^{n_{\eta_s}}\delta\Lambda^1_i\frac{\hat{N}^s_i(\eta_s)}{c},\quad c=\left.\frac{\partial{N^s_2(\xi_s)}}{\partial{\xi_s}}\right|_{\xi_s=0},
		      \end{split}
	      \end{equation}
	\item for the intersection $\xi_s=1$,
	      \begin{equation}
		      \begin{split}
			      \lambda_0^h&=\sum_{i=1}^{n_{\eta_s}}\Lambda^0_i\hat{N}^s_i(\eta_s),\quad \mu_0^h=\sum_{i=1}^{n_{\eta_s}}\delta\Lambda^0_i\hat{N}^s_i(\eta_s)\\
			      \lambda_1^h&=\sum_{i=1}^{n_{\eta_s}}\Lambda^1_i\frac{\hat{N}^s_i(\eta_s)}{c},\quad \mu_1^h=\sum_{i=1}^{n_{\eta_s}}\delta\Lambda^1_i\frac{\hat{N}^s_i(\eta_s)}{c},\quad c=\left.\frac{\partial{N^s_{n_{\xi_s}-1}(\xi_s)}}{\partial{\xi_s}}\right|_{\xi_s=1},
		      \end{split}
	      \end{equation}
	\item for the intersection $\eta_s=0$,
	      \begin{equation}
		      \begin{split}
			      \lambda_0^h&=\sum_{i=1}^{n_{\xi_s}}\Lambda^0_i\hat{N}^s_i(\xi_s),\quad \mu_0^h=\sum_{i=1}^{n_{\xi_s}}\delta\Lambda^0_i\hat{N}^s_i(\xi_s)\\
			      \lambda_1^h&=\sum_{i=1}^{n_{\xi_s}}\Lambda^1_i\frac{\hat{N}^s_i(\xi_s)}{c},\quad \mu_1^h=\sum_{i=1}^{n_{\xi_s}}\delta\Lambda^1_i\frac{\hat{N}^s_i(\xi_s)}{c},\quad c=\left.\frac{\partial{N^s_2(\eta_s)}}{\partial{\eta_s}}\right|_{\eta_s=0},
		      \end{split}
	      \end{equation}
	\item for the intersection $\eta_s=1$,
	      \begin{equation}
		      \begin{split}
			      \lambda_0^h&=\sum_{i=1}^{n_{\xi_s}}\Lambda^0_i\hat{N}^s_i(\xi_s), \quad  \mu_0^h=\sum_{i=1}^{n_{\xi_s}}\delta\Lambda^0_i\hat{N}^s_i(\xi_s)\\
			      \lambda_1^h&=\sum_{i=1}^{n_{\xi_s}}\Lambda^1_i\frac{\hat{N}^s_i(\xi_s)}{c},\quad \mu_1^h=\sum_{i=1}^{n_{\xi_s}}\delta\Lambda^1_i\frac{\hat{N}^s_i(\xi_s)}{c},\quad c=\left.\frac{\partial{N^s_{n_{\eta_s}-1}(\eta_s)}}{\partial{\eta_s}}\right|_{\eta_s=1}.
		      \end{split}
	      \end{equation}
\end{itemize}\par

By substituting the discretized displacement field and Lagrange multipliers into the bilinear form $a_b(\cdot,\cdot)$, $b_0(\cdot,\cdot)$ and $b_1(\cdot,\cdot)$, we obtain the following stiffness and constraint matrices
\begin{equation}
	\mathbf{V}^T\mathbf{K}_b\mathbf{U}=a_b(u^h,v^h)\quad\text{and}\quad
	\begin{bmatrix}
		\delta\mathbf{\Lambda}^0 \\
		\delta\mathbf{\Lambda}^1
	\end{bmatrix}^T\mathbf{B}_b\mathbf{U}=\begin{bmatrix}
		b_0(\mu_0^h,u^h) \\
		b_1(\mu_1^h,u^h)
	\end{bmatrix}.
\end{equation}
The remaining question is how to effectively construct the vector basis (in matrix form $\mathbf{C}_b$) of the null space of $\mathbf{B}_b$ such that the resulting basis functions (constructed by Equation~\eqref{eq:basis-null-space}) of $\mathcal{K}_b^h \subset \mathcal{K}_b$ have compact support and lead to a sparse stiffness matrix in~\eqref{eq:mortar-form-discretized}.

\section{Building a basis for the null space of $\mathbf{B}_b$}

\subsection{The two-patch case}

Recall from Section~\ref{sec:dual-mortar-method} that for a constraint matrix taking the form~\eqref{eq:constraint-form} the corresponding operator $\mathbf{C}$ can be constructed in an elegant manner via Equation~\eqref{eq:null-space}. In this section, we will show how to recover form~\eqref{eq:constraint-form} from the constraint matrix $\mathbf{B}_b$ via a simple linear transformation. We first classify the basis functions of $\mathcal{X}_b^h$ into five different types, depending on their proximity to an interface, as shown in Figure~\ref{fig:two-patch-mesh}:
\begin{enumerate}
	\item The basis functions $N^s_i$ such that $\supp(N^s_i)\bigcap\Gamma_{sm}=\emptyset$ and $\supp(\frac{\partial{}N^s_i}{\partial{\xi_s}})\bigcap\Gamma_{sm}\neq\emptyset$ (the second closest column of slave basis functions to the intersection $\Gamma_{sm}$), whose indices are denoted by the index set $I_\text{i}$. (denoted by blue triangles)
	\item The basis functions $N^s_i$ such that $\supp(N^s_i)\bigcap\Gamma_{sm}\neq\emptyset$ (the column of slave basis functions on the intersection $\Gamma_{sm}$), whose indices are denoted by the index set $I_\text{ii}$. (denoted by red squares)
	\item The basis functions $N^m_i$ such that $\supp(N^m_i)\bigcap\Gamma_{sm}\neq\emptyset$  (the column of master basis functions on the intersection $\Gamma_{sm}$), whose indices are denoted by the index set $I_\text{iii}$. (denoted by green circles)
	\item The basis functions $N^m_i$ such that $\supp(N^m_i)\bigcap\Gamma_{sm}=\emptyset$ and $\supp(\frac{\partial{}N^m_i}{\partial{\xi_s}})\bigcap\Gamma_{sm}\neq\emptyset$  (the second closest column of master basis functions to the intersection $\Gamma_{sm}$), whose indices are in the index set $I_\text{iv}$. (denoted by purple diamonds)
	\item The basis functions $N^m_i$ whose values and first order derivative values in the $\xi_s$ direction are zero on $\Gamma_{sm}$, whose indices are denoted by the index set $I_\text{v}$. (denoted by grey circles)
\end{enumerate}\par

\begin{figure}[ht]
	\center
	\begin{subfigure}[b]{0.45\textwidth}
		\includegraphics[scale=.15]{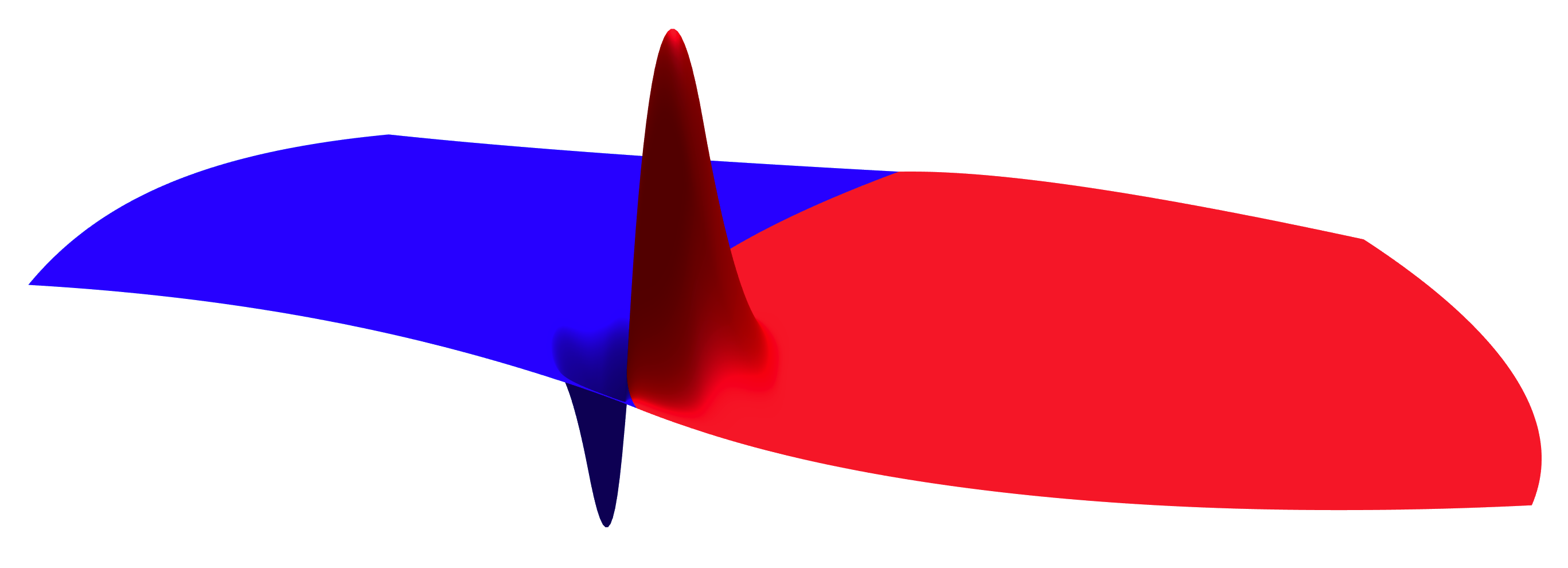}
	\end{subfigure}
	\hfill
	\begin{subfigure}[b]{0.45\textwidth}
		\includegraphics[scale=.15]{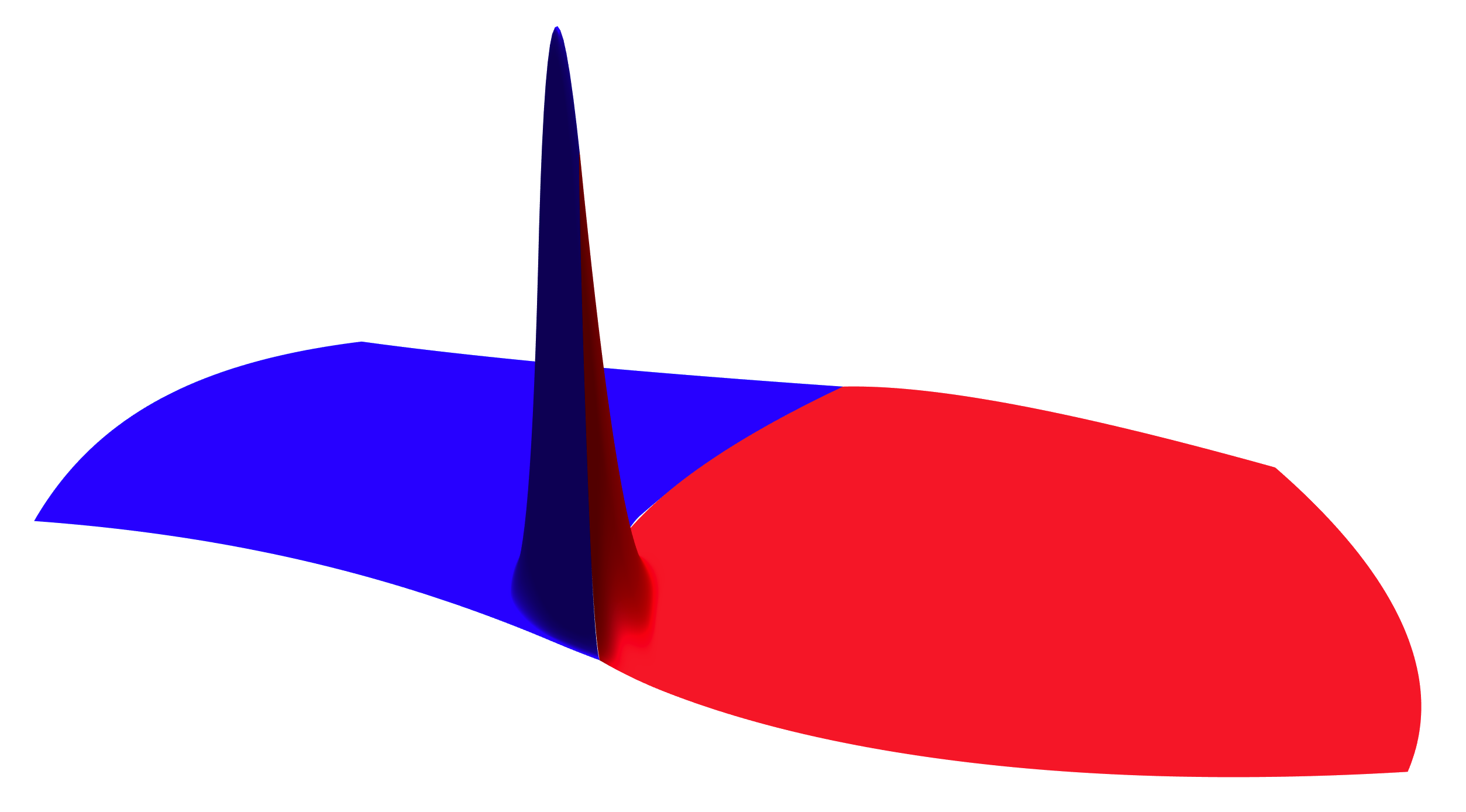}
	\end{subfigure}
	\caption{Two exemplary basis functions in the constrained space $\mathcal{K}_b^h$.}
	\label{fig:constrained_basis}
\end{figure}

Since the structure of the constraint matrix $\mathbf{B}_b$ depends on the index sets $I_s$ and $I_m$ and the ordering of the Lagrange multiplier basis functions we introduce two permutation matrices $\mathbf{P}_c$ and $\mathbf{P}_r$ (this step is not neccessary from the implementation point-of-view, but is helpful during the derivation, especially for multi-patch problems). We define the column-wise permutation matrix $\mathbf{P}_c$ as
\begin{equation}
	\begin{bmatrix}
		\mathbf{I}_\text{i}   \\
		\mathbf{I}_\text{ii}  \\
		\mathbf{I}_\text{iii} \\
		\mathbf{I}_\text{iv}  \\
		\mathbf{I}_\text{v}
	\end{bmatrix}=
	\mathbf{P}_c
	\begin{bmatrix}
		\mathbf{I}_s \\
		\mathbf{I}_m
	\end{bmatrix},
\end{equation}
where $\mathbf{I}_i$ is the vector form of the index set  $I_i$. We also define a row-wise permutation matrix $\mathbf{P}_r$ such that the permuted constraint matrix can be written in the partitioned form
\begin{equation}
	\mathbf{B}_p=\mathbf{P}_r\mathbf{B}_b\mathbf{P}_c^T=
	\begin{bmatrix}
		\mathbf{B}_1^1 & \mathbf{B}_1^2 & \mathbf{B}_1^3 & \mathbf{B}_1^4 & \mathbf{0} \\
		\mathbf{0}     & \mathbf{B}_2^2 & \mathbf{B}_2^3 & \mathbf{0}     & \mathbf{0}
	\end{bmatrix},
\end{equation}
where $\mathbf{B}_1^1$ is the contribution of the first type of B-spline basis function in the discretization of $b_1$ and $\mathbf{B}_2^2$ is the contribution of the second type of B-spline basis function in the discretization of $b_0$. Under the row-wise permutation matrix $\mathbf{P}_r$, $\mathbf{B}_1^1$ and $\mathbf{B}_2^2$ become identity submatrices. Under a rank-preserving transformation $\mathbf{T}$ we can eliminate the submatrix $\mathbf{B}_1^2$ such that
\begin{equation}
	\sbox0{$\begin{matrix}\mathbf{B}_1^3-\mathbf{B}_1^2\mathbf{B}_2^3 & \mathbf{B}_1^4 & \mathbf{0} \\ \mathbf{B}_2^3 & \mathbf{0} & \mathbf{0}\end{matrix}$}
	\mathbf{T}\mathbf{B}_p=\left[
		\begin{array}{c:c}
			\makebox[\wd0/3]{\large $\mathbf{I}$} & \usebox{0} \\
		\end{array}
		\right].\label{eq:simple_form}
\end{equation}
We may now take
\begin{equation}
	\mathbf{C}_p=
	\left[\begin{array}{ccc}
			\mathbf{B}_1^2\mathbf{B}_2^3-\mathbf{B}_1^3 & -\mathbf{B}_1^4 & \mathbf{0}        \\
			-\mathbf{B}_2^3                             & \mathbf{0}      & \mathbf{0}        \\ \hdashline[2pt/2pt]
			\multicolumn{3}{c}{\multirow{3}{*}{\raisebox{0mm}{\scalebox{1.5}{$\mathbf{I}$}}}} \\
			                                            &                 &
		\end{array}\right].
\end{equation}
The vector basis of the null space of $\mathbf{B}_b$ can now be obtained from
\begin{equation}
	\mathbf{C}_b=\mathbf{P}_c^T\mathbf{C}_p.\label{eq:permuted_back_nullspace}
\end{equation}
Examples of basis functions, represented by vectors of $\mathbf{C}_b$, are shown in Figure~\ref{fig:constrained_basis}.

\subsection{The multi-patch case}

\begin{figure}[ht]
	\center
	\includegraphics[scale=.4]{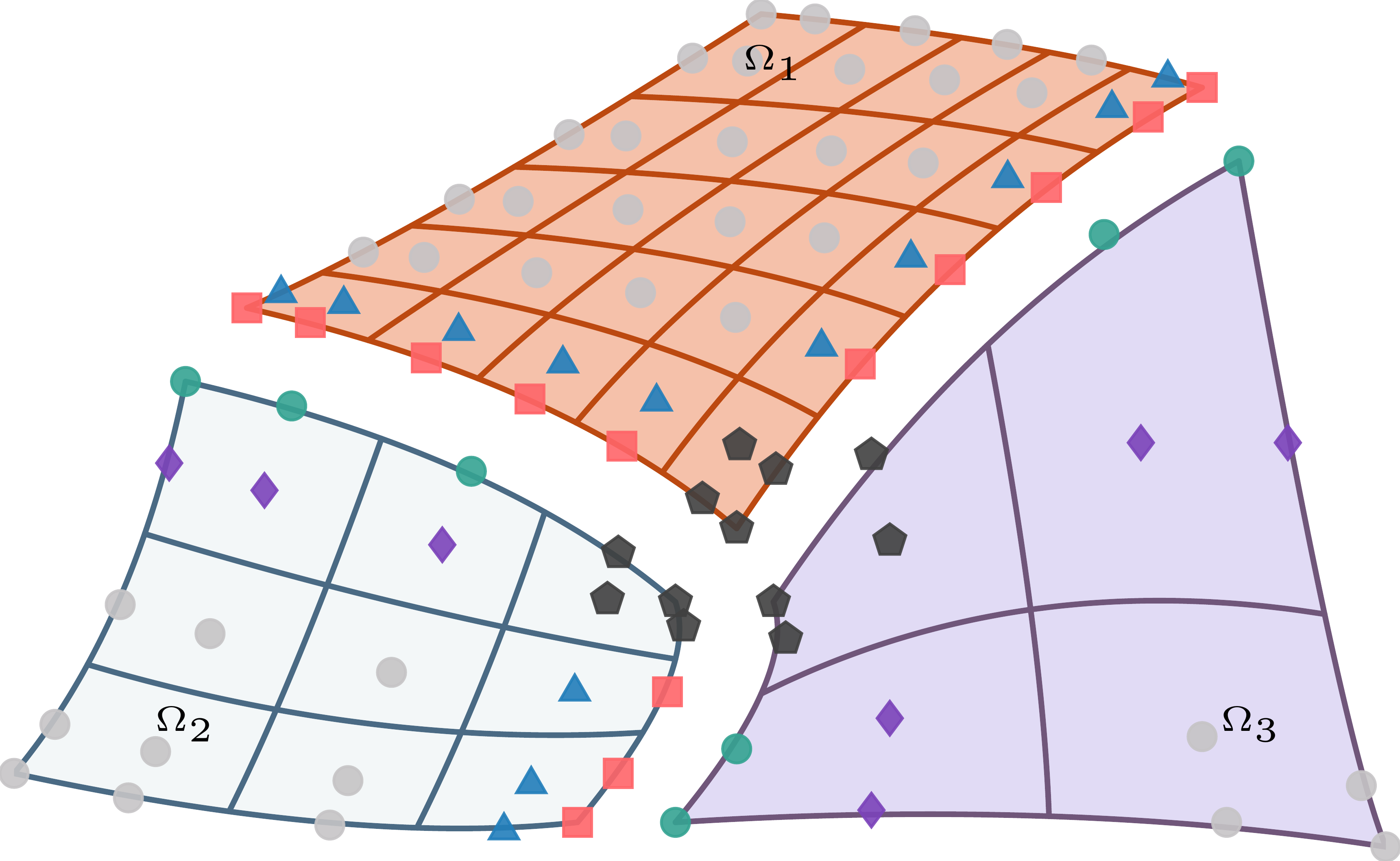}
	\caption{A three-patch planar domain $\Omega$ consisting of $\Omega_1$, $\Omega_2$ and $\Omega_3$.}
	\label{fig:three-patch-mesh}
\end{figure}

To extend this approach to more complex geometries requires the ability to stitch together multiple patches as shown in Figure~\ref{fig:three-patch-mesh}. To build the null space of a constraint matrix in the neighborhood of a vertex requires special care. From an implementation perspective, if we naively let the Lagrange multiplier spaces along all the interfaces adjacent to the vertex have the same dimension as the univariate spline basis along the slave side of the interface, there will be basis functions which serve as both slave and master. As an example, consider the the basis functions corresponding to the black pentagons of patch $\Omega_2$ in Figure~\ref{fig:three-patch-mesh}. As a result, there is no permutation under which the constraint matrix $\mathbf{B}$ can be modified to form~\eqref{eq:constraint-form} so that the basis of the null space can be found in a trivial way. Additionally, although the constraint matrices defined on each interface adjacent to a vertex are full row rank, the assembled constraint matrix $\mathbf{B}$ may not be full row rank. To overcome this problem, one may adopt matrix factorization techniques to solve the null space problem of a general constraint matrix $\mathbf{B}$, including LU, QR, SVD, etc. For example, a rank-revealing QR factorization of a rank-deficient constraint matrix $\mathbf{B}$ yields
\begin{equation}
	\mathbf{B}\mathbf{P}=
	\mathbf{Q}
	\begin{bmatrix}
		\mathbf{R}_1 & \mathbf{R}_2 \\
		\mathbf{0}   & \mathbf{0}
	\end{bmatrix}
\end{equation}
where $\mathbf{P}$ is a permutation matrix, $\mathbf{Q}$ is a unitary matrix, $\mathbf{R}_1$ is an upper triangular matrix and $\mathbf{R}_2$ is a rectangular matrix. The vector basis of the null space can then be taken to be
\begin{equation}
	\mathbf{C}=
	\mathbf{P}
	\begin{bmatrix}
		-\mathbf{R}_1^{-1}\mathbf{R}_2 \\
		\mathbf{I}
	\end{bmatrix}.\label{eq:direct_kernel}
\end{equation}
This type of global factorization has been utilized for patch coupling problems in~\cite{coox_robust_2017, coox2017flexible, dornisch2017dual}. However, it requires a global factorization of the entire constraint matrix $\mathbf{B}$, and fails to leverage the local properties of the dual basis. Moreover, the sparsity of the resulting constrained stiffness matrix might be negatively impacted since the inverse of $\mathbf{R}_1$ is a dense matrix. Additionally, \textit{inf-sup} stability may be violated and pathologies can be activated such as spurious oscillations and locking~\cite{barbosa1991finite}. For $2^\text{nd}$ order problems, one approach is to reduce the polynomial order of elements adjacent to vertices by one~\cite{bernardi_domain_1993,bernardi_basics_2005,belgacem_mortar_1998,brivadis_isogeometric_2015}. Then the modified Lagrange multiplier discretization is a subspace of the trace space of the slave patch of codimension $2$. By reducing the number of constraints, the basis functions in the neighborhood of vertices can now be considered masters. In addition, the modified Lagrange multiplier discretization is \textit{inf-sup} stable.\par



We introduce a sixth kind of B-spline basis function near a vertex $v\in \mathbf{V}$,
\begin{enumerate}
	\setcounter{enumi}{5}
	\item The basis function $N_i$ such that $\supp(N_i)\bigcap{v}\neq{}0$, or $\supp(\frac{\partial{}N_i}{\partial\xi})\bigcap{v}\neq{}0$,  or $\supp(\frac{\partial{}N_i}{\partial\eta})\bigcap{v}\neq{}0$ or $\supp(\frac{\partial^2{}N_i}{\partial\xi\partial\eta})\bigcap{v}\neq{}0$, whose indices are denoted by the index set $I_\text{vi}$ (denoted by black pentagons in Figure~\ref{fig:three-patch-mesh}).
\end{enumerate}
The definitions of the other five kinds of B-spline basis functions remain the same except that their intersection with the sixth kind are excluded, that is
\begin{equation}
	I_k=I_k-I_\text{vi}\bigcap{}I_k, \quad k\in\{\text{i},\text{ii},\dots,\text{v}\}.
\end{equation}

\subsubsection{Global dual basis multi-patch treatment}\label{sec:vertex_modification}
\begin{figure}[ht]
	\centering
	\includegraphics[width=.4\linewidth]{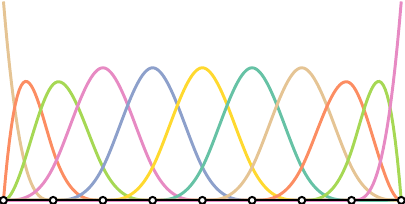}\\
	\includegraphics[width=.4\linewidth]{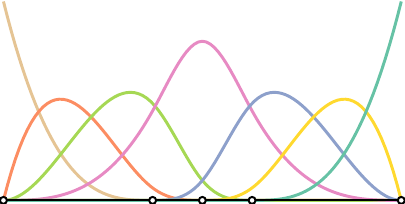}
	\caption{A coarsening procedure for cubic B-spline basis functions. Top: original cubic basis functions. Bottom: basis functions of the coasened mesh, where two knots on each end are removed.}\label{fig:boundary_modification}
\end{figure}

In this approach, we reduce the number of constraints on each interface such that all of the sixth kind of B-spline basis functions can be classified as masters. For the biharmonic problem this requires a reduction of four constraints per vertex per patch. We accomplish this by coarsening the mesh in the neighborhood of each vertex. Specifically, we remove the two knots adjacent to each vertex. An example of this coarsening procedure for cubic univariate B-spline basis functions is shown in Figure~\ref{fig:boundary_modification}. The corresponding global dual basis can then be constructed by using~\eqref{eq:global_dual}. For a set of B-spline basis functions $\left\{N_i\right\}_{i=1}^n$ and the corresponding coarsened global dual basis $\left\{\hat{N}^G_i\right\}_{i=1}^{n-4}$, the biorthogonality relation is then given as
\begin{equation}
	\int_\Gamma \hat{N}^G_iN_{j+2} d\Gamma = \delta_{ij}, \quad \forall 1\leq i,j-2\leq n-4.\label{eq:modify_biorthogonality}
\end{equation}
In other words, the biorthogonality relation holds for all but the two basis functions nearest the vertices.\par

A column-wise permutation matrix $\mathbf{P}_c$ which handles the multi-patch case can be defined as
\begin{equation}
	\begin{bmatrix}
		\mathbf{I}_\text{i}   \\
		\mathbf{I}_\text{ii}  \\
		\mathbf{I}_\text{iii} \\
		\mathbf{I}_\text{iv}  \\
		\mathbf{I}_\text{v}   \\
		\mathbf{I}_\text{vi}
	\end{bmatrix}=
	\mathbf{P}_c
	\begin{bmatrix}
		\mathbf{I}_1 \\
		\mathbf{I}_2 \\
		\vdots
	\end{bmatrix}.
\end{equation}
With the help of a row-wise permutation matrix, $\mathbf{P}_r$, the permuted constraint matrix can be written in the partitioned form
\begin{equation}
	\mathbf{B}^\text{mod}_p:=\mathbf{P}_r\mathbf{B}^\text{mod}_b\mathbf{P}_c^T=
	\begin{bmatrix}
		{\mathbf{B}}_1^1 & {\mathbf{B}}_1^2 & {\mathbf{B}}_1^3 & {\mathbf{B}}_1^4 & \mathbf{0} & {\mathbf{B}}_1^6 \\
		\mathbf{0}       & {\mathbf{B}}_2^2 & {\mathbf{B}}_2^3 & \mathbf{0}       & \mathbf{0} & {\mathbf{B}}_2^6
	\end{bmatrix},\label{eq:cross_point_constraint_structure}
\end{equation}
where $\mathbf{B}^\text{mod}_b$ is the constraint matrix constructed from the coarsened dual basis, $\mathbf{B}_1^1$ is the contribution of the first type of B-spline basis function in the discretization of $b_1$, and $\mathbf{B}_2^2$ is the contribution of the second type of B-spline basis function in the discretization of $b_0$. Under the row-wise permutation matrix $\mathbf{P}_r$, $\mathbf{B}_1^1$ and $\mathbf{B}_2^2$ become identity submatrices. Note that in the multi-patch case, basis functions in the neighborhood of in-domain vertices are excluded from the definitions of the first four types of B-spline basis functions. As a result, there is no basis function that will serve as both slave and master. As for the two patch case, the vector basis of the null space of $\mathbf{B}^\text{mod}_b$, denoted by $\mathbf{C}^\text{mod}_b$, can now be obtained from Equation~\eqref{eq:permuted_back_nullspace}, with
\begin{equation}
	\mathbf{C}^\text{mod}_p =
	\left[\begin{array}{cccc}
			{\mathbf{B}}_1^2{\mathbf{B}}_2^3-{\mathbf{B}}_1^3 & -{\mathbf{B}}_1^4 & \mathbf{0} & -{\mathbf{B}}_1^6 \\
			-{\mathbf{B}}_2^3                                 & \mathbf{0}        & \mathbf{0} & -{\mathbf{B}}_2^6 \\ \hdashline[2pt/2pt]
			\multicolumn{4}{c}{\multirow{3}{*}{\raisebox{0mm}{\scalebox{1.5}{$\mathbf{I}$}}}}                      \\
			                                                  &                   &            &
		\end{array}\right].\label{eq:nullspace_modification}
\end{equation}

\begin{figure}[ht]
	\centering
	\includegraphics[width=.7\linewidth]{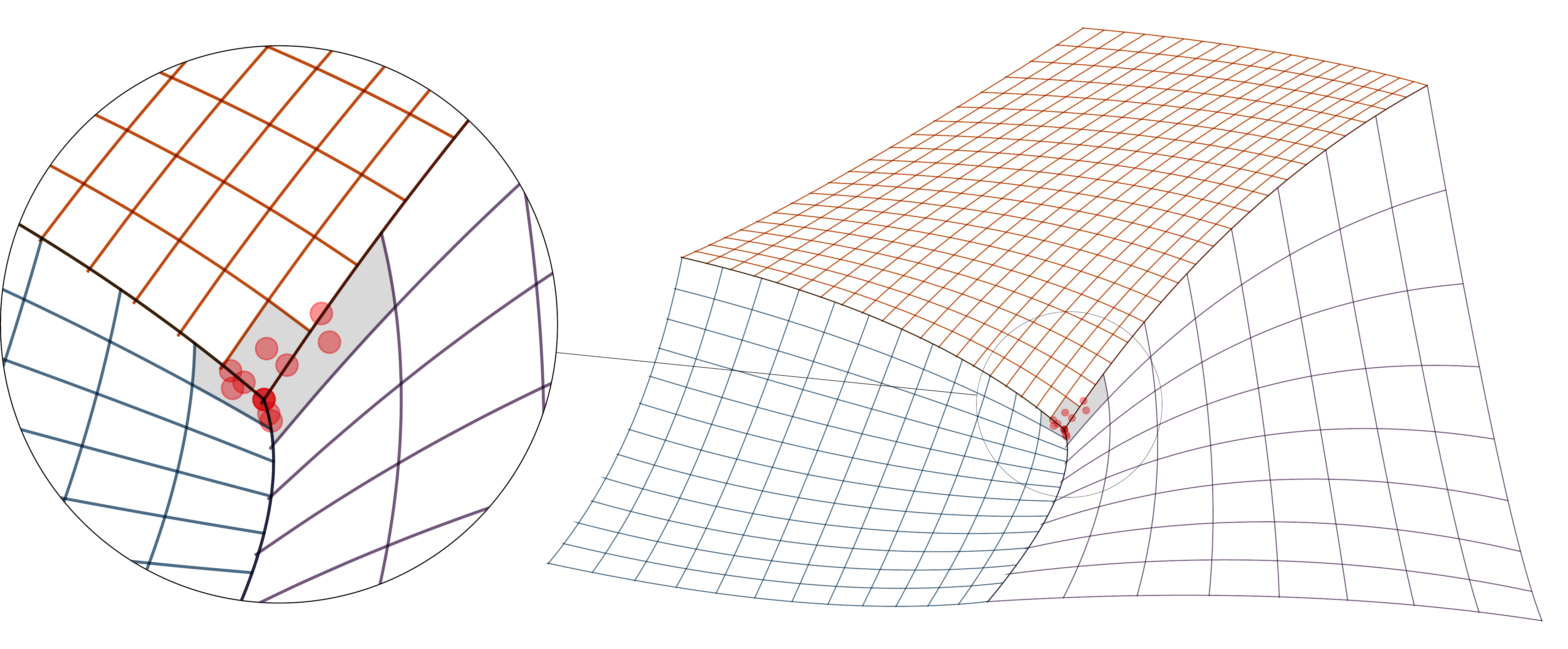}
	\caption{The degrees of freedom (red dots) involved in the constraint matrix $\mathbf{B}_v$ in section~\ref{sec:vertex_modification}.}\label{fig:cross_point_dof_modify}
\end{figure}

\begin{figure}[ht]
	\centering
	\includegraphics[width=.4\linewidth]{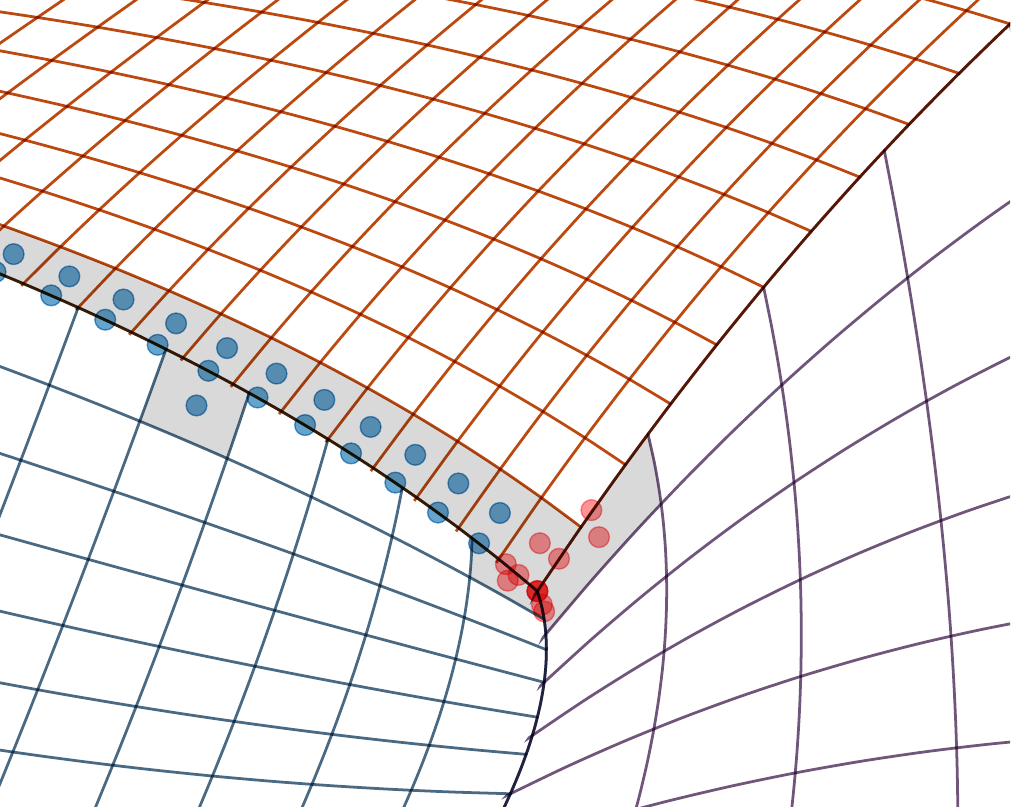}\qquad
	\includegraphics[width=.4\linewidth]{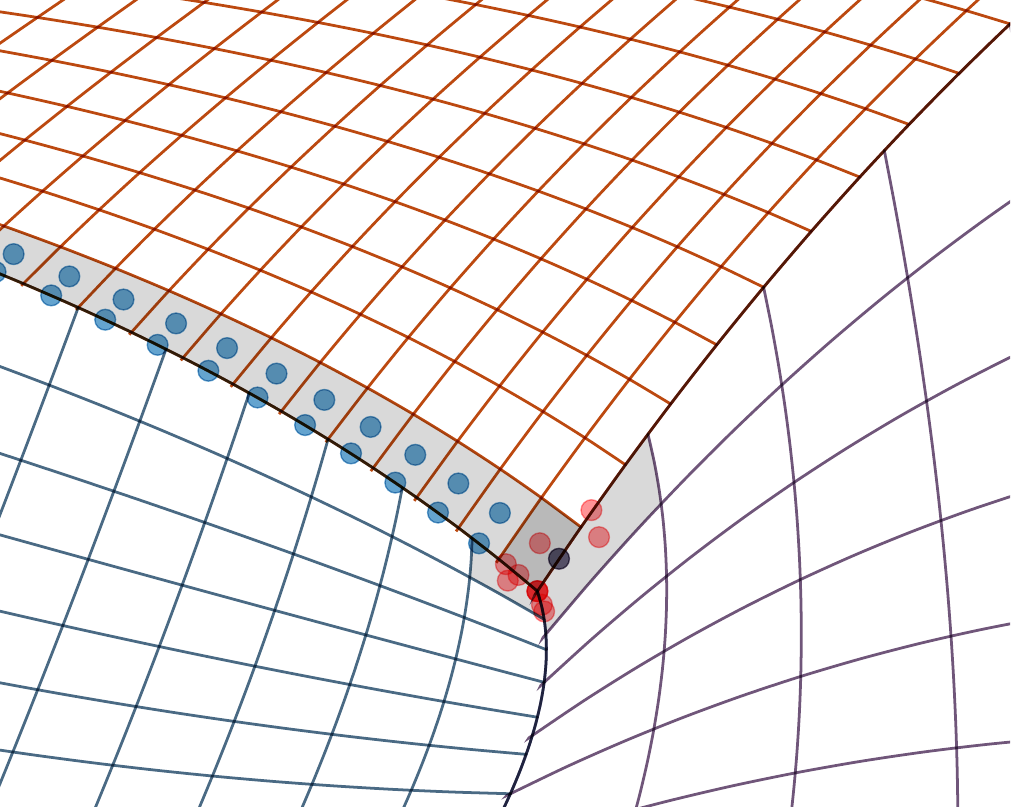}
	\caption{Activated degrees of freedom of the vector basis (blue) defined by the columns of $\mathbf{C}_b^\text{mod}$. Left: A vector basis classified as $\mathbf{C}_1$; Right: A vector basis classified as $\mathbf{C}_2$}\label{fig:overlap_nonoverlap_modify}
\end{figure}

In order to guarantee the well-posedness of the mortar formulation and to improve the approximation, we further apply two continuity constraints at each vertex
\begin{equation}
	\left\{\begin{split}
		&u_s(v) = u_m(v)\\
		&\frac{\partial u_s(v)}{\partial \xi_s} = \frac{\partial u_m(v)}{\partial \xi_s}\text{ if }\Gamma\parallel\eta_s\text{ or }\frac{\partial u_s(v)}{\partial \eta_s} = \frac{\partial u_m(v)}{\partial \eta_s}\text{ if }\Gamma\parallel\xi_s
	\end{split}\right.,\label{eq:vertex_constraint}
\end{equation}
where $v$ is the position of a vertex. The matrix form of the pointwise constraints~\eqref{eq:vertex_constraint} is denoted by $\mathbf{B}_v$. Hence, in the presence of in-domain vertices, the constraint matrix $\mathbf{B}_b$ is formed from both applying constraints weakly through the coarsened dual basis functions along each interface and by applying constraints strongly at each vertex as
\begin{equation}
	\mathbf{B}_b =
	\begin{bmatrix}
		\mathbf{B}_v \\
		\mathbf{B}^\text{mod}_b
	\end{bmatrix}.
\end{equation}
The null space of $\mathbf{B}_b$ is the intersection of the null space of $\mathbf{B}_v$ and the null space of $\mathbf{B}^\text{mod}_b$. As a result, $\mathbf{C}_b$ is the vector basis of the null space of $\mathbf{B}_v$ constructed from $\Ima\mathbf{C}^\text{mod}_b$. First, we split the column vectors of $C^\text{mod}_b$ into two matrices
\begin{equation}
	\begin{split}
		\mathbf{C}_1&:=\{\mathbf{v}\in\mathbf{C}^\text{mod}_b\colon{\mathbf{B}_v\mathbf{v}={0}}\},\\
		\mathbf{C}_2&:=\{\mathbf{v}\in\mathbf{C}^\text{mod}_b\colon{\mathbf{B}_v\mathbf{v}\neq{0}}\}.
	\end{split}
\end{equation}
An example of this split is given in Figure~\ref{fig:overlap_nonoverlap_modify}. $\mathbf{C}_1$ contains vectors of $\mathbf{C}^\text{mod}_b$ that are also in the null space of $\mathbf{B}_v$. Thus, they contribute to part of $\mathbf{C}_b$. The null space of $\mathbf{B}_v$ from $\Ima \mathbf{C}_2$ can be constructed as $\mathbf{C}_2\bar{\mathbf{C}}$, where $\bar{\mathbf{C}}$ is the vector basis of the null space of $\bar{\mathbf{B}} = \mathbf{B}_v\mathbf{C}_2$ and can be constructed through the factorization in~\eqref{eq:direct_kernel}. Since for each patch, only at most four degrees of freedom per vertex per patch are involved in the formulation of $\mathbf{B}_v$ (see Figure~\ref{fig:cross_point_dof_modify}), the number of vectors in $\mathbf{C}_2$ is very small and the cost of the factorization of $\bar{\mathbf{B}}$ is negligible. The vector basis of the null space of $\mathbf{B}_b$ can now be given as
\begin{equation}
	\mathbf{C}_b =
	\begin{bmatrix}
		\mathbf{C}_1 & \mathbf{C}_2\bar{\mathbf{C}}\label{eq:c_b_matrix}
	\end{bmatrix}.
\end{equation}

\begin{remark}
	The \Bezier dual basis can also be coarsened by replacing the three closest basis functions at each vertex by their summation. In other words, for a set of \Bezier dual basis functions $\{\hat{N}^B_i\}_{i=1}^n$, the coarsened \Bezier dual basis functions $\{\hat{N}^\text{mod}_i\}_{i=1}^{n-4}$ are defined as follows
	\begin{equation}
		\hat{N}^\text{mod}_i=
		\begin{cases}
			\hat{N}^B_1+\hat{N}^B_2+\hat{N}^B_3\quad           & i=1,              \\
			\hat{N}^B_{n}+\hat{N}^B_{n-1}+\hat{N}^B_{n-2}\quad & i=n-4,            \\
			\hat{N}^B_{i+3}\quad                               & \text{otherwise}.
		\end{cases}\label{eq:bezier_dual_modification}
	\end{equation}
	This modification preserves the biorthogonality relation~\eqref{eq:modify_biorthogonality}. However, as can be seen from numerical examples, the modified \Bezier dual basis demonstrate poor performance. Hence, in the next subsection, we will introduce a multi-patch treatment for \Bezier dual basis functions without resorting to a modification at each vertex.
\end{remark}

\subsubsection{\Bezier dual basis multi-patch treatment}\label{sec:original_dual_basis}

\begin{figure}[ht]
	\centering
	\includegraphics[width=.7\linewidth]{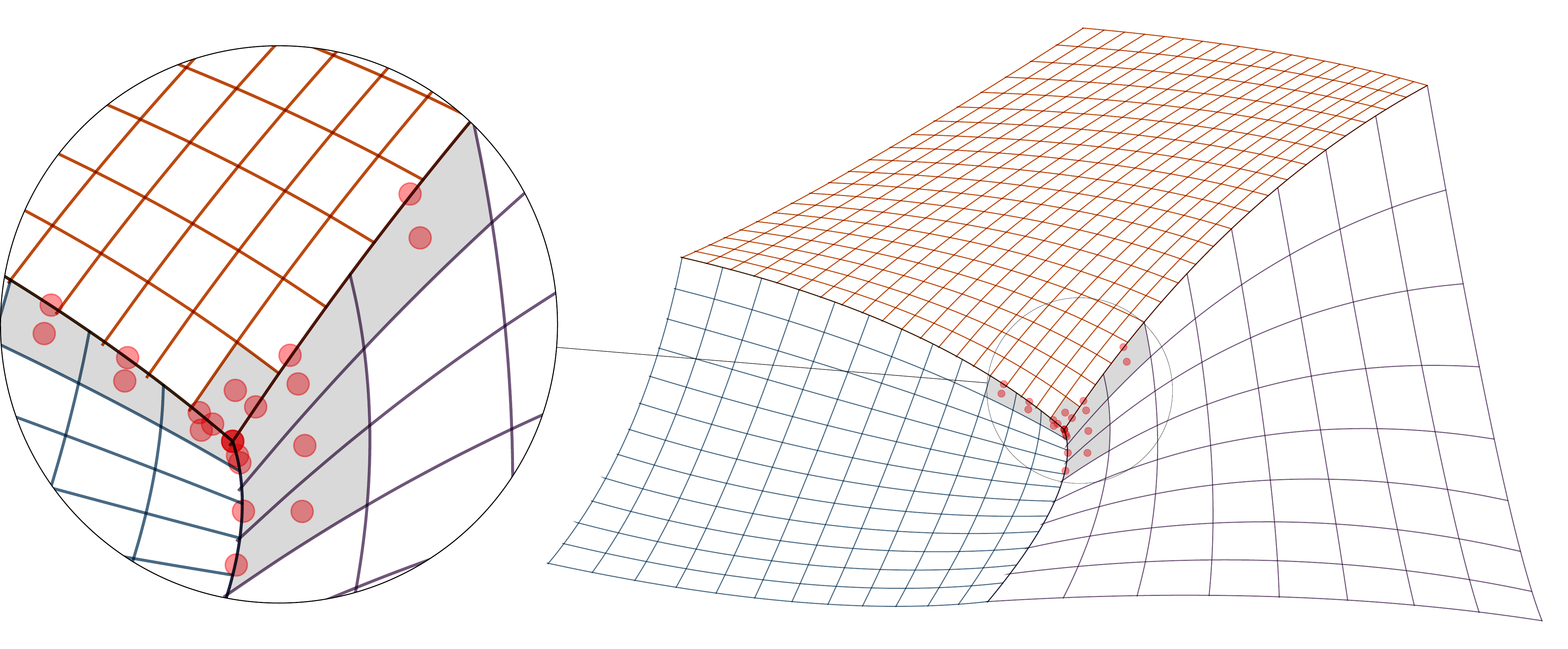}
	\caption{The degrees of freedom (red dots) involved in the constraint matrix $\mathbf{B}_v$ in section~\ref{sec:original_dual_basis} formed from quadratic B-splines.}\label{fig:cross_point_dof}
\end{figure}

\begin{figure}[ht]
	\centering
	\includegraphics[width=.4\linewidth]{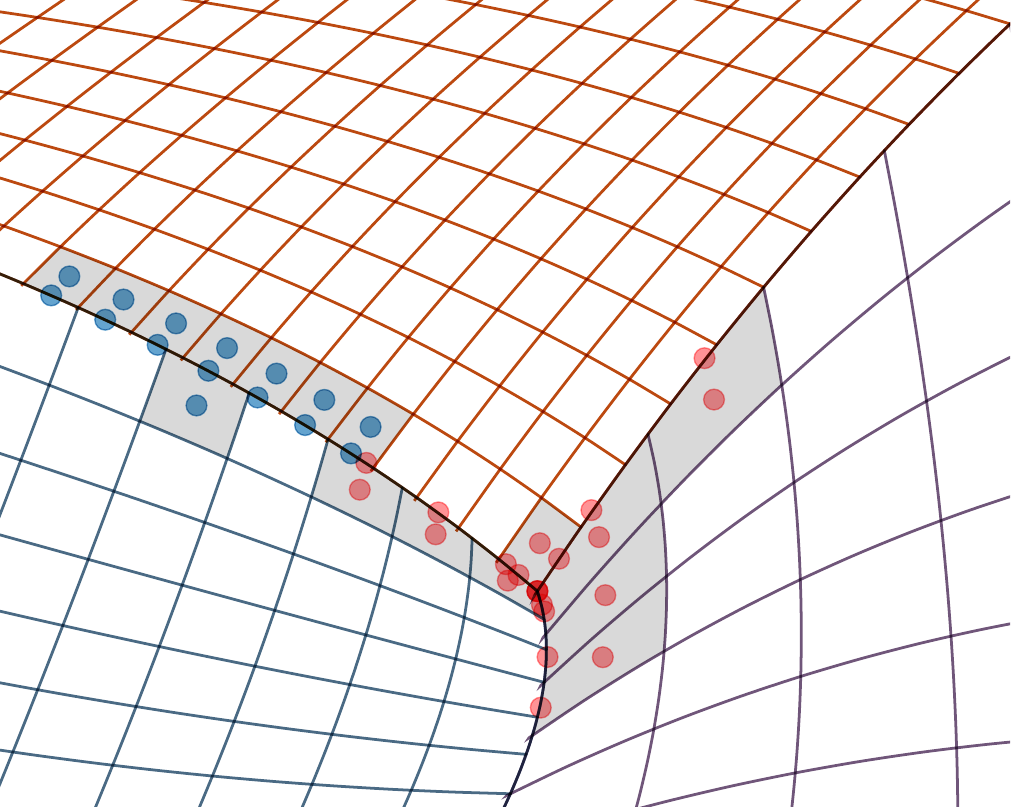}\qquad
	\includegraphics[width=.4\linewidth]{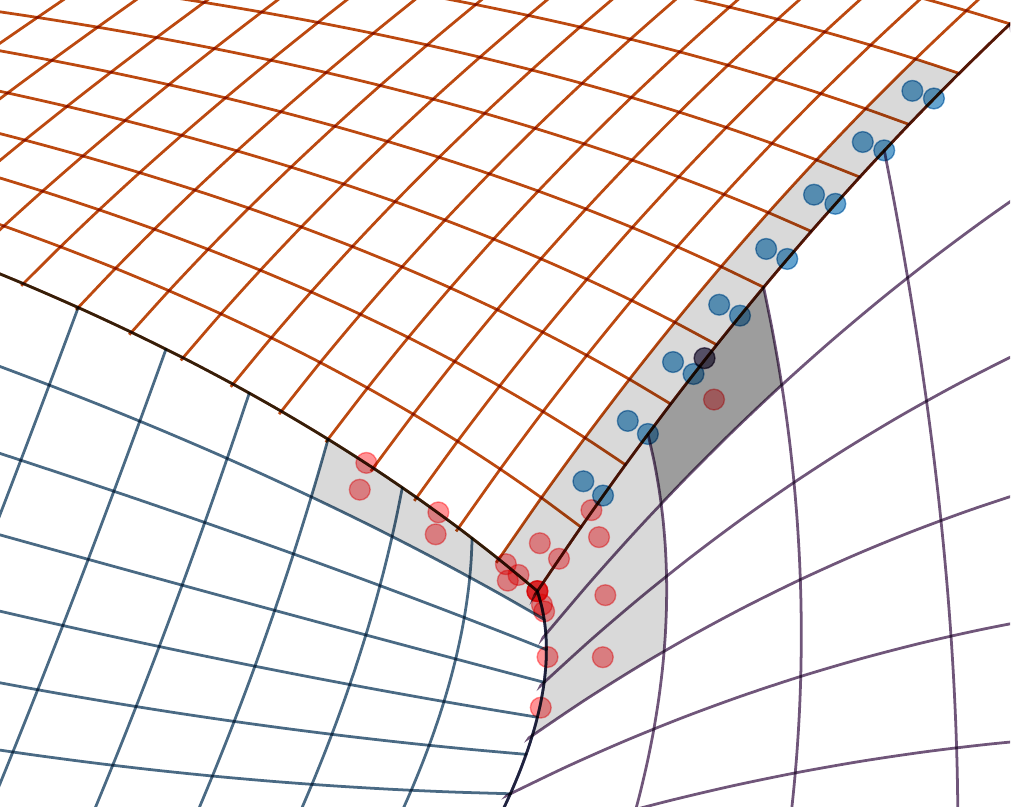}
	\caption{Activated degrees of freedom of the vector basis (blue) defined by columns of $\mathbf{C}_\text{inter}$. Left: A vector basis classified as $\mathbf{C}_1$; Right: A vector basis classified as $\mathbf{C}_2$}\label{fig:overlap_nonoverlap}
\end{figure}

In this approach, the basis functions in the discrete Lagrange multiplier spaces $\mathcal{M}^h_0$ and $\mathcal{M}^h_1$ are classified according to their proximity to a vertex. The basis functions in $\mathcal{M}^{h}_v$ are those in both $\mathcal{M}^h_0$ and $\mathcal{M}^h_1$ whose values and derivatives are non-zero at a vertex. The remaining basis functions are put in $\mathcal{M}^{h}_\text{inter}$. The constraint matrix can then be written as
\begin{equation}
	\mathbf{B}_b =
	\begin{bmatrix}
		\mathbf{B}_v \\
		\mathbf{B}_\text{inter}
	\end{bmatrix}
\end{equation}
where $\mathbf{B}_v$ is the matrix form of $b_0$ and $b_1$ restricted to the functions in $\mathcal{M}^{h}_v$ and $\mathbf{B}_\text{inter}$ is the matrix form of $b_0$ and $b_1$ restricted to the functions in $\mathcal{M}^{h}_\text{inter}$. Note that $\mathbf{B}_\text{inter}$ has the same structure as $\mathbf{B}_b^\text{mod}$ in Section~\ref{sec:vertex_modification}. The basis vectors of the null space of $\mathbf{B}_\text{int}$ can be constructed from~\eqref{eq:nullspace_modification} and is denoted by $\mathbf{C}_\text{inter}$. Using the approach outlined in Section~\ref{sec:vertex_modification}, we construct the basis vectors of the null space of $\mathbf{B}_v$ from $\Ima \mathbf{C}_\text{inter}$. $\mathbf{C}_\text{inter}$ can also be split into two submatrices $\mathbf{C}_1$ and $\mathbf{C}_2$. However, owing to the fact that $\mathbf{B}_v$ is constructed variationally, slightly more degrees of freedom are involved than in the approach described in Section~\ref{sec:vertex_modification} (see Figure~\ref{fig:cross_point_dof}). However, thanks to the locality of the \Bezier dual basis, the number of vectors in $\mathbf{C}_2$ is fixed. An example of this split is given in Figure~\ref{fig:overlap_nonoverlap}. Following the same procedure as described in Section~\ref{sec:vertex_modification}, we can construct $\mathbf{C}_b$ from Equation~\eqref{eq:c_b_matrix}.

\begin{remark}
	This approach can be directly extended to the global dual basis as well. However, as can be seen from numerical examples, the global dual basis without the coarsening procedure is not inf-sup stable.
\end{remark}

\section{Finite element error analysis}\label{sec:error_analysis}

In this section, we study the finite element approximation of~\eqref{eq:biharmonic_mixed}. Suppose that $\mathcal{X}_b^h\subset{\mathcal{X}_b}$, $\mathcal{K}_b^h\subset{\mathcal{K}_b}$, $\mathcal{M}_0^h\subset{\mathcal{M}_0}$ and $\mathcal{M}_1^h\subset{\mathcal{M}_1}$ are finite-dimensional linear subspaces of the spaces $\mathcal{X}_b$, $\mathcal{K}_b$, $\mathcal{M}_0$, and $\mathcal{M}_1$.
\begin{lemma}\label{aspt:bounded-operator}
	(bounded above) The bilinear functionals $a_b(\cdot,\cdot)$, $b_0(\cdot,\cdot)$ and $b_1(\cdot,\cdot)$ are all bounded; i.e., there exists positive constants $C_a$, $C_{b_0}$ and $C_{b_1}$ such that
	\begin{align}
		\begin{split}
			\vert{a_b{(u,v)}}\vert&\leq{C_a\|u\|_{H^2_*}\|v\|_{H^2_*}},\\
			\vert{b_0{(\mu_0,u)}}\vert&\leq{C_{b_0}\|\mu_0\|_{H^{-\frac{3}{2}}}\|u\|_{H^2_*}},\\
			\vert{b_1{(\mu_1,u)}}\vert&\leq{C_{b_1}\|\mu_1\|_{H^{-\frac{1}{2}}}\|u\|_{H^2_*}},
		\end{split}
		\quad\quad
		\begin{split}
			\forall{u,v}&\in{\mathcal{X}_b},\\
			\forall{\mu_0}\in\mathcal{M}_0, {u}&\in{\mathcal{X}_b},\\
			\forall{\mu_1}\in\mathcal{M}_1, {u}&\in{\mathcal{X}_b}.
		\end{split}
	\end{align}

	\begin{proof}
		The continuity of $a_b(\cdot,\cdot)$ follows from the Cauchy-Schwarz inequality, the continuities of $b_0(\cdot,\cdot)$ and $b_1(\cdot,\cdot)$ follow from the definition of the dual norm and the trace theorem.
	\end{proof}
\end{lemma}

\begin{lemma}\label{aspt:coercive}
	(bounded below) The bilinear functional $a_b(\cdot,\cdot)$ is coercive on the space $\mathcal{K}^h_b$ that is constructed from both the global dual basis functions and the \Bezier dual basis functions, i.e.,
	\begin{equation}
		\exists{c_a>0}\;\text{that is independent of the mesh size $h$ such that}\,\,\forall{v^h}\in\mathcal{K}_b^h,\, a_b(v^h,v^h)\geq{c_a\|{v^h}\|_{H^2_*}}.\label{eq:coercive}
	\end{equation}
	\begin{proof}
		We introduce the function space $\bar{\mathcal{X}}_b$ as all functions in $\mathcal{X}_b$ that satisfy the pointwise constraint~\eqref{eq:vertex_constraint} on each vertex and are also $C^\infty$ in each patch. Then we have the following inclusion, $\mathcal{K}^h_b\subset \bar{\mathcal{X}}_b \subset  \mathcal{X}_b$. For the case of the \Bezier dual basis, the inclusion comes from~\eqref{eq:boundary_interpolation}; for the case of the global dual basis, the inclusion comes from the pointwise constraints at the vertex. We consider the coercivity of $a_b(\cdot,\cdot)$ on $\bar{\mathcal{X}}_b$. Suppose that this statement is false, then there exists a sequence $\{v_n\}_{n=1}^{\infty}\in \bar{\mathcal{X}}_b$ such that
		\begin{equation}
			\sum_{k=1}^K\vert v_n \vert_{H^2(\Omega_k)}\leq \frac{1}{n}\;\text{ and }\;\sum_{k=1}^K\| v_n \|_{H^2(\Omega_k)}=1.\label{eq:countrary}
		\end{equation}
		For any $1\leq k \leq K$, the sequence $\left\{ v_n\vert_{\Omega_k} \right\}_{n=1}^\infty$ is bounded in $H^2(\Omega_k)$. Hence, there exists a subsequence $\left\{ (v_n)_m\vert_{\Omega_k} \right\}_{m=1}^\infty$ such that
		\begin{equation}
			(v_n)_m\vert_{\Omega_k} \rightharpoonup v\vert_{\Omega_k} \text{weakly in } H^2(\Omega_k).
		\end{equation}
		By the Rellich-Kondrachov theorem, $\left\{ (v_n)_m\vert_{\Omega_k} \right\}_{m=1}^\infty$ converges to $v\vert_{\Omega_k}$ strongly in $H^1(\Omega_k)$. In other words,
		\begin{equation}
			\lim_{m\to\infty}\| (v_n)_m \vert_{\Omega_k} -v \vert_{\Omega_k}\|_{H^1(\Omega_k)}=0.
		\end{equation}
		From $\vert v_n \vert_{H^2(\Omega_k)}\rightarrow 0$, we have that $\left\{ (v_n)_m\vert_{\Omega_k} \right\}_{m=1}^\infty$ is a Cauchy sequence in $H^2(\Omega_k)$. Thus, we have that $v \vert_{\Omega_k}\in H^2(\Omega_k)$ and
		\begin{equation}
			\lim_{m\to\infty}\| (v_n)_m \vert_{\Omega_k} -v \vert_{\Omega_k}\|_{H^2(\Omega_k)}=0.
		\end{equation}
		From the approximation theory,
		\begin{equation}
			\inf_{r^k_n\in\mathcal{P}^1(\Omega_k)}\|v_{n}-r^k_n\|_{H^1(\Omega_k)}\leq C \vert v_{n} \vert_{H^2(\Omega_k)} \leq \frac{C}{n}
		\end{equation}
		where $\mathcal{P}^1(\Omega_k)$ is the $1^\text{st}$ order polynomial space on $\Omega_k$. Hence,
		\begin{equation}
			\| (r^k_n)_m-v\|_{H^1(\Omega_k)} \leq \| (r^k_n)_m-(v_n)_m\|_{H^1(\Omega_k)}+\| (v_n)_m-v\|_{H^1(\Omega_k)} \rightarrow 0.
		\end{equation}
		This means $v\vert_{\Omega_k}$ is a linear function on $\Omega_k$. If $\Omega_k$ is a boundary element (i.e. $\partial \Omega_k\bigcap \partial\Omega\neq\emptyset$), the only linear functions that satisfies the boundary condition is the zero function. Hence,
		\begin{equation}
			\|(v_{n})_m\|_{H^2(\Omega_k)}=\|(v_{n})_m-v\|_{H^2(\Omega_k)}\rightarrow 0.
		\end{equation}
		If patch $\Omega_l$ is coupled with $\Omega_k$ along the intersection $\Gamma_{kl}\parallel \eta_k$, then on both ends $a,b$ of $\Gamma_{kl}$ we have
		\begin{equation}
			v(a) = 0,\; v(b) = 0,\; \frac{\partial v(a)}{\partial \xi_k} = 0,\;\text{ and }\frac{\partial v(b)}{\partial \xi_k} = 0.\label{eq:boundary_conditions}
		\end{equation}
		Hence, $v\vert_{\Omega_l}=0$. Similar arguments can be applied to all patches, leading to
		\begin{equation}
			\sum_{k=1}^K\| (v_n)_m \|_{H^2(\Omega_k)}\rightarrow 0.
		\end{equation}
		This is inconsistent with~\eqref{eq:countrary}. As a result, \eqref{eq:coercive} holds.
	\end{proof}
\end{lemma}

\begin{theorem}
	There exists a unique solution $u^h\in{\mathcal{K}_b^h}$ that satisfies~\eqref{eq:weak-form} for all $v^h\in{\mathcal{K}_b^h}$, with $\mathcal{K}_b^h$ constructed using either the \Bezier or global dual basis.
	\begin{proof}
		Thanks to Lemma~\ref{aspt:bounded-operator} and Lemma~\ref{aspt:coercive}, the well-posedness of problem~\eqref{eq:weak-form} in $\mathcal{K}_b^h$ follows from the Lax-Milgram theorem~\cite{brenner_mathematical_2007}.\par
	\end{proof}
\end{theorem}

\begin{lemma}
	(Strang's lemma) Let $u\in H^2_0(\Omega)$ satisfy problem~\eqref{eq:weak-form} for all $v\in H^2_0(\Omega)$, then the error between $u$ and $u^h$ is given by
	\begin{equation}
		\|u-u^h\|_{H^2_*}\leq{\left(1+\frac{C_a}{c_a}\right)}\inf_{v^h\in{\mathcal{K}_b^h}}\|u-v^h\|_{H^2_*}+\frac{1}{c_a}\sup_{w^h\in \mathcal{K}_b^h\backslash \{0\}}\frac{\vert a_b(u-u^h,w^h) \vert }{\|w^h\|_{H^2_*}},
	\end{equation}
	where the first term on the righthand side is often called the approximation error and the second term is often called the consistency error.
	\begin{proof}
		See~\cite{brenner_mathematical_2007}.
	\end{proof}
\end{lemma}
From Strang's lemma, we may now obtain the following result by expanding the term $a_b(u-u^h,w^h)$.
\begin{theorem}\label{thm:fea-approx}
	The error between $u$ and $u^h$ is given by
	\begin{equation}
		\begin{split}
			\|u-u^h\|_{H^2_*}\leq &{\left(1+\frac{C_a}{c_a}\right)}\inf_{v^h\in{\mathcal{K}_b^h}}\|u-v^h\|_{H^2_*}+\frac{c_b}{c_a}\sum_{\Gamma\in \mathbf{S}}B,
		\end{split}
		\label{eq:fem_approximation}
	\end{equation}
	with
	\begin{equation}
		B= \left\{\begin{split}
			\inf_{\lambda_0^h\in \mathcal{M}_0^h}\| \frac{\partial \Delta u}{\partial \mathbf{n}} + r_{\eta_s} - \lambda_0^h \|_{H^{-\frac{3}{2}}(\Gamma)}+ \inf_{\lambda_1^h\in \mathcal{M}_1^h} \| \Delta u \frac{\partial \xi_s}{\partial \mathbf{n}} - \lambda_1^h \|_{H^{-\frac{1}{2}}(\Gamma)}\quad\quad &\Gamma\parallel \eta_s,\\
			\inf_{\lambda_0^h\in \mathcal{M}_0^h}\| \frac{\partial \Delta u}{\partial \mathbf{n}} + r_{\xi_s} - \lambda_0^h \|_{H^{-\frac{3}{2}}(\Gamma)}+ \inf_{\lambda_1^h\in \mathcal{M}_1^h} \| \Delta u \frac{\partial \eta_s}{\partial \mathbf{n}} - \lambda_1^h \|_{H^{-\frac{1}{2}}(\Gamma)}\quad \quad&\Gamma\parallel \xi_s,
		\end{split}\right.
	\end{equation}
	where
	\begin{equation}
		r_{\eta_s} = \partial_{\eta_s}\left(\Delta u \frac{\partial \eta_s}{\partial \mathbf{n}}\vert \Gamma' \vert\right)\frac{1}{\vert \Gamma' \vert},\quad r_{\xi_s} = \partial_{\xi_s}\left(\Delta u \frac{\partial \xi_s}{\partial \mathbf{n}}\vert \Gamma' \vert\right)\frac{1}{\vert \Gamma' \vert},
	\end{equation}
	with
	\begin{equation}
		\begin{split}
			\vert \Gamma' \vert = \sqrt{\frac{\partial x}{\partial \eta_s}^2+\frac{\partial y}{\partial \eta_s}^2}\quad\quad &\Gamma\parallel \eta_s\\
			\vert \Gamma' \vert = \sqrt{\frac{\partial x}{\partial \xi_s}^2+\frac{\partial y}{\partial \xi_s}^2}\quad\quad &\Gamma\parallel \xi_s
		\end{split}
	\end{equation}
	and $c_b$ is a constant independent of the mesh size.
	\begin{proof}
		Here we only discuss the stituation where $\Gamma\parallel \eta_s$. By Green's theorem, we have that
		\begin{equation}
			\begin{split}
				a_b(u-u^h,w^h)&=a_b(u, w^h)-\langle f,w^h \rangle_{\Omega}\\
				&=\langle \Delta^2u,w^h \rangle_{\Omega}+\sum_{\Gamma\in \mathbf{S}}\left(\int_\Gamma\Delta u\left[\frac{\partial{w^h}}{\partial \mathbf{n}}\right]_\Gamma d\Gamma - \int_\Gamma\frac{\partial \Delta u}{\partial \mathbf{n}}\left[{w^h}\right]_\Gamma d\Gamma\right)-\langle f,w^h \rangle_{\Omega}\\
				&=\sum_{\Gamma\in \mathbf{S}}\left(\int_\Gamma\Delta u\left[\frac{\partial{w^h}}{\partial \mathbf{n}}\right]_\Gamma d\Gamma - \int_\Gamma\frac{\partial \Delta u}{\partial \mathbf{n}}\left[{w^h}\right]_\Gamma d\Gamma\right)\\
				&=\sum_{\Gamma\in \mathbf{S}}\left(\int_\Gamma\Delta u \frac{\partial \xi_s}{\partial \mathbf{n}}\left[\frac{\partial{w^h}}{\partial \xi_s}\right]_\Gamma d\Gamma+\int_\Gamma\Delta u \frac{\partial \eta_s}{\partial \mathbf{n}}\left[\frac{\partial{w^h}}{\partial \eta_s}\right]_\Gamma d\Gamma-\int_\Gamma\frac{\partial \Delta u}{\partial \mathbf{n}}\left[{w^h}\right]_\Gamma d\Gamma\right),
			\end{split}
		\end{equation}
		where the second term can be rewritten as
		\begin{equation}
			\begin{split}
				\int_\Gamma\Delta u \frac{\partial \eta_s}{\partial \mathbf{n}}\left[\frac{\partial{w^h}}{\partial \eta_s}\right]_\Gamma d\Gamma &= \int_0^1\Delta u \frac{\partial \eta_s}{\partial \mathbf{n}}\left[\frac{\partial{w^h}}{\partial \eta_s}\right]_\Gamma \vert \Gamma' \vert d\eta_s\\
				&=-\int_0^1\partial_{\eta_s}\left(\Delta u \frac{\partial \eta_s}{\partial \mathbf{n}}\vert \Gamma' \vert\right)\left[{w^h}\right]_\Gamma d\eta_s\\
				&=-\int_\Gamma r_{\eta_s}\left[{w^h}\right]_\Gamma d\Gamma.
			\end{split}
		\end{equation}
		Hence,
		\begin{equation}
			a_b(u-u^h,w^h)=\sum_{\Gamma\in \mathbf{S}}\left(\int_\Gamma\Delta u \frac{\partial \xi_s}{\partial \mathbf{n}}\left[\frac{\partial{w^h}}{\partial \xi_s}\right]_\Gamma d\Gamma-\int_\Gamma\left(\frac{\partial \Delta u}{\partial \mathbf{n}} + r_{\eta_s}\right)\left[{w^h}\right]_\Gamma d\Gamma\right).
		\end{equation}
		Next, using the constraints in the definition of the function space $\mathcal{K}^h_b$, we have that, for any $\lambda_0^h, \lambda_1^h \in \mathcal{M}_0^h\times\mathcal{M}_1^h$,
		\begin{equation}
			\begin{split}
				\sum_{\Gamma\in \mathbf{S}}\left(\int_\Gamma\left(\frac{\partial \Delta u}{\partial \mathbf{n}} + r_{\eta_s}  \right)\left[{w^h}\right]_\Gamma d\Gamma\right) &= \sum_{\Gamma\in \mathbf{S}}\left(\int_\Gamma\left(\frac{\partial \Delta u}{\partial \mathbf{n}} + r_{\eta_s} - \lambda_0^h \right)\left[{w^h}\right]_\Gamma d\Gamma\right)\\
				&\leq \sum_{\Gamma\in \mathbf{S}} \| \frac{\partial \Delta u}{\partial \mathbf{n}} + r_{\eta_s} - \lambda_0^h \|_{H^{-\frac{3}{2}}(\Gamma)}\left(\|w^h_s\|_{H^{\frac{3}{2}}(\Gamma)}+\|w^h_m\|_{H^{\frac{3}{2}}(\Gamma)}\right),
			\end{split}
		\end{equation}
		and
		\begin{equation}
			\begin{split}
				\sum_{\Gamma\in \mathbf{S}}\left(\int_\Gamma\Delta u \frac{\partial \xi_s}{\partial \mathbf{n}}\left[\frac{\partial{w^h}}{\partial \xi_s}\right]_\Gamma d\Gamma\right) &= \sum_{\Gamma\in \mathbf{S}}\left(\int_\Gamma\left(\Delta u \frac{\partial \xi_s}{\partial \mathbf{n}} - \lambda_1^h \right)\left[\frac{\partial{w^h}}{\partial \xi_s}\right]_\Gamma d\Gamma\right)\\
				&\leq \sum_{\Gamma\in \mathbf{S}} \| \Delta u \frac{\partial \xi_s}{\partial \mathbf{n}} - \lambda_1^h \|_{H^{-\frac{1}{2}}(\Gamma)}\left(\|\frac{\partial w^h_s}{\partial \xi_s}\|_{H^{\frac{1}{2}}(\Gamma)}+\|\frac{\partial w^h_m}{\partial \xi_s}\|_{H^{\frac{1}{2}}(\Gamma)}\right).
			\end{split}
		\end{equation}
		From the trace theorem, we obtain the following conclusion:
		\begin{equation}
			\sup_{w^h\in \mathcal{K}_b^h\backslash \{0\}}\frac{\vert a_b(u-u^h,w^h) \vert }{\|w^h\|_{H^2_*(\Omega)}}\leq c_b\sum_{\Gamma\in \mathbf{S}}\left( \inf_{\lambda_0^h\in \mathcal{M}_0^h}\| \frac{\partial \Delta u}{\partial \mathbf{n}} + r_{\eta_s} - \lambda_0^h \|_{H^{-\frac{3}{2}}(\Gamma)}+ \inf_{\lambda_1^h\in \mathcal{M}_1^h} \| \Delta u \frac{\partial \xi_s}{\partial \mathbf{n}} - \lambda_1^h \|_{H^{-\frac{1}{2}}(\Gamma)}\right).\label{eq:consistency_error}
		\end{equation}
	\end{proof}
\end{theorem}
Hence, the error of finite element approximations in the broken $H^2_*(\Omega)$ norm are bounded by the best $H^2(\Omega)$ approximation of $u^h\in\mathcal{K}_b^h$ and the best approximations of $\mu_0^h\in{\mathcal{M}_0^h}$, $\mu_1^h\in{\mathcal{M}_1^h}$ in $H^{-\frac{3}{2}}(\Gamma)$ and $H^{-\frac{1}{2}}(\Gamma)$, respectively.

Now, let $P_0$ and $P_1$ be $L^2$ projection operators in $\mathcal{M}^h_0$ and $\mathcal{M}^h_1$. Then the approximation of $u$ in the fractional Sobolev space $H^{s-\frac{1}{2}}(\Gamma)$ is given by~\cite{lamichhane2006higher, bernardi_domain_1993}
\begin{equation}
	\|u-P_iu\|_{L^2(\Gamma)}\leq Ch^{\min\{s,(p_i) + 1\}-\frac{1}{2}}\|u\|_{H^{s-\frac{1}{2}}(\Gamma)}, \quad \; i\in\left\{ 0, 1 \right\},\label{eq:estimates}
\end{equation}
where $p_0$ and $p_1$ are the polynomial order that $\mathcal{M}^h_0$ or $\mathcal{M}^h_1$ reproduce, respectively.

Now, recalling the standard Aubin-Nitsche duality argument~\cite{tagliabue2014isogeometric, strang1973analysis} and applying estimates~\eqref{eq:estimates} and the trace theorem, we get
\begin{equation}
	\begin{cases}
		\|u-P_0u\|_{H^{-\frac{3}{2}}(\partial \Omega_k)}\leq C h_k^{\frac{3}{2}}\|u-P_0u\|_{L^2(\partial \Omega_k)}\leq Ch_k^{\min\{s,p_0+1\}+1}\|u\|_{H^{s-\frac{1}{2}}(\partial \Omega_k)}\leq Ch_k^{\min\{s,p_0+1\}+1}\|u\|_{H^{s}(\Omega_k)}, \\
		\|u-P_1u\|_{H^{-\frac{1}{2}}(\partial \Omega_k)}\leq C h_k^{\frac{1}{2}}\|u-P_1u\|_{L^2(\partial \Omega_k)}\leq Ch_k^{\min\{s,p_1+1\}}\|u\|_{H^{s-\frac{1}{2}}(\partial \Omega_k)}\leq Ch_k^{\min\{s,p_1+1\}}\|u\|_{H^{s}(\Omega_k)}.
	\end{cases}
\end{equation}
Although the approximation power of $\mathcal{K}_b^h$ remains unknown, the ability of $\mathcal{X}_b^h$ to approximate functions $u\in H^s(\Omega)$ is given by
\begin{equation}
	\inf_{v^h\in{\mathcal{X}^h_b}} \|u-v^h\|_{H^l(\Omega_k)}\leq{Ch_k^{\min\{s,p+1\}-l}}\|u\|_{H^{s}(\Omega_k)},
\end{equation}
where $p$ is the polynomial order that $\mathcal{X}^h_b$ reproduces. In order to find the best approximation error of $\mathcal{K}_b^h$, we need the following assumption:
\begin{assumption}\label{aspt:inf-sup}
	(inf-sup) Assume that the bilinear functionals $b_0(\cdot,\cdot)$ and $b_1(\cdot,\cdot)$ are inf-sup stable in the discretized formulation, i.e., there exist positive constants $\beta_0$ and $\beta_1$ independent of the mesh size such that
	\begin{align}
		\inf_{\mu_0^h\in\mathcal{M}_0^h}\sup_{u^h\in\mathcal{X}_b^h\backslash \{0\}}\frac{{b_0\left({\mu_0^h,u}\right)}}{\|u^h\|_{H^2_*}\|\mu_0^h\|_{H^{-\frac{3}{2}}}}\geq{\beta_0}, \\
		\inf_{\mu_1^h\in\mathcal{M}_1^h}\sup_{u^h\in\mathcal{X}_b^h\backslash \{0\}}\frac{{b_1\left({\mu_0^h,u}\right)}}{\|u^h\|_{H^2_*}\|\mu_1^h\|_{H^{-\frac{1}{2}}}}\geq{\beta_1}.
	\end{align}
\end{assumption}
Now we may bound the best approximation error of $\mathcal{K}_b^h$ by the best approximation error of $\mathcal{X}_b^h$ via the following result.
\begin{theorem}
	Under Lemma~\ref{aspt:bounded-operator} and Assumption~\ref{aspt:inf-sup}, we have that, for any $u\in\mathcal{K}_b$,
	\begin{equation}
		\inf_{v^h\in{\mathcal{K}_b^h}}\|{u-v^h}\|_{H^2_*}\leq\left({1+\frac{C_b}{\beta}}\right)\inf_{w^h\in{\mathcal{X}_b^h}}\|{u-w^h}\|_{H^2_*}
	\end{equation}
	where $\beta=\min\left({\beta_{0},\beta_{1}}\right)$, $C_b=\max\left({C_{b_0},C_{b_1}}\right)$.
	\begin{proof}
		See~\cite{brenner_mathematical_2007} or~\cite{boffi_mixed_2013}.
	\end{proof}
\end{theorem}

The optimality of $u^h\in\mathcal{K}_b^h$ in $H^2_*$ requires the \textit{inf-sup} stability of the bilinear functionals $b_0$ and $b_1$. The analytical study of the \textit{inf-sup} stability of these functionals is beyond the scope of this paper. Instead, we demonstrate the approximation ability of $\mathcal{K}_b^h$ by directly conducting $H^2_*$ projection in different numerical examples. We may now give the final estimate:
\begin{theorem}\label{thm:approximation-of-bezier-formulation}
	Given Assumption~\ref{aspt:inf-sup}, we have that, on a smooth discretization i.e. $\mathbf{F}_i\in \left(C^\infty(\Omega_i)\right)^2$, for any $u\in H^s(\Omega)$,
	\begin{equation}
		\|u-u^h\|_{H^2_*(\Omega)}^2 \leq C \sum_{k=1}^K h_k^{2\sigma}\|u\|^2_{H^s(\Omega_k)},
	\end{equation}
	where $\sigma=\min\{s-2,p-1,p_0+2,p_1+1\}$.
\end{theorem}
Hence, for a smooth solution $u$, the optimality of the finite element approximation of the proposed method requires that $p_0\geq p-3$ and $p_1\geq p-2$.

\begin{remark}
	Theorem~\ref{thm:fea-approx} gives meaning to the Lagrange multipliers. In other words,
	\begin{equation}
		\left\{\begin{split}
			\lambda_0 &= \frac{\partial \Delta u}{\partial \mathbf{n}} + r_{\eta_s},\\
			\lambda_1 &= \Delta u \frac{\partial \xi_s}{\partial \mathbf{n}},
		\end{split}\right.\quad \Gamma\parallel \eta_s
		\quad\quad\text{ or }\quad\quad
		\left\{\begin{split}
			\lambda_0 &= \frac{\partial \Delta u}{\partial \mathbf{n}} + r_{\xi_s},\\
			\lambda_1 &= \Delta u \frac{\partial \eta_s}{\partial \mathbf{n}},
		\end{split}\right.\quad \Gamma\parallel \xi_s.
	\end{equation}
	If we apply the conventional $C^1$ constraints $\left[u\right]_\Gamma=0$ and $\left[\frac{\partial u}{\partial \mathbf{n}}\right]_\Gamma=0$, then
	\begin{equation}
		\left\{\begin{split}
			\lambda_0 &= \frac{\partial \Delta u}{\partial \mathbf{n}},\\
			\lambda_1 &= \Delta u.
		\end{split}\right.
	\end{equation}
\end{remark}

\section{Numerical examples}\label{sec:numerical_examples}

In this section, we investigate the performance of the proposed method using both global and \Bezier dual basis functions on several challenging benchmark problems. The first example is a two-patch coupling problem, where different discretizations and parameterizations are studied. To investigate the influence of the number of in-domain vertices, a three-patch coupling and a five-patch coupling problem are studied. The approximation errors of each benchmark problem are studied through the use of $H^2_*$ projection. To also demonstrate the advantages of the proposed coupling method for $2^\text{nd}$ order problems, we consider the transverse vibrations of an elastic membrane on a nine-patch square domain. For the two-patch cases, results computed using the global dual basis are labeled G-$Q_i$ while results computed using the \Bezier dual basis are labeled B-$Q_i$. The subscript $i$ denotes the degree. For the multi-patch cases, the methods tested are summarized in Table~\ref{tab:methods}. All problems are solved with the conjugate gradient module in Eigen~\cite{eigenweb}.

\begin{table}[h]
	\setlength{\tabcolsep}{10pt}
	\def\arraystretch{1.5}
	\caption{ Number of turns and distance between top and bottom.}
	\begin{tabularx}{\textwidth}{>{\hsize=.1\textwidth}X>{\hsize=.6\textwidth}X>{\hsize=.2\textwidth}X}
		\hline
		Method   & Description                                                                                                                                                                           & Optimality                                                    \\
		\hline
		MG-$Q_i$ & Using the global dual basis with the multi-patch treatment described in Section~\ref{sec:vertex_modification}.                                                                        & Yes                                                           \\
		MB-$Q_i$ & Using the \Bezier dual basis with the multi-patch treatment described in Section~\ref{sec:vertex_modification}. The coarsening procedure follows~\eqref{eq:bezier_dual_modification}. & No, early domination of the consistency error.                \\
		OG-$Q_i$ & Using the global dual basis with the multi-patch treatment described in Section~\ref{sec:original_dual_basis}.                                                                        & No, the formulation is not \textit{inf-sup} stable.           \\
		OB-$Q_i$ & Using the \Bezier dual basis with the multi-patch treatment described in Section~\ref{sec:original_dual_basis}.                                                                       & No, but the domination of the consistency error is postponed. \\
		\hline
	\end{tabularx}
	\label{tab:methods}
\end{table}

\subsection{The approximation power of dual basis functions}
\begin{figure}[ht]
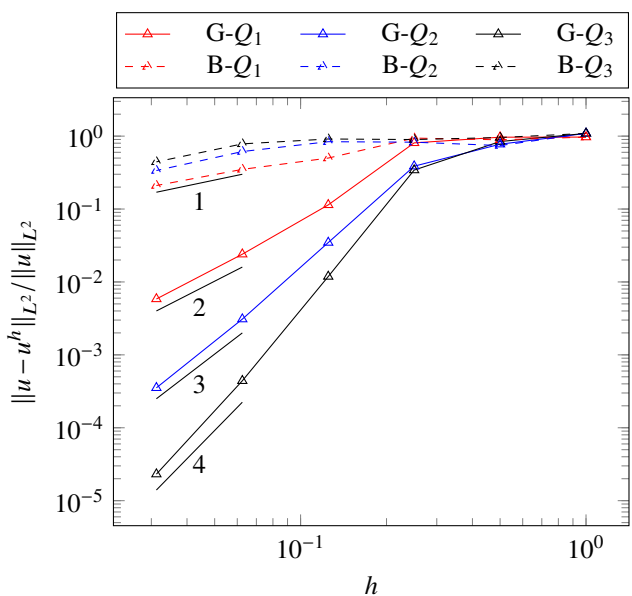

	\center
	\includestandalone[scale=1]{l_2_error}
	\caption{A comparison of the convergence rates for global and \Bezier dual basis functions for the $L^2$ projection of a sine function $u(x)=\sin (4\pi x)$ onto the domain $\left[ 0, 1\right]$.}
	\label{fig:l_2_error}
\end{figure}

We study the approximation power of both global and \Bezier dual basis functions by considering the $L^2$ projection of a sinusoid function $u(x)=\sin (4\pi x)$ onto the domain $\left[0, 1\right]$. The $L^2$-error of the \Bezier dual basis is compared to that of the global dual basis in Figure~\ref{fig:l_2_error}. It can be seen that the global dual basis functions, as predicted, converge optimally while the convergence of \Bezier dual basis functions is sub-optimal and the error is $\mathcal{O}(h)$ for all tested polynomial degrees. \par

\begin{figure}[ht]
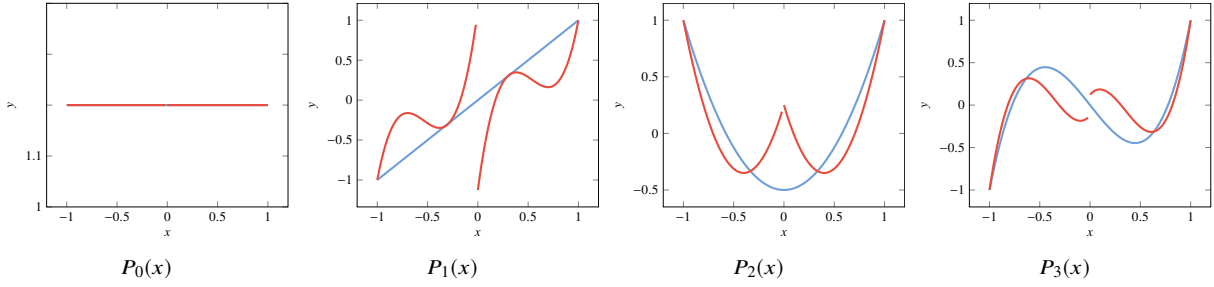

	\captionsetup[subfigure]{labelformat=empty, font = footnotesize}
	\centering
	\begin{subfigure}[b]{0.24\textwidth}
		\centering
		\includestandalone[scale=.5]{Q3_p=0}
		\caption{$P_0(x)$}
	\end{subfigure}
	\begin{subfigure}[b]{0.24\textwidth}
		\centering
		\includestandalone[scale=.5]{Q3_p=1}
		\caption{$P_1(x)$}
	\end{subfigure}
	\begin{subfigure}[b]{0.24\textwidth}
		\centering
		\includestandalone[scale=.5]{Q3_p=2}
		\caption{$P_2(x)$}
	\end{subfigure}
	\begin{subfigure}[b]{0.24\textwidth}
		\centering
		\includestandalone[scale=.5]{Q3_p=3}
		\caption{$P_3(x)$}
	\end{subfigure}
	\caption{The Legendre polynomials (\protect\blueline) and the corresponding best $L^2$ approximations (\protect\redline) by $3^\text{rd}$ order B\'ezier dual basis functions defined on a two element domain. B\'ezier dual basis functions can only replicate the constant function.}
	\label{fig:polynomial_completeness}
\end{figure}

To better understand the cause of the poor approximation of the \Bezier dual basis, we consider the polynomial reproduction properties of the \Bezier dual basis. The domain $\Omega$ is uniformly partitioned into two elements since the B\'ezier dual basis is equivalent to the global dual basis on a one element domain. The $L^2$ approximation of the $n^{th}$ order Legendre polynomial is evaluated. From this test, we can see that the \Bezier dual basis functions can only reproduce the constant function. The result for the $3^\text{rd}$ order \Bezier dual basis is shown in Figure~\ref{fig:polynomial_completeness}. There are significant discrepancies between the dual approximations and the corresponding polynomials for all Legendre polynomials except the constant function. From the approximation theory, to achieve $p^\text{th}$ order convergence rates requires the reproduction of polynomials up to $p-1^\text{th}$ order. This explains the sub-optimality of the \Bezier dual basis in $L^2$ projection. Moreover, for the same mesh, the error increases as the polynomial order increases.\par

The sub-optimality of the \Bezier dual basis may deteriorate the finite element approximation. From Theorem~\ref{thm:approximation-of-bezier-formulation} and an Aubin-Nitsche duality argument, the expected convergence rates of the proposed method with the \Bezier dual basis are $1$ in the $H^2_*$ norm and $2$ in the $L^2$ norm. Although the poor approximation power is currently a flaw in the \Bezier dual basis, its local support and straightforward construction make it an appealing choice in practical use. Additionally, the authors have recently developed a technique that restores the optimal approximation of the \Bezier dual basis without appreciably changing the simplicity of construction. These results will be reported in a forthcoming paper.
\FloatBarrier

\subsection{The biharmonic problem on a two-patch domain}\label{sec:two_patch}

\begin{figure}[ht]
	\centering
	\begin{subfigure}[t]{0.3\textwidth}
		\includegraphics[width=\textwidth]{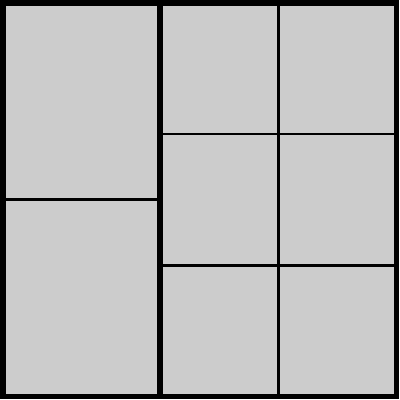}
		\caption{Simple non-conforming mesh}\label{fig:two_patch_biharmonic_problem_basic}
	\end{subfigure}
	\hfill
	\begin{subfigure}[t]{0.3\textwidth}
		\includegraphics[width=\textwidth]{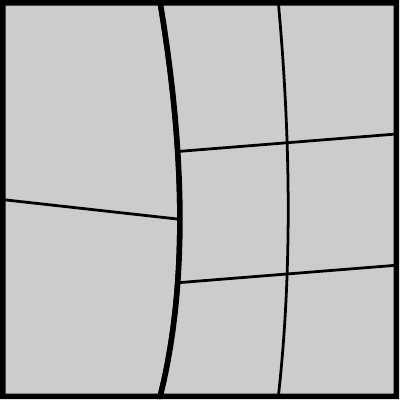}
		\caption{Distorted non-conforming mesh}\label{fig:two_patch_biharmonic_problem_distorted}
	\end{subfigure}
	\hfill
	\begin{subfigure}[t]{0.3\textwidth}
		\includegraphics[width=\textwidth]{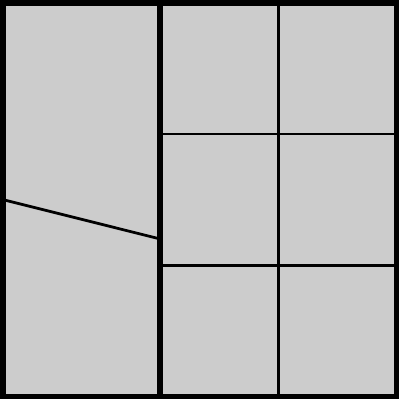}
		\caption{non-conforming mesh with mismatched parameterizations}\label{fig:two_patch_biharmonic_problem_nonmatch}
	\end{subfigure}\\
	\begin{subfigure}[t]{0.35\textwidth}
		\includegraphics[width=\textwidth]{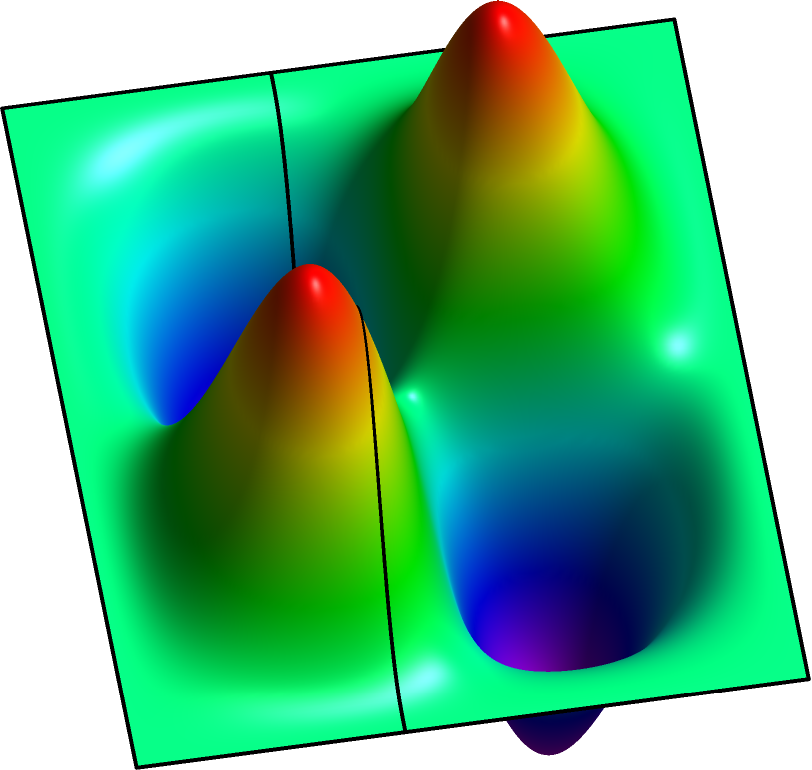}

		\caption{The manufactured solution\vspace{3mm}}\label{fig:two_patch_biharmonic_problem_solution-plot}
	\end{subfigure}
	\caption{The discretizations of the domain $\Omega$ ((a) - (c)) and the manufactured solution (d) with the property $u=\frac{\partial{u}}{\partial{\mathbf{n}}}=0$ on $\partial{\Omega}$ for the problem in Section~\ref{sec:two_patch}.}\label{fig:two_patch_biharmonic_problem}
\end{figure}

We now solve the biharmonic problem $\Delta^2{}u=f$ on a square domain $\Omega={(0,1)\times(0,1)}$. A manufactured solution is given by
\begin{equation}
	u(x,y)=\sin(2\pi{x})\sin(2\pi{y})(xy(x-1)(y-1))^2.
\end{equation}
This solution satisfies the homogeneous Dirichlet boundary condition ($u=\frac{\partial{u}}{\partial{\mathbf{n}}}=0$) and is shown in Figure~\ref{fig:two_patch_biharmonic_problem_solution-plot}. The domain $\Omega$ is decomposed into two patches $\Omega_1={(0,0.4)\times(0,1)}$ and $\Omega_2={(0.4,1)\times(0,1)}$, as shown in Figure~\ref{fig:two_patch_biharmonic_problem}. The right-hand side function $f$ can be obtained by applying the biharmonic operator to $u$.\par

\begin{figure}[ht]
	\centering
	\begin{subfigure}[b]{0.33\textwidth}
		\includegraphics[width=\textwidth]{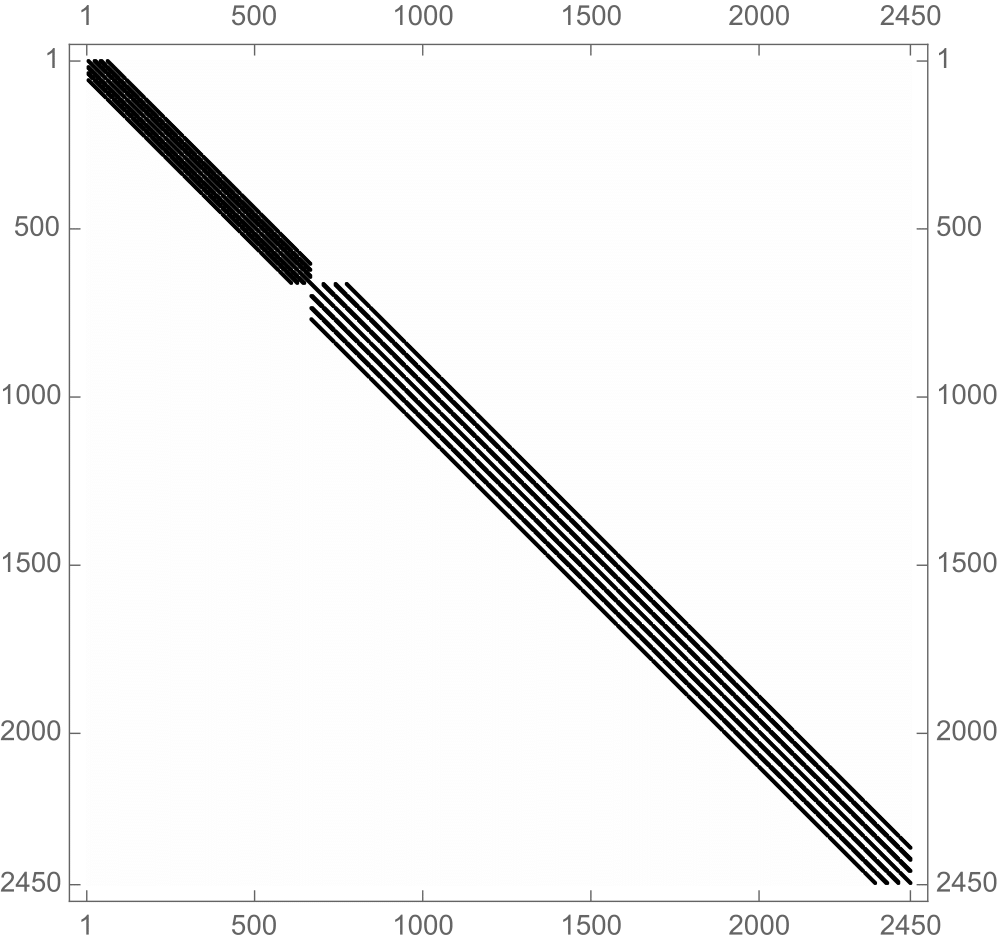}
		\caption{}
	\end{subfigure}
	\begin{subfigure}[b]{0.33\textwidth}
		\includegraphics[width=\textwidth]{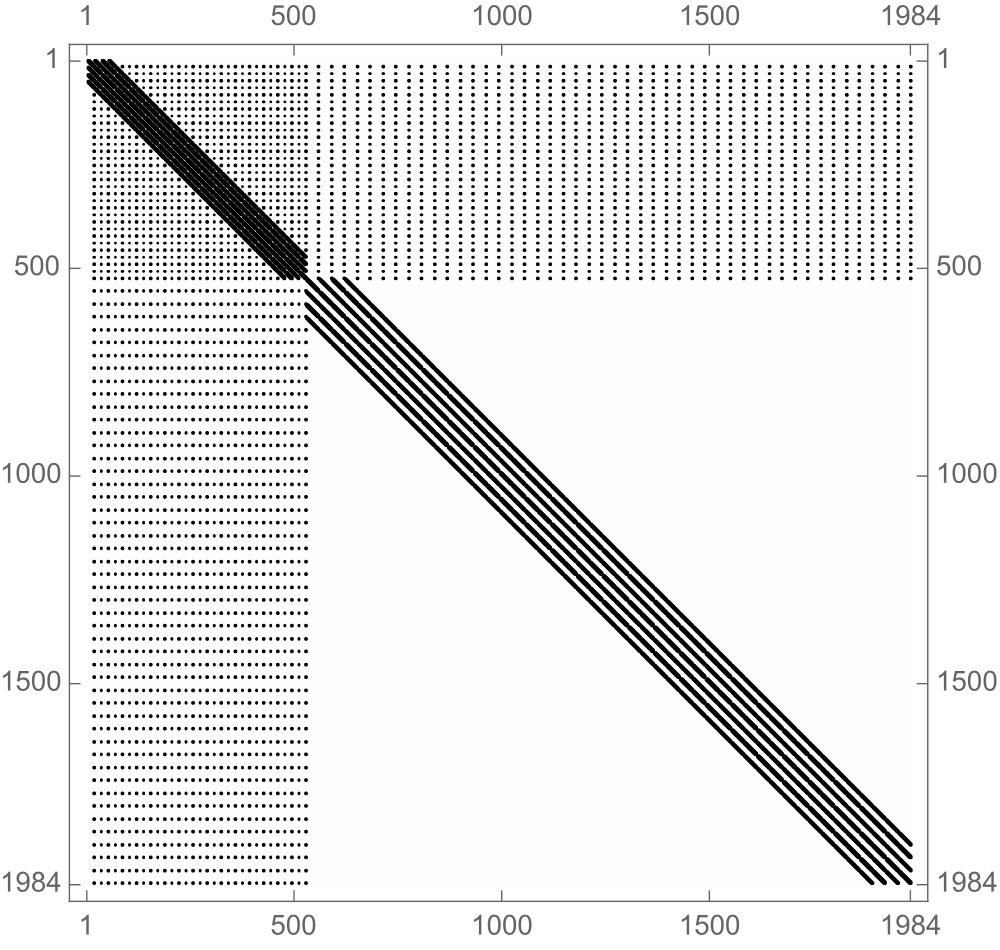}
		\caption{}
	\end{subfigure}
	\begin{subfigure}[b]{0.33\textwidth}
		\includegraphics[width=\textwidth]{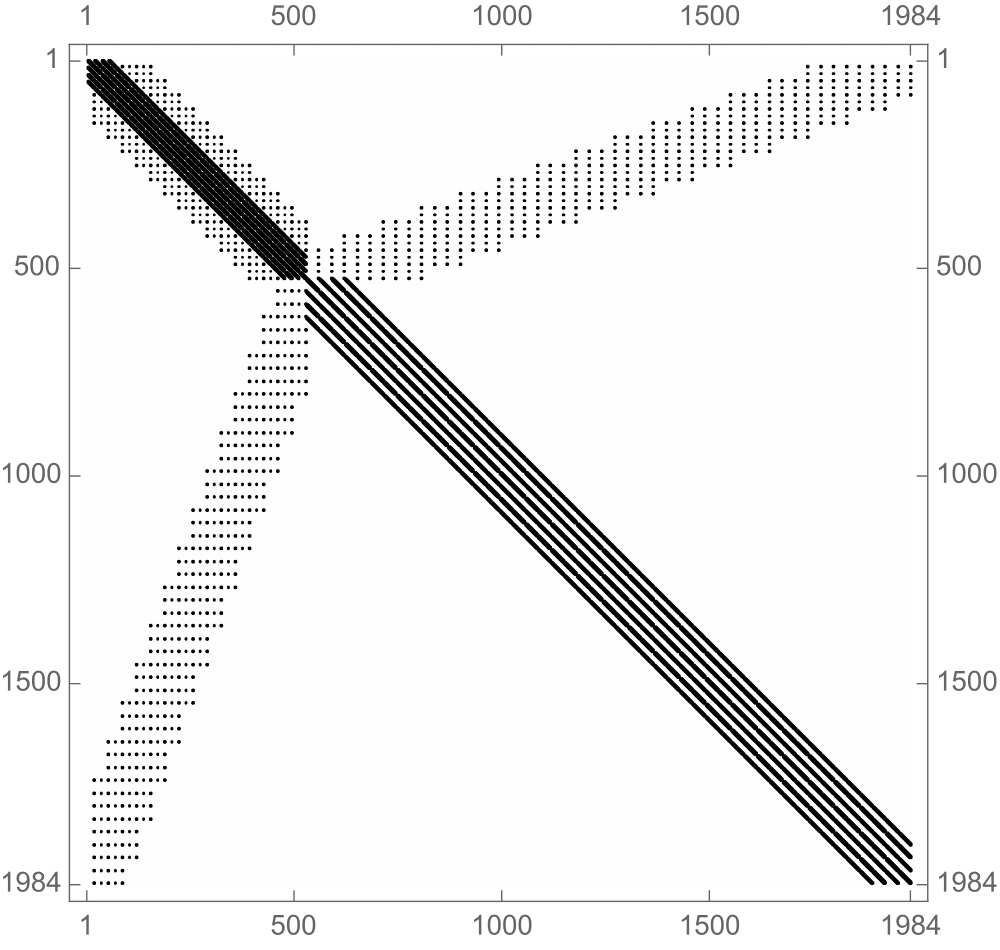}
		\caption{}
	\end{subfigure}
	\caption{Stiffness matrix sparsity patterns for (a) the uncoupled linear system, (b) the coupled linear system using the global dual basis, and (c) the coupled linear system using the \Bezier dual basis for the problem in Section~\ref{sec:two_patch}. The stiffness matrices are computed from the two-patch domain in Figure~\ref{fig:two_patch_biharmonic_problem_basic} after $4$ levels of refinement.}\label{fig:sparsity_pattern}
\end{figure}

The sparsity patterns for the stiffness matrices corresponding to the uncoupled problem, the coupled problem using the global dual basis, and the coupled problem using the \Bezier dual basis are shown in Figure~\ref{fig:sparsity_pattern}. Note that the matrix constructed using the global dual basis is denser than the matrix constructed using the \Bezier dual basis.\par

We conduct convergence studies for $p=2,3,4,5$ in both the $L^2$ and $H^2_*$ norms for the mesh shown in Figure~\ref{fig:two_patch_biharmonic_problem_basic}. The results are shown in Figure~\ref{fig:two_patc_biharmonic_convergence_basic}. Notice that despite the poor approximability of the \Bezier dual basis it performs surprisingly well in practice. As can be seen, both the global and B\'ezier dual basis obtain optimal convergence rates in both norms for all polynomial degrees. In fact, the convergence plots are almost identical between the global and \Bezier dual basis. The influence of the consistency error of the \Bezier dual basis cannot be observed for all tested polynomial degrees. We conjecture that for biharmonic problems, the coefficient $c_b$ in~\eqref{eq:consistency_error} is so small that the contribution of the consistency error in the finite element approximation is negligible.\par

\begin{figure}[ht]
	\centering
	\begin{subfigure}[b]{0.45\textwidth}
		\includestandalone[scale=.8]{two_patch_biharmonic_basic}
	\end{subfigure}
	\hfill
	\begin{subfigure}[b]{0.45\textwidth}
		\includestandalone[scale=.8]{two_patch_biharmonic_basic_H2}
	\end{subfigure}
	\caption{Convergence plots for the mesh shown in Figure~\ref{fig:two_patch_biharmonic_problem_basic}. Left: error measured in the $L^2$ norm. Right: error measured in the $H^2_*$ norm.}\label{fig:two_patc_biharmonic_convergence_basic}
\end{figure}

To study the performance of the proposed method in the presence of mesh distortion and mismatched parameterizations, we consider the meshes shown in Figure~\ref{fig:two_patch_biharmonic_problem_distorted} and~\ref{fig:two_patch_biharmonic_problem_nonmatch}. For the distorted mesh, shown in Figure~\ref{fig:two_patch_biharmonic_problem_distorted}, the proposed method with both global and B\'ezier dual basis functions perform similarly with optimal convergence rates being achieved in all cases as shown in Figure~\ref{fig:two_patc_biharmonic_convergence_distorted}.

For the mesh with mismatched parameterizations, shown in Figure~\ref{fig:two_patch_biharmonic_problem_nonmatch}, the convergence behavior of the B\'ezier dual basis, though optimal, deterioriates relative to the global dual basis as shown in Figure~\ref{fig:two_patc_biharmonic_convergence_nonmatch}. This indicates that the B\'ezier dual basis is more sensitive to mesh distortion than the global dual basis. Interestingly, as the mesh is refined, the results obtained using the $5^{th}$ order global dual basis become sub-optimal. We speculate that this is caused by an \textit{inf-sup} instability in this specific problem.

\begin{figure}[ht]
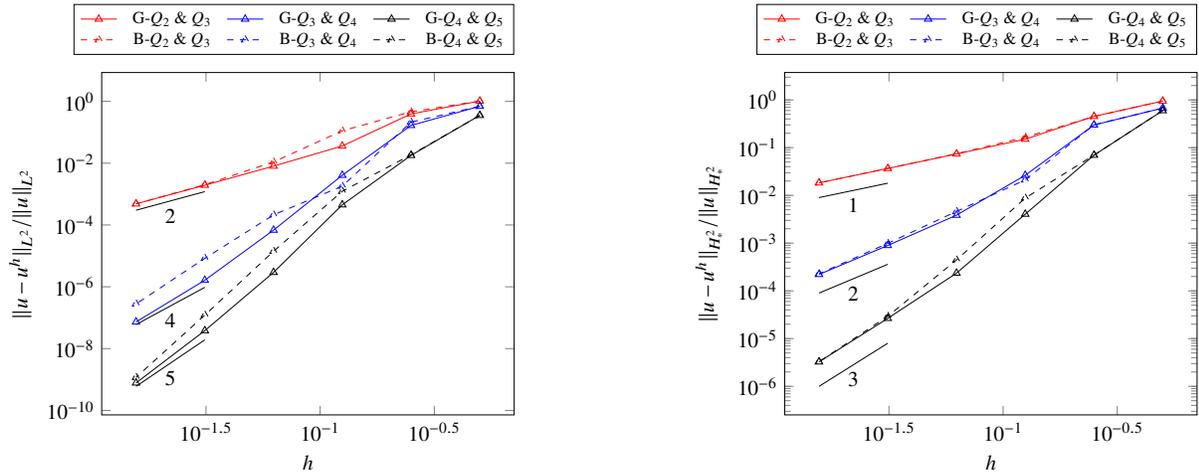

	\centering
	\begin{subfigure}[b]{0.45\textwidth}
		\includestandalone[scale=.8]{two_patch_biharmonic_diff_degree}
	\end{subfigure}
	\hfill
	\begin{subfigure}[b]{0.45\textwidth}
		\includestandalone[scale=.8]{two_patch_biharmonic_diff_degree_H2}
	\end{subfigure}
	\caption{Convergence plots for the mesh shown in Figure~\ref{fig:two_patch_biharmonic_problem_basic} with mismatched degrees. Left: error measured in the $L^2$ norm. Right: error measured in the $H^2_*$ norm.}\label{fig:two_patc_biharmonic_convergence_diff_degree}
\end{figure}

For the degree mismatched case shown in Figure~\ref{fig:two_patch_biharmonic_problem_basic}, the convergence rates are between $p_\text{left}+1$ and $p_\text{left}+2$ in the $L^2$ norm, and between $p_\text{left}-1$ and $p_\text{left}$ in the $H^2_*$ norm as expected for all cases and shown in Figure~\ref{fig:two_patc_biharmonic_convergence_diff_degree}. \par

\begin{figure}[ht]
	\centering
	\begin{subfigure}[b]{0.45\textwidth}
		\includestandalone[scale=.8]{two_patch_biharmonic_distorted}
	\end{subfigure}
	\hfill
	\begin{subfigure}[b]{0.45\textwidth}
		\includestandalone[scale=.8]{two_patch_biharmonic_distorted_H2}
	\end{subfigure}
	\caption{Convergence plots for the mesh shown in Figure~\ref{fig:two_patch_biharmonic_problem_distorted}. Left: error measured in the $L^2$ norm. Right: error measured in the $H^2_*$ norm.}\label{fig:two_patc_biharmonic_convergence_distorted}
\end{figure}

\begin{figure}[ht]
	\centering
	\begin{subfigure}[b]{0.45\textwidth}
		\includestandalone[scale=.8]{two_patch_biharmonic_nonmatch}
	\end{subfigure}
	\hfill
	\begin{subfigure}[b]{0.45\textwidth}
		\includestandalone[scale=.8]{two_patch_biharmonic_nonmatch_H2}
	\end{subfigure}
	\caption{Convergence plots for the mesh shown in Figure~\ref{fig:two_patch_biharmonic_problem_nonmatch}. Left: error measured in the $L^2$ norm. Right: error measured in the $H^2_*$ norm.}\label{fig:two_patc_biharmonic_convergence_nonmatch}
\end{figure}

Although a functional analysis of the contribution of the consistency error in the finite element approximation error is beyond the scope of this paper and postponed for future work, here we study the influence of the consistency error in a numerical manner. Since the finite element error is composed of the approximation error and the consistency error, the effect of the consistency error can be demonstrated by a comparison between the finite element error and the approximation error. The approximation error is the best $H^2_*$ approximation of $u$ in the discretized weak $C^1$ space $\mathcal{K}_b^h$, which is given as: find $u\in\mathcal{K}_b^h$ such that
\begin{equation}
	\langle{v^h,u^h}\rangle_{H^2_*}= \langle{v^h,u}\rangle_{H^2_*}\quad\quad\forall{v^h\in\mathcal{K}_b^h}.
\end{equation}

\begin{figure}[ht]
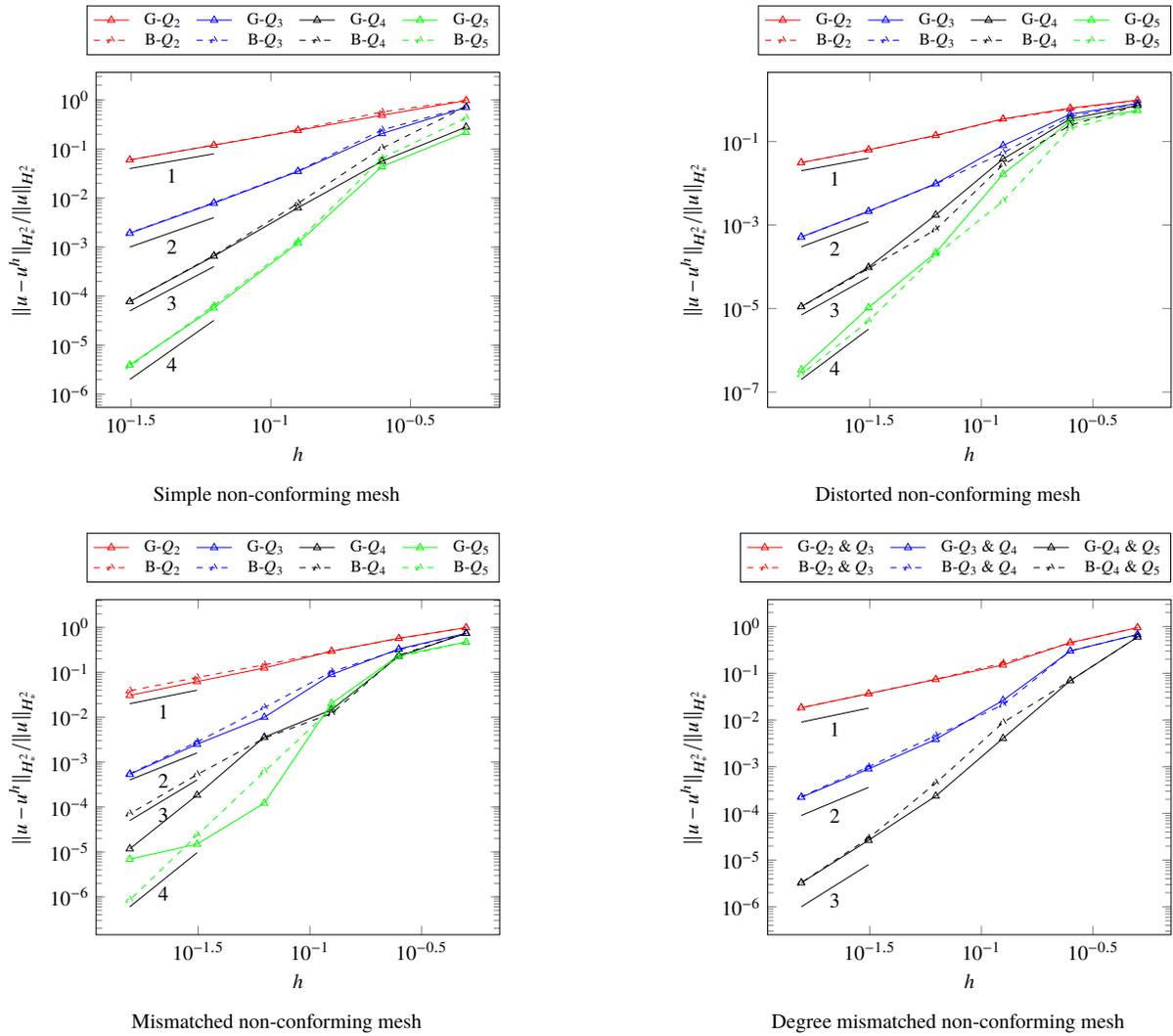

	\captionsetup[subfigure]{labelformat=empty, font = footnotesize, justification=centering}
	\centering
	\begin{subfigure}[b]{0.45\textwidth}
		\includestandalone[scale=.8]{two_patch_biharmonic_basic_H2}
		\caption{Simple non-conforming mesh}
		\vspace*{3mm}
	\end{subfigure}
	\hfill
	\begin{subfigure}[b]{0.45\textwidth}
		\includestandalone[scale=.8]{two_patch_biharmonic_distorted_H2}
		\caption{Distorted non-conforming mesh}
		\vspace*{3mm}
	\end{subfigure}
	\begin{subfigure}[b]{0.45\textwidth}
		\includestandalone[scale=.8]{two_patch_biharmonic_nonmatch_H2}
		\caption{Mismatched non-conforming mesh}
	\end{subfigure}
	\hfill
	\begin{subfigure}[b]{0.45\textwidth}
		\includestandalone[scale=.8]{two_patch_biharmonic_diff_degree_H2}
		\caption{Degree mismatched non-conforming mesh}
	\end{subfigure}
	\caption{Convergence plots of the approximation error for the two-patch coupling problem shown in Figure~\ref{fig:two_patch_biharmonic_problem} and described in Section~\ref{sec:two_patch}.}\label{fig:two_patch_best_approximation}
\end{figure}

Plots of the approximation error for the proposed method for the meshes shown in Figure~\ref{fig:two_patch_biharmonic_problem} are shown in Figure~\ref{fig:two_patch_best_approximation}. As can be seen, the convergence plots of the approximation error are identical to those of the finite element error in the $H^2_*$ norm. The approximation errors for all cases are no more than $1\%$ smaller than their finite element counterparts, which confirms our conjecture that the contribution of the best approximation error for the Lagrange multipliers (the consistency error) are negligible for the problems we tested. In addition, the approximation error plots for the global dual basis also demonstrate wavy and less asymptotic behavior. For the $p=5$ mismatched non-conforming mesh case it also suffers from reduced convergence rates, which confirms that the main cause of this phenomena in the finite element approximation is due to the \textit{inf-sup} instability.\par

\subsection{The biharmonic problem on multi-patch domains}
\subsubsection{The biharmonic problem on a three-patch domain}\label{sec:three-patch}
\begin{figure}[ht]
	\centering
	\begin{subfigure}[b]{0.33\textwidth}
		\includegraphics[width=\textwidth]{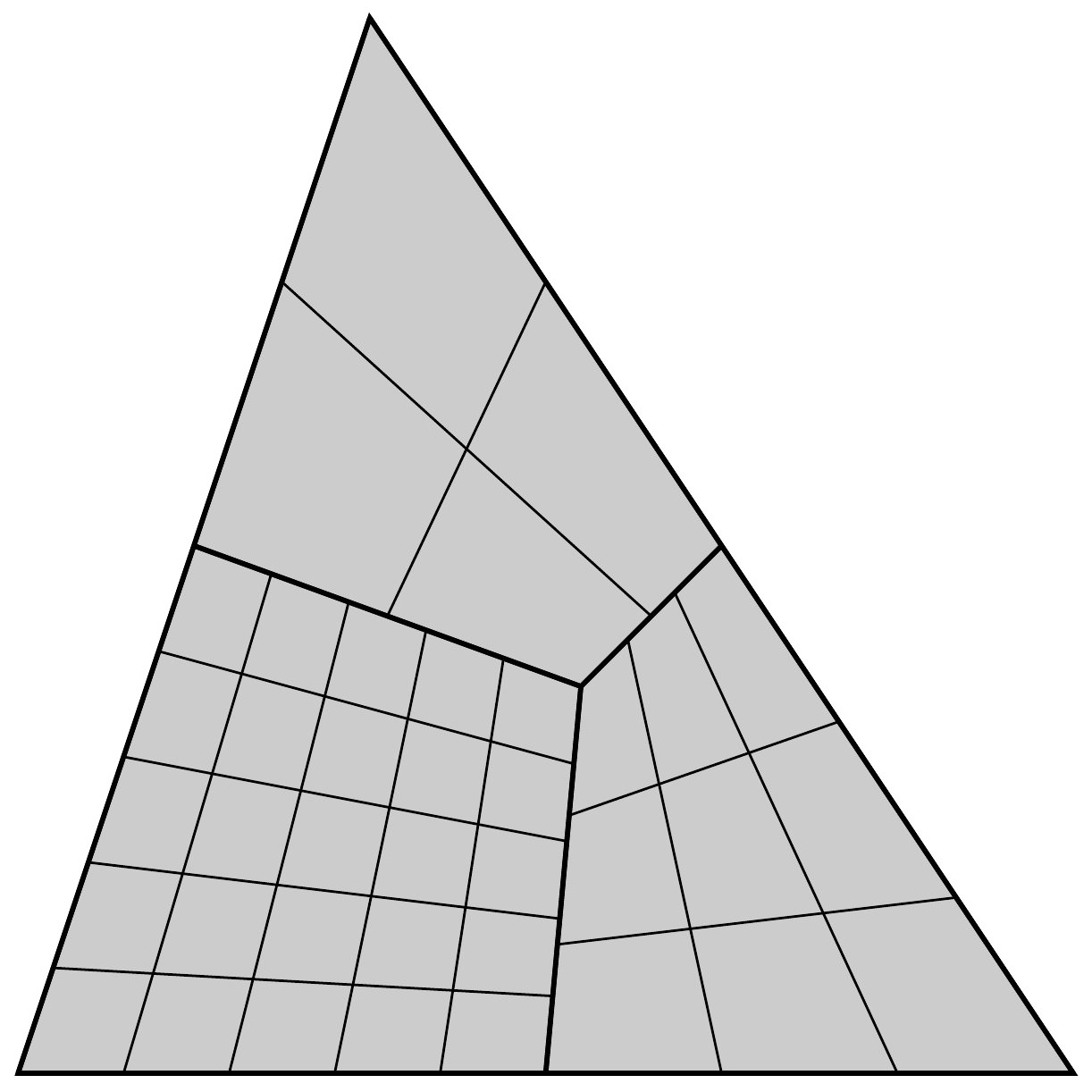}
		\caption{Non-conforming mesh}
	\end{subfigure}
	\begin{subfigure}[b]{0.41\textwidth}
		\includegraphics[width=\textwidth]{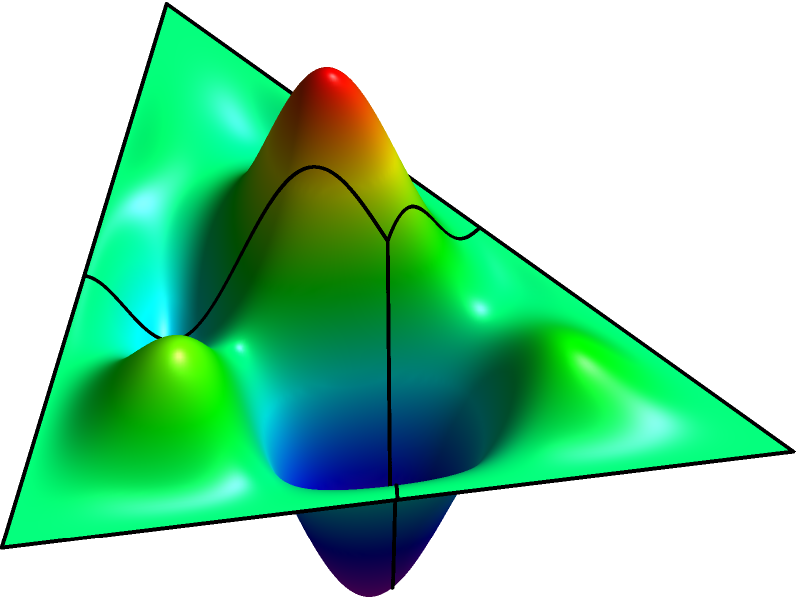}
		\caption{The manufactured solution}
	\end{subfigure}
	\caption{The three-patch domain and the manufactured solution for the problem in Section~\ref{sec:three-patch}.}\label{fig:three_patch_biharmonic_problem}
\end{figure}

We now examine the proposed method for multi-patch coupling. We first solve a biharmonic problem with the manufactured solution
\begin{equation}
	u(x,y)=sin(2\pi{x})sin(2\pi{y})\left({y(3x-y)(3x+2y-9)}\right)^2,
\end{equation}
on the triangular domain decomposed into three patches as shown in Figure~\ref{fig:three_patch_biharmonic_problem}. Both multi-patch treatments are tested in this problem.

\begin{figure}[ht]
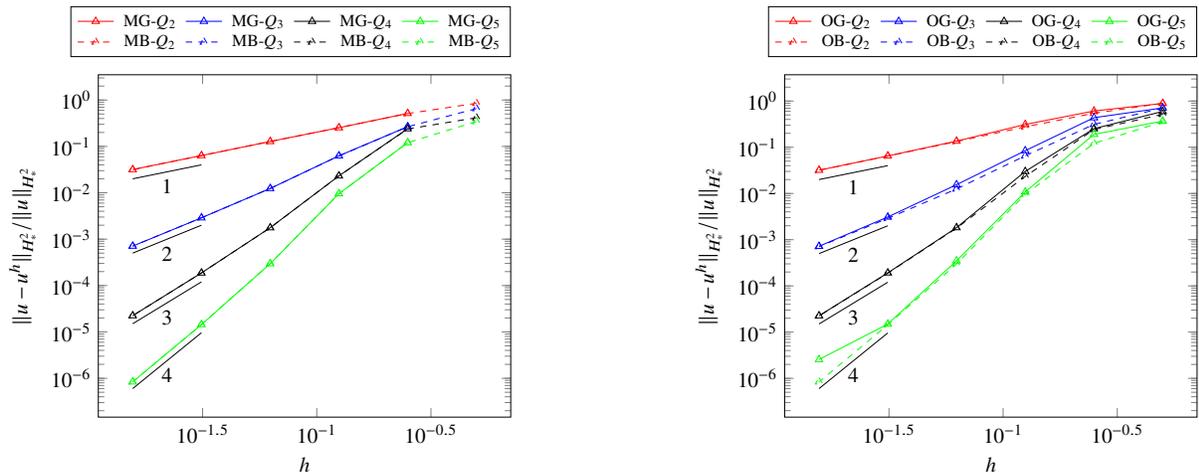

	\centering
	\begin{subfigure}[b]{0.45\textwidth}
		\includestandalone[scale=.8]{three_patch_modify_projection_basic_H2}
	\end{subfigure}
	\hfill
	\begin{subfigure}[b]{0.45\textwidth}
		\includestandalone[scale=.8]{three_patch_projection_basic_H2}
	\end{subfigure}
	\caption{Convergence plots of the approximation error for three-patch coupling in Section.~\ref{sec:three-patch}.}\label{fig:three_patch_approximation}
\end{figure}

The results are shown in Figure.~\ref{fig:three_patc_biharmonic_convergence}. As can be seen, the MG method yields optimal convergence rates in both measures for all tested polynomial orders. The MB method, on the other hand, yields $\mathcal{O}(h)$ convergence in $H^2_*$ norm and $\mathcal{O}(h^2)$ convergence in $L^2$ norm. The poor performance of the MB method is due to the consistency error, which can be verified by the optimal convergence rates in the approximation error (see Figure~\ref{fig:three_patch_approximation}). Moreover, for $p=3,4,5$, the finite element error increases as the polynomial order increases, which is consistent with the approximation power of the \Bezier dual basis (see Figure~\ref{fig:l_2_error}). Nevertheless, despite the sub-optimal convergence for certain circumstances, the MB method still converges asymptotically for all tested cases. The OG method yields optimal results for all tested polynomial orders except $p=5$. The sub-optimal rate for $p=5$ in both the $L^2$ and $H^2_*$ norms is due to the approximation error (see Figure~\ref{fig:three_patch_approximation}). Although the consistency error still influences the finite element approximation of the OB method, the error level at which the consistency error dominates is much lower ($10^{4}$ times lower in the $L^2$ norm and $10^{3}$ times lower in the $H^2_*$ norm) than that of the MB method. As a result, for $p=2,3,4$, the results obtained from the OB method demonstrate optimal convergence with sub-optimal convergence only occurring at the finest mesh for $p=5$. The approximation error for both the MB and OB methods are optimal for all tested polynomial orders, which indicates the \textit{inf-sup} stability of the proposed method with the \Bezier dual basis.

\begin{figure}[ht]
	\centering
	\begin{subfigure}[b]{0.45\textwidth}
		\includestandalone[scale=.8]{three_patch_modify_biharmonic_basic}
	\end{subfigure}
	\hfill
	\begin{subfigure}[b]{0.45\textwidth}
		\includestandalone[scale=.8]{three_patch_modify_biharmonic_basic_H2}
	\end{subfigure}
	\begin{subfigure}[b]{0.45\textwidth}
		\includestandalone[scale=.8]{three_patch_biharmonic_basic}
	\end{subfigure}
	\hfill
	\begin{subfigure}[b]{0.45\textwidth}
		\includestandalone[scale=.8]{three_patch_biharmonic_basic_H2}
	\end{subfigure}
	\caption{Convergence plots for the three-patch problem in Section~\ref{sec:three-patch}. Left: error measured in $L^2$ norm. Right: error measured in $H^2_*$ norm.}\label{fig:three_patc_biharmonic_convergence}
\end{figure}

\FloatBarrier

\subsubsection{The biharmonic problem on a five-patch domain}\label{sec:five-patch}

\begin{figure}[ht]
	\centering
	\begin{subfigure}[b]{0.33\textwidth}
		\includegraphics[width=\textwidth]{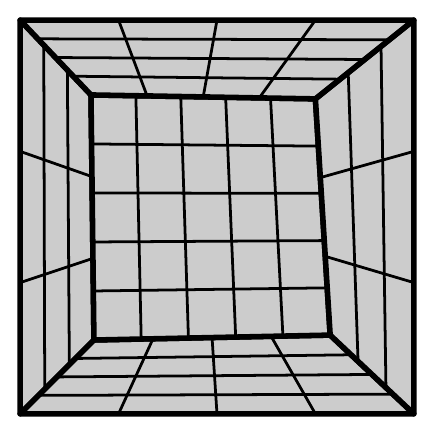}
		\caption{Non-conforming mesh}
	\end{subfigure}
	\begin{subfigure}[b]{0.34\textwidth}
		\includegraphics[width=\textwidth]{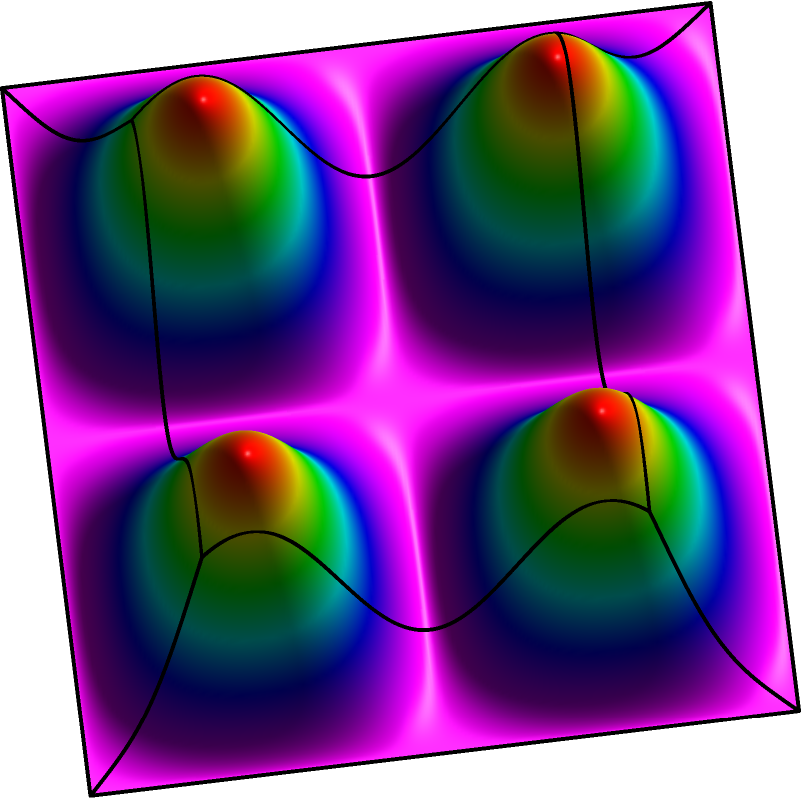}
		\caption{The manufactured solution}
	\end{subfigure}
	\caption{The five-patch domain parameterization and the manufactured solution in Section~\ref{sec:five-patch}.}\label{fig:five_patch_biharmonic_problem}
\end{figure}

To further study the effect of in-domain vertices, we solve a biharmonic problem on a five-patch domain, as shown in Figure~\ref{fig:five_patch_biharmonic_problem}, with the manufactured solution
\begin{equation}
	u(x,y)=\sin(2\pi{x})^2\sin(2\pi{y})^2.
\end{equation}

\begin{figure}[ht]
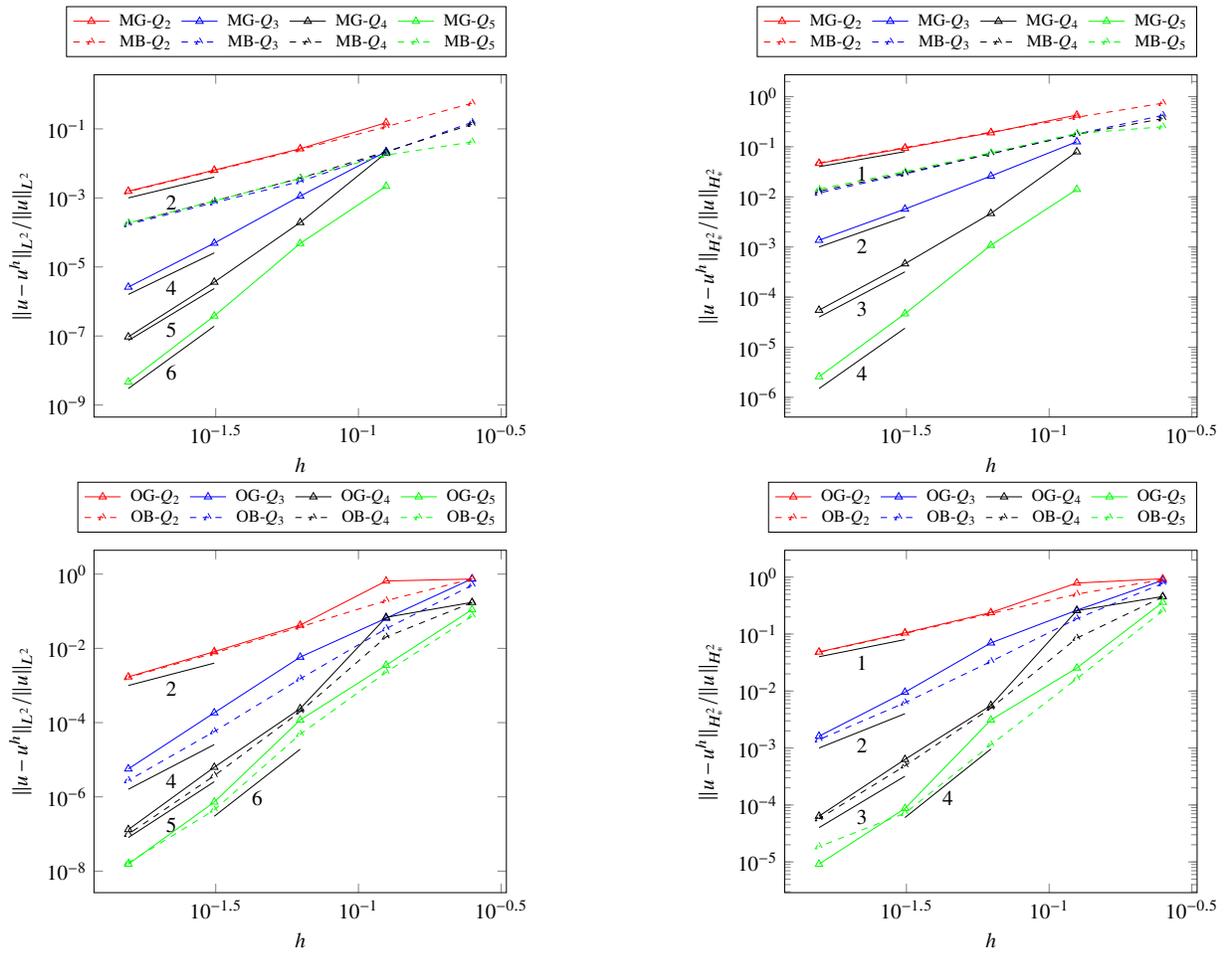

	\centering
	\begin{subfigure}[b]{0.45\textwidth}
		\includestandalone[scale=.8]{five_patch_modify_biharmonic_basic}
	\end{subfigure}
	\hfill
	\begin{subfigure}[b]{0.45\textwidth}
		\includestandalone[scale=.8]{five_patch_modify_biharmonic_basic_H2}
	\end{subfigure}

	\begin{subfigure}[b]{0.45\textwidth}
		\includestandalone[scale=.8]{five_patch_biharmonic_basic}
	\end{subfigure}
	\hfill
	\begin{subfigure}[b]{0.45\textwidth}
		\includestandalone[scale=.8]{five_patch_biharmonic_basic_H2}
	\end{subfigure}
	\caption{Convergence plots for the five-patch problem in Section~\ref{sec:three-patch}. Left: error measured in $L^2$ norm. Right: error measured in $H^2_*$ norm.}\label{fig:five_patc_biharmonic_convergence}
\end{figure}

\begin{figure}[ht]
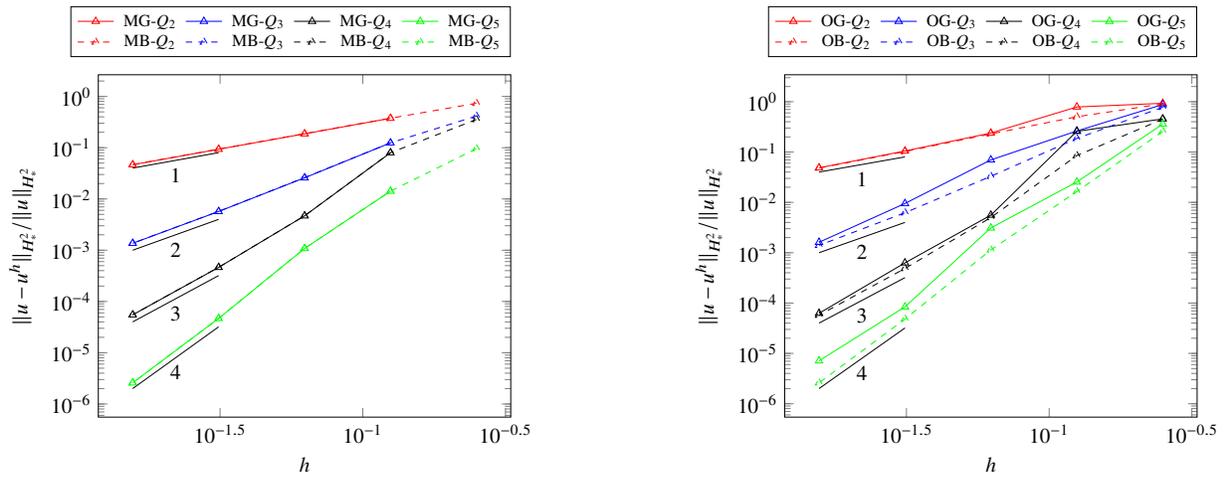

	\centering
	\begin{subfigure}[b]{0.45\textwidth}
		\includestandalone[scale=.8]{five_patch_modify_projection_basic_H2}
	\end{subfigure}
	\hfill
	\begin{subfigure}[b]{0.45\textwidth}
		\includestandalone[scale=.8]{five_patch_projection_basic_H2}
	\end{subfigure}
	\caption{Convergence plots of the approximation error for five-patch coupling in Section.~\ref{sec:five-patch}.}\label{fig:five_patch_approximation}
\end{figure}

The convergence behavior for all methods is shown in Figure~\ref{fig:five_patc_biharmonic_convergence}. The approximation error plots are shown in Figure~\ref{fig:five_patch_approximation}. The results are similar to that of the three-patch case. For the global dual basis, the MG method gives optimal convergence rates for both the $L^2$ and $H^2$ norms, while the OG method has an \textit{inf-sup} instability, which disrupts convergence rates for fine meshes. The OB method postpones the domination of the consistency error (i.e., postponed from $10^{-3}$ and $10^{-1}$ to $10^{-7}$ and $10^{-4}$ in the $L^2$ and $H^2$ norm, respectively) and significantly improves convergence rates.

\FloatBarrier

\subsection{Transverse vibrations of an elastic membrane}\label{sec:eigenvalue}

\begin{figure}[ht]
	\centering
	\includegraphics[width=.3\textwidth]{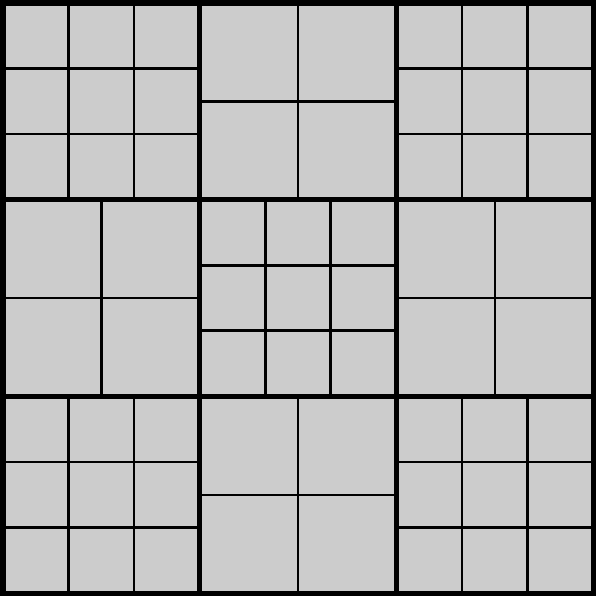}
	\caption{The nine-patch domain parameterization for the eigenvalue problem in Section~\ref{sec:eigenvalue}.}\label{fig:eigenvalue_mesh}
\end{figure}

We now study the OB method in the context of a $2^\text{nd}$-order eigenvalue problem. Of particular interest is the behavior of the highest frequencies in the system since they govern the critical timestep size in an explicit time-stepping scheme. In particular, we consider the transverse vibration of a square, elastic membrane on the domain $\left[{0,3}\right] \times \left[{0,3}\right]$ with the nonconforming discretization shown in Figure~\ref{fig:eigenvalue_mesh}. The natural frequencies and modes are governed by
\begin{equation}
	\left\{\begin{split}
		\nabla^2 u(x,y) +\omega^2 u(x,y) = 0, & \text{ in } \Omega,\\
		u(x, y) = 0 & \text{ on } \partial \Omega,
	\end{split}\right.
\end{equation}
where $\omega$ is the natural frequency. The exact natural frequencies are

\begin{equation}
	\omega_{mn} = \pi \sqrt{\left(\frac{m}{L}\right)^2+\left(\frac{n}{L}\right)^2},\quad\quad m,n = 1,2,3,\dots,
\end{equation}
where $L$ is the length of the boundary.\par

\begin{figure}[ht]
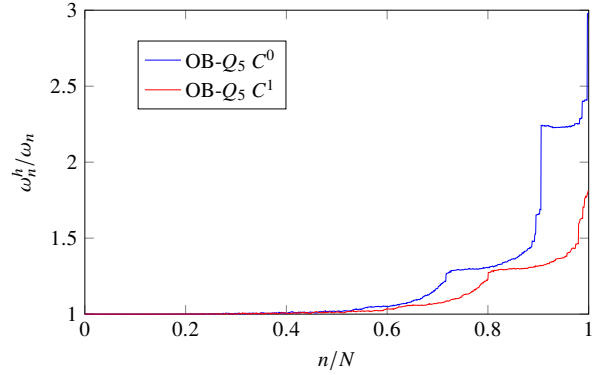

	\centering
	\includestandalone[scale=.8]{spectrum}
	\caption{Normalized discrete spectra using $p=5$. Results are obtained after three uniform refinements.}\label{fig:spectra}
\end{figure}

The highest computed eigenvalues for the OB method are given in Table~\ref{tab:eigenvalues}. As can be seen, the weak $C^1$ coupling dramatically reduces the highest eigenvalues for all tested cases. The effect becomes more significant as the polynomal degree is increased. For $p=5$, the highest eigenvalue obtained through weak $C^1$ coupling is $\frac{3}{5}$ of that obtained through weak $C^0$ coupling. The normalized discrete spectra for $p=5$ are shown in Figure~\ref{fig:spectra}. In this case, weak $C^1$ coupling improves the behavior of the entire spectra.

\begin{table}[!htbp]
	\centering
	\caption{The highest eigenvalues obtained by solving the eigenvalue problem for the square domain (see Figure~\ref{fig:eigenvalue_mesh})}
	\begin{tabular*}{\textwidth}{c @{\extracolsep{\fill}} SSSSSSSS}
		\toprule
		Refine &  \multicolumn{2}{c}{OB-$Q_2$} & \multicolumn{2}{c}{OB-$Q_3$} & \multicolumn{2}{c}{OB-$Q_4$} & \multicolumn{2}{c}{OB-$Q_5$}\\\cmidrule(lr{.5em}){2-3}\cmidrule(lr{.5em}){4-5}\cmidrule(lr{.5em}){6-7}\cmidrule(lr{.5em}){8-9}
		{} & \multicolumn{1}{c}{$C^0$} & \multicolumn{1}{c}{$C^1$} & \multicolumn{1}{c}{$C^0$} & \multicolumn{1}{c}{$C^1$} & \multicolumn{1}{c}{$C^0$} & \multicolumn{1}{c}{$C^1$} & \multicolumn{1}{c}{$C^0$} & \multicolumn{1}{c}{$C^1$}\\
		\midrule
		0   &  16.12  & 11.47  & 24.76  & 16.36 & 35.38  & 22.56  & 48.48  & 29.43\\
		1   &  20.97  & 17.41  & 29.39  & 20.86 & 41.41  & 26.91  & 54.60  & 34.18\\
		2   &  31.30  & 28.04  & 43.60  & 30.39 & 60.49  & 37.35  & 79.57  & 46.82\\
		3   &  54.04  & 52.61  & 78.60  & 51.09 & 108.47 & 61.75  & \cellcolor{hred!30}141.44 & \cellcolor{hred!30}81.53\\
		\bottomrule
	\end{tabular*}
	\label{tab:eigenvalues}
\end{table}

\FloatBarrier

\section{Conclusion}\label{sec:conclusion}

In this paper, we present a dual mortar formulation for the biharmonic problem and investigate its properties analytically and numerically. With the help of the dual mortar suitable $C^1$ constraint, the biorthogonality between the dual basis functions and the corresponding primal spline basis functions can be extended to the discretized $C^1$ constraint matrix. Hence, the condensed stiffness matrix can be formed efficiently without the need to solve linear systems associated with each intersection. Furthermore, the condensed stifness matrix remains sparse if the dual basis functions are compactly supported, which is the case for the \Bezier dual basis.\par

Due to the presence of in-domain vertices, some control points may serve as both slave and master. To overcome this we propose two solutions. The first method localizes the constraints to the neighborhood of each vertex and solves for the null space of a localized linear system. The second method reduces the number of constraints at each vertex by reducing the number of degrees of freedom of the dual basis. For both cases, when the \Bezier dual basis is used, the resulting linear systems are sparse. A suite of numerical experiments demonstrate the effectiveness of the proposed approach. \par

Areas of future research include the extension of the approach to Kirchhoff-Love shells and the development of compactly supported dual basis functions with optimal approximation power.

\clearpage
\bibliography{mybibfile}

\end{document}

%% file: tikz_scripts/geometry.tex
\begin{tikzpicture}[use Hobby shortcut,scale=1.2]
	\filldraw[fill=gray!40!white, draw=black, ultra thick] (1,0) -- (4,0)
    -- (4, 4)
    -- (0, 4)
	-- (0,1)
    arc (90:0:1);
\end{tikzpicture}

%% file: tikz_scripts/trimmed_geometry.tex
\begin{tikzpicture}[use Hobby shortcut,scale=1.2]
	\filldraw[fill=gray!40!white, draw=black, ultra thick] (1,0) -- (4,0)
    -- (4, 4)
    -- (0, 4)
	-- (0,1)
    arc (90:0:1);
    \draw[name path=curve1,draw=black, ultra thick] (0,4) -- (0,1) arc (90:0:1) -- (4,0);
    \foreach \i in {1,...,10}
    {
        \path [name path=line\i](.4*\i,4)--(.4*\i,0);
        \draw [name intersections={of=curve1 and line\i}, black, thick] (.4* \i,4) -- (intersection-1);
    }
    \foreach \i in {1,...,10}
    {
        \path [name path=line\i](4,.4*\i)--(0,.4*\i);
        \draw [name intersections={of=curve1 and line\i}, black, thick] (4,.4*\i) -- (intersection-1);
    }
\end{tikzpicture}

%% file: tikz_scripts/reparameterized_geometry.tex
\begin{tikzpicture}[use Hobby shortcut,scale=1.2]
	\filldraw[fill=gray!40!white, draw=black, ultra thick] (1,0) -- (4,0)
    -- (4, 4)
    -- (0, 4)
	-- (0,1)
    arc (90:0:1);
    
    \foreach \f in {1,2,...,10}{\draw[thick] plot[draw=black] file {tikz_scripts/data/curv_u_\f.txt};}
    \foreach \f in {1,2,...,10}{\draw[thick] plot[draw=black] file {tikz_scripts/data/curv_v_\f.txt};}
\end{tikzpicture}

%% file: tikz_scripts/basisfunctions.tex
\begin{tikzpicture}
\begin{axis}[axis lines=none,xmin=-0.03,xmax=1.03,ymin=-0.1,ymax=1, scale only axis,
        width=0.6\textwidth,
        height=0.2\textwidth]

\addplot+[no markers,ultra thick,o1, opacity=0.2] table [y=N0, x=x]{basis.dat};
\addplot+[no markers,ultra thick,restrict x to domain=0:.2,gray, opacity=0.2] table [y=N1, x=x]{basis.dat};
\addplot+[no markers,ultra thick,restrict x to domain=.2:.4,gray, opacity=0.2] table [y=N1, x=x]{basis.dat};
\addplot+[no markers,ultra thick,restrict x to domain=0:.2,magenta, opacity=0.2] table [y=N2, x=x]{basis.dat};
\addplot+[no markers,ultra thick,restrict x to domain=.2:.6,magenta, opacity=0.2] table [y=N2, x=x]{basis.dat};
\addplot+[no markers,ultra thick,b1,solid] table [y=N3, x=x]{basis.dat};
\addplot+[no markers,ultra thick,restrict x to domain=.4:.8,red, opacity=0.2,solid] table [y=N4, x=x]{basis.dat};
\addplot+[no markers,ultra thick,restrict x to domain=.8:1,solid,red, opacity=0.2] table [y=N4, x=x]{basis.dat};
\addplot+[no markers,ultra thick,restrict x to domain=.6:.8,solid,g1, opacity=0.2] table [y=N5, x=x]{basis.dat};
\addplot+[no markers,ultra thick,restrict x to domain=.8:1,solid,g1, opacity=0.2] table [y=N5, x=x]{basis.dat};
\addplot+[no markers,ultra thick,restrict x to domain=.8:.9999,y1, opacity=0.2] table [y=N6, x=x]{basis.dat};
\addplot+[mark=*,mark options={fill=white},mark size=2,color=black,ultra thick] plot coordinates {
        (0,    0)
        (1.0/5,   0)
        (2.0/5,   0)
        (3.0/5,    0)
        (4.0/5,    0)
        (1,    0)};

\end{axis}
\end{tikzpicture}

%% file: tikz_scripts/dual_global_one.tex
\begin{tikzpicture}
\begin{axis}[axis lines=none,xmin=-0.03,xmax=1.03,ymin=-0.777778,ymax=2, scale only axis,
	width=0.6\textwidth,
        height=0.2\textwidth]

\addplot+[no markers, smooth,ultra thick, solid, b1] table [y=N4, x=x]{dualbasis_global.dat};

\addplot+[mark=*, solid, mark options={fill=white},mark size=2,color=black,ultra thick] plot coordinates {
			        (0,    0)
			        (1.0/5,   0)
			        (2.0/5,   0)
			        (3.0/5,    0)
			        (4.0/5,    0)
			        (1,    0)};

\end{axis}
\end{tikzpicture}

%% file: tikz_scripts/bezier_dual_one.tex
	\begin{tikzpicture}
		\begin{axis}[axis lines=none,xmin=-0.03,xmax=1.03,ymin=-6,ymax=5, scale only axis,
			width=0.6\textwidth,
        height=0.2\textwidth]

		\addplot+[no markers, smooth,ultra thick, solid, b1,select coords between index={20}{39}] table [y=N4, x=x]{dualbasis.txt};
		\addplot+[no markers, smooth,ultra thick, solid, b1,select coords between index={40}{59}] table [y=N4, x=x]{dualbasis.txt};
		\addplot+[no markers, smooth,ultra thick, solid, b1,select coords between index={60}{79}] table [y=N4, x=x]{dualbasis.txt};

		\addplot+[mark=*, solid, mark options={fill=white},mark size=2,color=black,ultra thick] plot coordinates {
							(0,    0)
							(1.0/5,   0)
							(2.0/5,   0)
							(3.0/5,    0)
							(4.0/5,    0)
							(1,    0)};

		\end{axis}
	\end{tikzpicture}

%% file: tikz_scripts/patch_partition.tex
\begin{tikzpicture}[scale=1.5]
  \begin{scope}[shift={(-.2,-.2)}]
    \draw[->] (0,0) -- (.5,0) node[right]{$x$};
    \draw[->] (0,0) -- (0,.5) node[left]{$y$};
  \end{scope}

    \coordinate(Origin) at ($ (2,{2/3*sqrt(3)}) $);

    \draw[thick] (0,0) -- (4,0) coordinate[midway](M1)-- ($ (2,{2*sqrt(3)}) $)
      coordinate[midway](M2)-- cycle  coordinate[midway](M3);
    \draw [purples1] (Origin) -- (M1) node [midway, left, black] {$\Gamma_{12}$}  node [midway, right, black] {$\Gamma_{21}$};
    \draw [purples1] (Origin) -- (M2) node [pos=.45, above=.1, black] {$\Gamma_{23}$}  node [pos=.55, below=.1, black] {$\Gamma_{32}$};
    \draw [purples1] (Origin) -- (M3) node [pos=.45, above=.1, black] {$\Gamma_{31}$}  node [pos=.55, below=.1, black] {$\Gamma_{13}$};

    \node [below right] at (Origin) {$v_1$};
    \node[mark size=3pt,color=purples4] at (Origin) {\pgfuseplotmark{*}};
  
    
    \node (Omega1) at (1.2,.8) {$\Omega_1$};
    \node (Omega2)  at (2.8,.8) {$\Omega_2$};
    \node at (2,2) {$\Omega_3$};
    \node at (1,2.2) {$\partial\Omega$};

  \begin{scope}[shift={(-4,2)},rotate=10]
    \draw (0,0) rectangle (1.6,1.6);
    \draw[purples1,ultra thick] (1.6,0) -- (1.6,1.6);
    \begin{scope}[shift={(.2,.2)}]
      \draw[->] (0,0) -- (.5,0) node[above]{$\xi^{(1)}$};
      \draw[->] (0,0) -- (0,.5) node[above]{$\eta{(1)}$};
    \end{scope}
    \node (hOmega1) at (.8,.8) {$\hat{\Omega}_1$} ;
    \node (hedge1) at (1.6,.8) {};
  \end{scope}
  \draw[->,>=stealth,thick] (hOmega1) to [out=30,in=150] node[pos = .5, right=.4]{$\mathbf{F}_{1}$} (Omega1);

  \begin{scope}[shift={(-2.5,0)},rotate=-10]
    \draw (0,0) rectangle (1.6,1.6);
    \draw[purples1,ultra thick] (0,0) -- (0,1.6);
    \begin{scope}[shift={(.2,.2)}]
      \draw[->] (0,0) -- (.5,0) node[above]{$\xi^{(2)}$};
      \draw[->] (0,0) -- (0,.5) node[above right=-.15]{$\eta{(2)}$};
    \end{scope}
    \node (hOmega2) at (.8,.8) {$\hat{\Omega}_2$} ;
    \node (hedge2) at (0,.8) {};
  \end{scope}
  \draw[->,>=stealth,thick] (hOmega2) to [out=-55,in=-125] node[pos = .5, below]{$\mathbf{F}_{2}$} (Omega2);

  \draw[->,>=stealth,thick] (hedge1) to [out=-50,in=120] node[pos = .5, right]{$\mathbf{E}_{12}$} (hedge2);  
\end{tikzpicture}

%% file: tikz_scripts/l_2_error.tex
\begin{tikzpicture}
    \begin{loglogaxis}[
        legend columns=3,
    	legend style={at={(1,1.2)}, nodes={scale=1, transform shape}, column sep=.4cm},
        xlabel=$h$,
        ylabel=${\|u-u^{h}\|_{L^2}}/{\|u\|_{L^2}}$ 
    ]

    \addplot [color=red,mark=triangle] plot coordinates {

        (1, .97)
        (.5, .967)    
        (.25,   .808)
        (.125,   .114)
        (.0625,   .024)
        (0.03125,  .00584)
    };

    \addplot [color=blue,mark=triangle] plot coordinates {

        (1, 1.09)
        (.5, .775)    
        (.25,   .387)
        (.125,   .0347)
        (.0625,   .00309)
        (0.03125,  .000353)
    };

    \addplot [color=black,mark=triangle] plot coordinates {

        (1, 1.09)
        (.5, .843)    
        (.25,   0.342314)
        (.125,   0.0118696)
        (.0625,   0.000440564)
        (0.03125,  2.3112e-05)
    };

    \addplot [color=red,mark=triangle, dashed] plot coordinates {

        (1, .97)
        (.5, .89)    
        (.25,   .94)
        (.125,   .50)
        (.0625,   .35)
        (0.03125,  .21)
    };

    \addplot [color=blue,mark=triangle, dashed] plot coordinates {

        (1, 1.090343)
        (.5, 0.742515)    
        (.25,   0.825903)
        (.125,   0.838765)
        (.0625,   0.614495)
        (0.03125,  0.337256)
    };

    \addplot [color=black,mark=triangle, dashed] plot coordinates {

        (1, 1.08983)
        (.5, 0.956762)    
        (.25,   0.897757)
        (.125,   0.913494)
        (.0625,   0.782978)
        (0.03125,  0.443544)
    };

    \addplot [black] plot coordinates {

        (.0625,    0.3)
        (.03125,   .17)
    } node[pos=0.5,anchor=north]{$1$};

    \addplot [black] plot coordinates {

        (.0625,    0.016)
        (.03125,   .004)
    } node[pos=0.5,anchor=north]{$2$};

    \addplot [black] plot coordinates {

        (.0625,    0.002)
        (.03125,   .00025)
    } node[pos=0.5,anchor=north]{$3$};

    \addplot [black] plot coordinates {

        (.0625,    2.2400e-04)
        (.03125,   1.4e-5)
    } node[pos=0.5,anchor=north]{$4$};
    \legend{G-$Q_1$\\G-$Q_2$\\G-$Q_3$\\B-$Q_1$\\B-$Q_2$\\B-$Q_3$\\}
    \end{loglogaxis}
\end{tikzpicture}

%% file: tikz_scripts/Q3_p=0.tex
\begin{tikzpicture}
        \begin{axis}[xlabel=$x$, ylabel=$y$, legend style={at={(1,.2)}},
        width=8cm,
        height=7cm,
        yticklabels={.9,1,1.1},
        ymin = .8,
        ymax = 1.2]
            \addplot+[mark=none, color = hblue, ultra thick] table [y=Y2, x=X]{Q3_p=0.txt}; 
            \addplot+[restrict x to domain=-1:-.01, mark=none, color = hred, ultra thick] table [y=Y1, x=X]{Q3_p=0.txt};
            \addplot+[restrict x to domain=0:1, mark=none, color = hred, ultra thick] table [y=Y1, x=X]{Q3_p=0.txt};
        \end{axis}
    \end{tikzpicture}

%% file: tikz_scripts/Q3_p=1.tex
\begin{tikzpicture}
        \begin{axis}[xlabel=$x$, ylabel=$y$, legend style={at={(1,.2)}},
            width=8cm,
            height=7cm]
            \addplot+[mark=none, color = hblue, ultra thick] table [y=Y2, x=X]{Q3_p=1.txt}; 
            \addplot+[restrict x to domain=-1:-.01, mark=none, color = hred, ultra thick] table [y=Y1, x=X]{Q3_p=1.txt};
            \addplot+[restrict x to domain=0:1, mark=none, color = hred, ultra thick] table [y=Y1, x=X]{Q3_p=1.txt};
        \end{axis}
    \end{tikzpicture}

%% file: tikz_scripts/Q3_p=2.tex
\begin{tikzpicture}
        \begin{axis}[xlabel=$x$, ylabel=$y$, legend style={at={(1,.2)}},
            width=8cm,
            height=7cm]
            \addplot+[mark=none, color = hblue, ultra thick] table [y=Y2, x=X]{Q3_p=2.txt}; 
            \addplot+[restrict x to domain=-1:-.01, mark=none, color = hred, ultra thick] table [y=Y1, x=X]{Q3_p=2.txt};
            \addplot+[restrict x to domain=0:1, mark=none, color = hred, ultra thick] table [y=Y1, x=X]{Q3_p=2.txt};
        \end{axis}
    \end{tikzpicture}

%% file: tikz_scripts/Q3_p=3.tex
\begin{tikzpicture}
        \begin{axis}[xlabel=$x$, ylabel=$y$, legend style={at={(1,.2)}},
        width=8cm,
        height=7cm]
            \addplot+[mark=none, color = hblue, ultra thick] table [y=Y2, x=X]{Q3_p=3.txt}; 
            \addplot+[restrict x to domain=-1:-.01, mark=none, color = hred, ultra thick] table [y=Y1, x=X]{Q3_p=3.txt};
            \addplot+[restrict x to domain=0:1, mark=none, color = hred, ultra thick] table [y=Y1, x=X]{Q3_p=3.txt};
        \end{axis}
    \end{tikzpicture}

%% file: tikz_scripts/two_patch_biharmonic_basic.tex
\begin{tikzpicture}
    \begin{loglogaxis}[
        legend columns=4,
    	legend style={at={(1,1.2)}, nodes={scale=0.8, transform shape}, column sep=.2cm},
        xlabel=$h$,
        ylabel=${\|u-u^{h}\|_{L^2}}/{\|u\|_{L^2}}$ 
    ]

    \addplot [color=red,mark=triangle] plot coordinates {

        (.5,    1.05039)
        (.25,   0.317211 )
        (.125,   0.0645939)
        (.0625,   0.0159494)
        (0.03125,  0.0039786)
    };

    \addplot [color=blue,mark=triangle] plot coordinates {

        (.5,    0.532167 )
        (.25,   0.0630832 )
        (.125,   0.00163679 )
        (.0625,   7.62415e-05)
        (0.03125,   4.40098e-06 )
    };

    \addplot [color=black,mark=triangle] plot coordinates {

        (.5,    0.0907098)
        (.25,   0.00947519)
        (.125,   0.000280315   )
        (.0625,   5.2263e-06   )
        (0.03125,  1.38275e-07 )
    };

    \addplot [color=green,mark=triangle] plot coordinates {

        (.5,    0.0761197)
        (.25,   0.00983706)
        (.125,   5.79501e-05)
        (.0625,   7.27108e-07)
        (0.03125,    9.66976e-09 )
    };

    \addplot [color=red,mark=triangle, dashed] plot coordinates {

        (.5,    1.10986 )
        (.25,    0.546794)
        (.125,   0.0971345)
        (.0625,   0.0171764)
        (0.03125,   0.00401722)
    };

    \addplot [color=blue,mark=triangle, dashed] plot coordinates {

        (.5,    0.545609)
        (.25,   0.101467 )
        (.125,   0.0021806)
        (.0625,   0.000105985)
        (0.03125,   5.74047e-06 )
    };

    \addplot [color=black,mark=triangle, dashed] plot coordinates {

        (.5,    0.510789)
        (.25,   0.0334686)
        (.125,   0.000515799)
        (.0625,   6.78752e-06)
        (0.03125,  1.47147e-07)
    };

    \addplot [color=green,mark=triangle, dashed] plot coordinates {

        (.5,    0.227382 )
        (.25,   0.0169498   )
        (.125,   6.32472e-05  )
        (.0625,   7.13693e-07)
        (0.03125,    9.91923e-09)
    };

    \addplot [black] plot coordinates {

        (.0625,    0.0080)
        (0.03125,    0.002)
    } node[pos=0.5,anchor=north]{$2$};

    \addplot [black] plot coordinates {

        (.0625,    4.8000e-05)
        (0.03125,    3e-6)
    } node[pos=0.5,anchor=north]{$4$};

    \addplot [black] plot coordinates {

        (.0625,    2.5600e-06)
        (0.03125,    8e-8)
    } node[pos=0.5,anchor=north]{$5$};

    \addplot [black] plot coordinates {

        (.0625,    4.4800e-07)
        (0.03125,    7e-9)
    } node[pos=0.5,anchor=north]{$6$};

    \legend{G-$Q_2$\\G-$Q_3$\\G-$Q_4$\\G-$Q_5$\\B-$Q_2$\\B-$Q_3$\\B-$Q_4$\\B-$Q_5$\\}
    \end{loglogaxis}
\end{tikzpicture}

%% file: tikz_scripts/two_patch_biharmonic_basic_H2.tex
\begin{tikzpicture}
    \begin{loglogaxis}[
        legend columns=4,
    	legend style={at={(1,1.2)}, nodes={scale=0.8, transform shape}, column sep=.2cm},
        xlabel=$h$,
        ylabel=${\|u-u^{h}\|_{H^2_*}}/{\|u\|_{H^2_*}}$ 
    ]

    \addplot [color=red,mark=triangle] plot coordinates {

        (.5,    0.971535)
        (.25,   0.486982)
        (.125,   0.239064)
        (.0625,   0.120088)
        (0.03125,   0.0601225)
    };

    \addplot [color=blue,mark=triangle] plot coordinates {

        (.5,    0.703071)
        (.25,   0.20652)
        (.125,   0.0349283)
        (.0625,   0.00788988)
        (0.03125,   0.001924 )
    };

    \addplot [color=black,mark=triangle] plot coordinates {

        (.5,    0.280586)
        (.25,   0.0564605)
        (.125,   0.00632958)
        (.0625,  0.000652806)
        (0.03125,   7.73383e-05)
    };

    \addplot [color=green,mark=triangle] plot coordinates {

        (.5,    0.217478)
        (.25,   0.0440545)
        (.125,   0.00121006)
        (.0625,   5.74852e-05)
        (0.03125,    3.95598e-06)
    };

    \addplot [color=red,mark=triangle, dashed] plot coordinates {

        (.5,    0.984451 )
        (.25,   0.569577)
        (.125,   0.247429)
        (.0625,   0.120389)
        (0.03125,   0.0601331)
    };

    \addplot [color=blue,mark=triangle, dashed] plot coordinates {

        (.5,    0.723664)
        (.25,   0.244872)
        (.125,   0.0356371)
        (.0625,   0.0081103)
        (0.03125,   0.00194256)
    };

    \addplot [color=black,mark=triangle, dashed] plot coordinates {

        (.5,    0.728375)
        (.25,   0.104051)
        (.125,   0.00776761)
        (.0625,   0.000673919)
        (0.03125,  7.77689e-05)
    };

    \addplot [color=green,mark=triangle, dashed] plot coordinates {

        (.5,    0.428472)
        (.25,   0.065016)
        (.125,   0.001320065)
        (.0625,   6.40474e-05)
        (0.03125,    3.70118e-06)
    };

    \addplot [black] plot coordinates {

        (.0625,    0.0800)
        (0.03125,    .04)
    } node[pos=0.5,anchor=north]{$1$};

    \addplot [black] plot coordinates {

        (.0625,    0.004)
        (0.03125,    .001)
    } node[pos=0.5,anchor=north]{$2$};

    \addplot [black] plot coordinates {

        (.0625,    4.0000e-04)
        (0.03125,    5e-5)
    } node[pos=0.5,anchor=north]{$3$};

    \addplot [black] plot coordinates {

        (.0625,    3.2000e-05)
        (0.03125,    2e-6)
    } node[pos=0.5,anchor=north]{$4$};

    \legend{G-$Q_2$\\G-$Q_3$\\G-$Q_4$\\G-$Q_5$\\B-$Q_2$\\B-$Q_3$\\B-$Q_4$\\B-$Q_5$\\}
    \end{loglogaxis}
\end{tikzpicture}

%% file: tikz_scripts/two_patch_biharmonic_distorted.tex
\begin{tikzpicture}
    \begin{loglogaxis}[
        legend columns=4,
    	legend style={at={(1,1.2)}, nodes={scale=0.8, transform shape}, column sep=.2cm},
        xlabel=$h$,
        ylabel=${\|u-u^{h}\|_{L^2}}/{\|u\|_{L^2}}$ 
    ]

    \addplot [color=red,mark=triangle] plot coordinates {

        (.5,    1.08115)
        (.25,   0.792129)
        (.125,   0.322907)
        (.0625,   0.0307356)
        (0.03125,   0.00444968)
        (0.015625,   0.00106043)
    };

    \addplot [color=blue,mark=triangle] plot coordinates {

        (.5,    0.857518)
        (.25,   0.385025)
        (.125,   0.0225406)
        (.0625,   0.00041709)
        (0.03125,   1.09308e-05)
        (0.015625,   4.26084e-07) 
    };

    \addplot [color=black,mark=triangle] plot coordinates {

        (.5,    0.595125)
        (.25,   0.265113)
        (.125,   0.00797683)
        (.0625,   0.000106357)
        (0.03125,   7.13069e-07)
        (0.015625,  8.78607e-09 )
    };

    \addplot [color=green,mark=triangle] plot coordinates {

        (.5,    0.435502 )
        (.25,   0.20598 )
        (.125,   0.00330122)
        (.0625,   5.16129e-06)
        (0.03125,   1.30976e-07)
        (0.015625,   8.33014e-10)
    };

    \addplot [color=red,mark=triangle, dashed] plot coordinates {

        (.5,    1.08115)
        (.25,   0.834959)
        (.125,   0.318748)
        (.0625,   0.0277705)
        (0.03125,  0.00436393)
        (0.015625,   0.0010596 )
    };

    \addplot [color=blue,mark=triangle, dashed] plot coordinates {

        (.5,    0.857518)
        (.25,   0.359827)
        (.125,   0.00807084 )
        (.0625,   0.000358953)
        (0.03125,  1.07221e-05)
        (0.015625,  4.09114e-07)
    };

    \addplot [color=black,mark=triangle, dashed] plot coordinates {

        (.5,    0.595125)
        (.25,   0.125729)
        (.125,   0.0036115)
        (.0625,   8.02087e-06)
        (0.03125,   2.21844e-07)
        (0.015625,  5.87874e-09 )
    };

    \addplot [color=green,mark=triangle, dashed] plot coordinates {

        (.5,    0.359201 )
        (.25,   0.0937847 )
        (.125,   0.000346784)
        (.0625,   4.10663e-06 )
        (0.03125,  2.11333e-08 )
        (0.015625,   2.31939e-10 )
    };

    \addplot [black] plot coordinates {

        (0.03125,       0.008 )
        (0.015625,      .002 )
    } node[pos=0.5,anchor=south]{$2$};

    \addplot [black] plot coordinates {

        (0.03125,       1.2800e-05)
        (0.015625,      8e-7 )
    } node[pos=0.5,anchor=south]{$4$};

    \addplot [black] plot coordinates {

        (0.03125,       3.2000e-07)
        (0.015625,      1e-8 )
    } node[pos=0.5,anchor=south]{$5$};

    \addplot [black] plot coordinates {

        (0.03125,       5.7600e-09)
        (0.015625,      9e-11 )
    } node[pos=0.5,anchor=north]{$6$};

    \legend{G-$Q_2$\\G-$Q_3$\\G-$Q_4$\\G-$Q_5$\\B-$Q_2$\\B-$Q_3$\\B-$Q_4$\\B-$Q_5$\\}
    \end{loglogaxis}
\end{tikzpicture}

%% file: tikz_scripts/two_patch_biharmonic_distorted_H2.tex
\begin{tikzpicture}
    \begin{loglogaxis}[
        legend columns=4,
    	legend style={at={(1,1.2)}, nodes={scale=0.8, transform shape}, column sep=.2cm},
        xlabel=$h$,
        ylabel=${\|u_{}-u^{h}\|_{H^2_*}}/{\|u_{}\|_{H^2_*}}$ 
    ]

    \addplot [color=red,mark=triangle] plot coordinates {

        (.5,    0.970233)
        (.25,   0.630331)
        (.125,   0.347793)
        (.0625,   0.140633)
        (0.03125,   0.063464)
        (0.015625,   0.0313622)
    };

    \addplot [color=blue,mark=triangle] plot coordinates {

        (.5,    0.812838)
        (.25,   0.451057)
        (.125,   0.0804161)
        (.0625,   0.00981443)
        (0.03125,   0.00214694)
        (0.015625,   0.000516027) 
    };

    \addplot [color=black,mark=triangle] plot coordinates {

        (.5,    0.727478)
        (.25,   0.347299)
        (.125,   0.0382459)
        (.0625,   0.00174571)
        (0.03125,   9.75601e-05)
        (0.015625,  1.10561e-05)
    };

    \addplot [color=green,mark=triangle] plot coordinates {

        (.5,     0.5517)
        (.25,   0.310648)
        (.125,   0.0166288)
        (.0625,   0.000217162)
        (0.03125,   1.05336e-05)
        (0.015625,   3.38835e-07)
    };

    \addplot [color=red,mark=triangle, dashed] plot coordinates {

        (.5,    0.970233)
        (.25,   0.598473)
        (.125,   0.341423)
        (.0625,   0.137625)
        (0.03125,  0.0633613)
        (0.015625,  0.0313625)
    };

    \addplot [color=blue,mark=triangle, dashed] plot coordinates {

        (.5,    0.812838)
        (.25,   0.395827)
        (.125,   0.0533701 )
        (.0625,   0.00965695)
        (0.03125,  0.00211437)
        (0.015625,  0.000513982)
    };

    \addplot [color=black,mark=triangle, dashed] plot coordinates {

        (.5,    0.727478)
        (.25,   0.245455)
        (.125,    0.0276997)
        (.0625,    0.000779448)
        (0.03125,  9.08055e-05)
        (0.015625,  1.10167e-05 )
    };

    \addplot [color=green,mark=triangle, dashed] plot coordinates {

        (.5,    0.548646 )
        (.25,   0.204253 )
        (.125,   0.0037194)
        (.0625,   0.000188883)
        (0.03125,  4.98227e-06)
        (0.015625,  2.65114e-07)
    };

    \addplot [black] plot coordinates {

        (0.03125,       0.0400 )
        (0.015625,      .02 )
    } node[pos=0.5,anchor=north]{$1$};

    \addplot [black] plot coordinates {

        (0.03125,       0.0012 )
        (0.015625,      .0003 )
    } node[pos=0.5,anchor=north]{$2$};

    \addplot [black] plot coordinates {

        (0.03125,       5.6000e-05 )
        (0.015625,      7e-6 )
    } node[pos=0.5,anchor=north]{$3$};

    \addplot [black] plot coordinates {

        (0.03125,       3.2000e-06 )
        (0.015625,      2e-7 )
    } node[pos=0.5,anchor=north]{$4$};

    \legend{G-$Q_2$\\G-$Q_3$\\G-$Q_4$\\G-$Q_5$\\B-$Q_2$\\B-$Q_3$\\B-$Q_4$\\B-$Q_5$\\}
    \end{loglogaxis}
\end{tikzpicture}

%% file: tikz_scripts/two_patch_biharmonic_nonmatch.tex
\begin{tikzpicture}
    \begin{loglogaxis}[
        legend columns=4,
    	legend style={at={(1,1.2)}, nodes={scale=0.8, transform shape}, column sep=.2cm},
        xlabel=$h$,
        ylabel=${\|u_{}-u^{h}\|_{L^2}}/{\|u_{}\|_{L^2}}$ 
    ]

    \addplot [color=red,mark=triangle] plot coordinates {

        (.5,    1.1098)
        (.25,   0.508243)
        (.125,   0.179264)
        (.0625,   0.0494358 )
        (0.03125,   0.0128027)
        (0.015625,  0.00170833 )
    };

    \addplot [color=blue,mark=triangle] plot coordinates {

        (.5,    0.738315 )
        (.25,   0.274125 )
        (.125,  0.0438196)
        (.0625,   0.00100768)
        (0.03125,   7.91825e-05)
        (0.015625,  3.00839e-06 )
    };

    \addplot [color=black,mark=triangle] plot coordinates {

        (.5,    0.709611 )
        (.25,  0.243391)
        (.125,  0.00390737 )
        (.0625,  0.000299061)
        (0.03125,  3.69165e-06 )
        (0.015625,  3.6413e-08 )
    };

    \addplot [color=green,mark=triangle] plot coordinates {

        (.5,   0.277756 )
        (.25,   0.14141 )
        (.125,   0.00375236)
        (.0625,   4.49226e-06)
        (0.03125,   1.91177e-07 )
        (0.015625,  2.37417e-08 )
    };

    \addplot [color=red,mark=triangle, dashed] plot coordinates {

        (.5,    1.1098 )
        (.25,   0.486595)
        (.125,   0.551984)
        (.0625,   0.279431 )
        (0.03125,   0.0899244 )
        (0.015625,  0.0248471  )
    };

    \addplot [color=blue,mark=triangle, dashed] plot coordinates {

        (.5,    0.738315 )
        (.25,   0.388806)
        (.125,   0.124982 )
        (.0625,   0.0072821 )
        (0.03125,   0.000314343 )
        (0.015625,  1.20043e-05 )
    };

    \addplot [color=black,mark=triangle, dashed] plot coordinates {

        (.5,    0.709611 )
        (.25,   0.273513 )
        (.125,  0.00201785)
        (.0625,   0.000301985 )
        (0.03125,  1.31055e-05 )
        (0.015625, 4.46325e-07 )
    };

    \addplot [color=green,mark=triangle, dashed] plot coordinates {

        (.5,    0.277614)
        (.25,   0.111033)
        (.125,  0.00225552)
        (.0625,  3.09458e-05 )
        (0.03125,  2.99167e-07 )
        (0.015625, 2.68575e-09  )
    };

    \addplot [black] plot coordinates {

        (0.03125,     0.04)
        (0.015625,    0.01)
    } node[pos=0.5,anchor=north]{$2$};

    \addplot [black] plot coordinates {

        (0.03125,     3.2000e-04)
        (0.015625,    2e-5)
    } node[pos=0.5,anchor=south]{$4$};

    \addplot [black] plot coordinates {

        (0.03125,     6.4000e-06)
        (0.015625,    2e-7)
    } node[pos=0.5,anchor=north]{$5$};

    \addplot [black] plot coordinates {

        (0.03125,     6.4000e-08)
        (0.015625,    1e-9)
    } node[pos=0.5,anchor=north]{$6$};

    \legend{G-$Q_2$\\G-$Q_3$\\G-$Q_4$\\G-$Q_5$\\B-$Q_2$\\B-$Q_3$\\B-$Q_4$\\B-$Q_5$\\}
    \end{loglogaxis}
\end{tikzpicture}

%% file: tikz_scripts/two_patch_biharmonic_nonmatch_H2.tex
\begin{tikzpicture}
    \begin{loglogaxis}[
        legend columns=4,
    	legend style={at={(1,1.2)}, nodes={scale=0.8, transform shape}, column sep=.2cm},
        xlabel=$h$,
        ylabel=${\|u_{}-u^{h}\|_{H^2_*}}/{\|u_{}\|_{H^2_*}}$ 
    ]

    \addplot [color=red,mark=triangle] plot coordinates {

        (.5,    0.984462)
        (.25,   0.57064)
        (.125,   0.29491)
        (.0625,  0.124681)
        (0.03125,   0.0622175)
        (0.015625,  0.0306618)
    };

    \addplot [color=blue,mark=triangle] plot coordinates {

        (.5,    0.73507)
        (.25,   0.330168)
        (.125,  0.0893942)
        (.0625,   0.00996365)
        (0.03125,    0.00251387)
        (0.015625,  0.000536017)
    };

    \addplot [color=black,mark=triangle] plot coordinates {

        (.5,     0.737524)
        (.25,  0.227844)
        (.125,  0.0146347)
        (.0625,  0.00359733)
        (0.03125,  0.000183658)
        (0.015625,  1.17426e-05)
    };

    \addplot [color=green,mark=triangle] plot coordinates {

        (.5,   0.467302)
        (.25,  0.237741)
        (.125,   0.020301)
        (.0625,   0.000121442)
        (0.03125,   1.49368e-05 )
        (0.015625,  6.91432e-06 )
    };

    \addplot [color=red,mark=triangle, dashed] plot coordinates {

        (.5,    0.984462)
        (.25,   0.566638)
        (.125,   0.294139)
        (.0625,   0.146838)
        (0.03125,   0.076735)
        (0.015625,  0.0388408)
    };

    \addplot [color=blue,mark=triangle, dashed] plot coordinates {

        (.5,    0.73507 )
        (.25,   0.316245)
        (.125,   0.102614)
        (.0625,   0.0163837)
        (0.03125,   0.00280983 )
        (0.015625,  0.00055078 )
    };

    \addplot [color=black,mark=triangle, dashed] plot coordinates {

        (.5,    0.737524 )
        (.25,   0.239446 )
        (.125,  0.0120767)
        (.0625,   0.00340827 )
        (0.03125,  0.000525168 )
        (0.015625, 7.06666e-05 )
    };

    \addplot [color=green,mark=triangle, dashed] plot coordinates {

        (.5,    0.467299)
        (.25,   0.224219)
        (.125,  0.0154553)
        (.0625,  0.000614266 )
        (0.03125, 2.33943e-05)
        (0.015625, 8.5008e-07)
    };

    \addplot [black] plot coordinates {

        (0.03125,     0.04)
        (0.015625,    0.02)
    } node[pos=0.5,anchor=north]{$1$};

    \addplot [black] plot coordinates {

        (0.03125,     0.0016)
        (0.015625,    0.0004)
    } node[pos=0.5,anchor=north]{$2$};

    \addplot [black] plot coordinates {

        (0.03125,     4e-4)
        (0.015625,    5e-5)
    } node[pos=0.5,anchor=north]{$3$};

    \addplot [black] plot coordinates {

        (0.03125,     9.6000e-06)
        (0.015625,    6e-7)
    } node[pos=0.5,anchor=north]{$4$};

    \legend{G-$Q_2$\\G-$Q_3$\\G-$Q_4$\\G-$Q_5$\\B-$Q_2$\\B-$Q_3$\\B-$Q_4$\\B-$Q_5$\\}
    \end{loglogaxis}
\end{tikzpicture}

%% file: tikz_scripts/two_patch_biharmonic_diff_degree.tex
\begin{tikzpicture}
    \begin{loglogaxis}[
        legend columns=3,
    	legend style={at={(1,1.2)}, nodes={scale=0.8, transform shape}, column sep=.2cm},
        xlabel=$h$,
        ylabel=${\|u_{}-u^{h}\|_{L^2}}/{\|u_{}\|_{L^2}}$ 
    ]

    \addplot [color=red,mark=triangle] plot coordinates {

        (.5,    1.02271)
        (.25,   0.393272 )
        (.125,  0.0360061 )
        (.0625,   0.00802147 )
        (0.03125,   0.00193664 )
        (0.015625,  0.00048126 )
    };

    \addplot [color=blue,mark=triangle] plot coordinates {

        (.5,    0.694232)
        (.25,   0.166681)
        (.125,   0.0040645)
        (.0625,   6.81819e-05 )
        (0.03125,   1.6304e-06)
        (0.015625,  7.29714e-08)
    };

    \addplot [color=black,mark=triangle] plot coordinates {

        (.5,    0.351814)
        (.25,   0.0178482)
        (.125,      0.000449695)
        (.0625,     2.92727e-06)
        (0.03125,   3.81391e-08)
        (0.015625,  7.7457e-10 )
    };

    \addplot [color=red,mark=triangle, dashed] plot coordinates {

        (.5,    1.02271)
        (.25,   0.461445)
        (.125,   0.111545)
        (.0625,  0.0109794)
        (0.03125,   0.00201306)
        (0.015625,  0.000483481 )
    };

    \addplot [color=blue,mark=triangle, dashed] plot coordinates {

        (.5,    0.694232)
        (.25,   0.213121)
        (.125,      0.00180538)
        (.0625,     0.000212459 )
        (0.03125,   8.31712e-06)
        (0.015625,  2.7879e-07  )
    };

    \addplot [color=black,mark=triangle, dashed] plot coordinates {

        (.5,    0.351814 )
        (.25,   0.0186012)
        (.125,      0.00123195)
        (.0625,     1.40843e-05)
        (0.03125,   1.21223e-07)
        (0.015625,  1.17654e-09 )
    };

    \addplot [black] plot coordinates {

        (0.03125,    0.0012)
        (0.015625,    .0003)
    } node[pos=0.5,anchor=north]{$2$};

    \addplot [black] plot coordinates {

        (0.03125,    9.6000e-07)
        (0.015625,   6e-8)
    } node[pos=0.5,anchor=north]{$4$};

    \addplot [black] plot coordinates {

        (0.03125,    1.9200e-08)
        (0.015625,   6e-10)
    } node[pos=0.5,anchor=north]{$5$};

    \legend{G-$Q_2$ $\&$ $Q_3$\\G-$Q_3$ $\&$ $Q_4$\\G-$Q_4$ $\&$ $Q_5$\\B-$Q_2$ $\&$ $Q_3$\\B-$Q_3$ $\&$ $Q_4$\\B-$Q_4$ $\&$ $Q_5$\\}
    \end{loglogaxis}
\end{tikzpicture}

%% file: tikz_scripts/two_patch_biharmonic_diff_degree_H2.tex
\begin{tikzpicture}
    \begin{loglogaxis}[
        legend columns=3,
    	legend style={at={(1,1.2)}, nodes={scale=0.8, transform shape}, column sep=.2cm},
        xlabel=$h$,
        ylabel=${\|u_{}-u^{h}\|_{H^2_*}}/{\|u_{}\|_{H^2_*}}$ 
    ]

    \addplot [color=red,mark=triangle] plot coordinates {

        (.5,    0.947573)
        (.25,   0.452374)
        (.125,  0.150636)
        (.0625,   0.0738313)
        (0.03125,   0.0366665 )
        (0.015625,  0.0183142 )
    };

    \addplot [color=blue,mark=triangle] plot coordinates {

        (.5,    0.675618)
        (.25,   0.29897)
        (.125,      0.0263881)
        (.0625,     0.00385068)
        (0.03125,   0.000894183)
        (0.015625,  0.000221005)
    };

    \addplot [color=black,mark=triangle] plot coordinates {

        (.5,    0.597089)
        (.25,   0.069544)
        (.125,      0.00404805)
        (.0625,     0.00023493)
        (0.03125,   2.65706e-05)
        (0.015625,  3.25804e-06)
    };

    \addplot [color=red,mark=triangle, dashed] plot coordinates {

        (.5,    0.947573)
        (.25,   0.448626)
        (.125,      0.165961)
        (.0625,     0.0745943)
        (0.03125,   0.0366899)
        (0.015625,  0.0183149)
    };

    \addplot [color=blue,mark=triangle, dashed] plot coordinates {

        (.5,    0.675618)
        (.25,   0.300972)
        (.125,      0.0211338)
        (.0625,     0.0046142 )
        (0.03125,   0.000991629)
        (0.015625,  0.0002282)
    };

    \addplot [color=black,mark=triangle, dashed] plot coordinates {

        (.5,    0.597089)
        (.25,   0.0697155)
        (.125,      0.00879059)
        (.0625,     0.000444331)
        (0.03125,   2.93938e-05)
        (0.015625,  3.28279e-06)
    };

    \addplot [black] plot coordinates {

        (0.03125,    0.0180)
        (0.015625,   .009)
    } node[pos=0.5,anchor=north]{$1$};

    \addplot [black] plot coordinates {

        (0.03125,    3.6000e-04)
        (0.015625,   0.00009)
    } node[pos=0.5,anchor=north]{$2$};

    \addplot [black] plot coordinates {

        (0.03125,    8.0000e-06)
        (0.015625,   1e-6)
    } node[pos=0.5,anchor=north]{$3$};

    \legend{G-$Q_2$ $\&$ $Q_3$\\G-$Q_3$ $\&$ $Q_4$\\G-$Q_4$ $\&$ $Q_5$\\B-$Q_2$ $\&$ $Q_3$\\B-$Q_3$ $\&$ $Q_4$\\B-$Q_4$ $\&$ $Q_5$\\}
    \end{loglogaxis}
\end{tikzpicture}

%% file: tikz_scripts/three_patch_modify_biharmonic_basic.tex
\begin{tikzpicture}
    \begin{loglogaxis}[
        legend columns=4,
    	legend style={at={(1,1.2)}, nodes={scale=0.8, transform shape}, column sep=.1cm},
        xlabel=$h$,
        ylabel=${\|u-u^{h}\|_{L^2}}/{\|u\|_{L^2}}$ 
    ]

    \addplot [color=red,mark=triangle] plot coordinates {

        (.25,   0.327832)
        (.125,   0.0710077)
        (.0625,   0.0176199)
        (0.03125,  0.00438157)
        (0.015625,  0.00109072)
    };

    \addplot [color=blue,mark=triangle] plot coordinates {

        (.25,   .1294 )
        (.125,      .009544 )
        (.0625,     .000317 )
        (0.03125,   1.623e-5 )
        (0.015625,  9.595e-7 )
    };

    \addplot [color=black,mark=triangle] plot coordinates {

        (.25,   .1295)
        (.125,      .003024   )
        (.0625,     3.335e-5   )
        (0.03125,   7.118e-07  )
        (0.015625,  2.0136e-08  )
    };

    \addplot [color=green,mark=triangle] plot coordinates {

        (.25,   0.04626)
        (.125,   0.001240)
        (.0625,   5.994e-06)
        (0.03125,   5.92657e-08 )
        (0.015625,  9.16817e-10  )
    };

    \addplot [color=red,mark=triangle, dashed] plot coordinates {

        (.5,        0.757134)
        (.25,       0.33531)
        (.125,      0.0713054)
        (.0625,     0.0176182)
        (0.03125,   0.0043803)
        (0.015625,  0.001091)
    };

    \addplot [color=blue,mark=triangle, dashed] plot coordinates {

        (.5,        0.583603)
        (.25,       0.116352)
        (.125,      0.00952435)
        (.0625,     0.000547887)
        (0.03125,   0.000180265 )
        (0.015625,  5.13425e-05 )
    };

    \addplot [color=black,mark=triangle, dashed] plot coordinates {

        (.5,        0.336308)
        (.25,       0.125011)
        (.125,      0.00314383)
        (.0625,     0.000707717)
        (0.03125,   0.000228204)
        (0.015625,  6.18458e-05 )
    };

    \addplot [color=green,mark=triangle, dashed] plot coordinates {

        (.5,        0.399626 )
        (.25,       0.0522817   )
        (.125,      0.00211541  )
        (.0625,     0.000852918)
        (0.03125,       0.000254254)
        (0.015625,      6.74393e-05 )
    };

    \addplot [black] plot coordinates {

        (0.03125,    0.0032)
        (0.015625,    .0008)
    } node[pos=0.5,anchor=north]{$2$};

    \addplot [black] plot coordinates {

        (0.03125,    1.1200e-05)
        (0.015625,   7e-7)
    } node[pos=0.5,anchor=north]{$4$};

    \addplot [black] plot coordinates {

        (0.03125,    3.200e-07)
        (0.015625,   1e-8)
    } node[pos=0.5,anchor=north]{$5$};

    \addplot [black] plot coordinates {

        (0.03125,    3.2000e-08)
        (0.015625,   5e-10)
    } node[pos=0.5,anchor=north]{$6$};

    \legend{MG-$Q_2$\\MG-$Q_3$\\MG-$Q_4$\\MG-$Q_5$\\MB-$Q_2$\\MB-$Q_3$\\MB-$Q_4$\\MB-$Q_5$\\}
    \end{loglogaxis}
\end{tikzpicture}

%% file: tikz_scripts/three_patch_modify_biharmonic_basic_H2.tex
\begin{tikzpicture}
    \begin{loglogaxis}[
        legend columns=4,
    	legend style={at={(1,1.2)}, nodes={scale=0.8, transform shape}, column sep=.1cm},
        xlabel=$h$,
        ylabel=${\|u_{}-u^{h}\|_{H^2_*}}/{\|u_{}\|_{H^2_*}}$ 
    ]

    \addplot [color=red,mark=triangle] plot coordinates {

        (.25,   0.510891)
        (.125,   0.251056)
        (.0625,   0.126542)
        (0.03125,   0.0633989)
        (0.015625,  0.0316347 )
    };

    \addplot [color=blue,mark=triangle] plot coordinates {

        (.25,       0.262475)
        (.125,      0.0632935)
        (.0625,     0.0124056)
        (0.03125,   0.00289077)
        (0.015625,  0.0007086345 )
    };

    \addplot [color=black,mark=triangle] plot coordinates {

        (.25,       0.241842)
        (.125,      0.0231632)
        (.0625,     0.00178264)
        (0.03125,   0.000200438)
        (0.015625,  2.35884e-05 )
    };

    \addplot [color=green,mark=triangle] plot coordinates {

        (.25,       0.114182)
        (.125,      0.00956557)
        (.0625,     0.000349277)
        (0.03125,   1.94526e-05)
        (0.015625,  1.21051e-06 )
    };

    \addplot [color=red,mark=triangle, dashed] plot coordinates {

        (.5,        0.845793 )
        (.25,       0.518268)
        (.125,      0.251447)
        (.0625,     0.126541)
        (0.03125,   0.0634613)
        (0.015625,  0.0316795 )
    };

    \addplot [color=blue,mark=triangle, dashed] plot coordinates {

        (.5,        0.678273)
        (.25,       0.271455)
        (.125,      0.0646432)
        (.0625,     0.0140928)
        (0.03125,   0.00497191)
        (0.015625,  0.00230988 )
    };

    \addplot [color=black,mark=triangle, dashed] plot coordinates {

        (.5,    0.480135)
        (.25,   0.249006)
        (.125,      0.0292573)
        (.0625,     0.00963198)
        (0.03125,   0.00526865)
        (0.015625,  0.00277059)
    };

    \addplot [color=green,mark=triangle, dashed] plot coordinates {

        (.5,        0.62434)
        (.25,       0.142799)
        (.125,      0.025262)
        (.0625,     0.0118577)
        (0.03125,   0.00631867)
        (0.015625,  0.00326862 )
    };

    \addplot [black] plot coordinates {

        (.03125,    0.0400)
        (0.015625,    .02)
    } node[pos=0.5,anchor=north]{$1$};

    \addplot [black] plot coordinates {

        (.03125,    0.0020)
        (0.015625,    .0005)
    } node[pos=0.5,anchor=north]{$2$};

    \addplot [black] plot coordinates {

        (.03125,    1.2000e-04)
        (0.015625,  1.5e-5)
    } node[pos=0.5,anchor=north]{$3$};

    \addplot [black] plot coordinates {

        (.03125,    1.4400e-05)
        (0.015625,  9e-7)
    } node[pos=0.5,anchor=north]{$4$};

    \legend{MG-$Q_2$\\MG-$Q_3$\\MG-$Q_4$\\MG-$Q_5$\\MB-$Q_2$\\MB-$Q_3$\\MB-$Q_4$\\MB-$Q_5$\\}
    \end{loglogaxis}
\end{tikzpicture}

%% file: tikz_scripts/three_patch_biharmonic_basic.tex
\begin{tikzpicture}
    \begin{loglogaxis}[
        legend columns=4,
    	legend style={at={(1,1.2)}, nodes={scale=0.8, transform shape}, column sep=.1cm},
        xlabel=$h$,
        ylabel=${\|u-u^{h}\|_{L^2}}/{\|u\|_{L^2}}$ 
    ]

    \addplot [color=red,mark=triangle] plot coordinates {

        (.5,    0.867622)
        (.25,   0.550937)
        (.125,   0.123841)
        (.0625,   0.0204693)
        (0.03125,  0.00453748)
        (0.015625,  0.00109953)
    };

    \addplot [color=blue,mark=triangle] plot coordinates {

        (.5,    0.709561 )
        (.25,   0.293877 )
        (.125,      0.0144029 )
        (.0625,     0.000576479)
        (0.03125,   2.17777e-05 )
        (0.015625,  1.05547e-06 )
    };

    \addplot [color=black,mark=triangle] plot coordinates {

        (.5,    0.533089)
        (.25,   0.127303)
        (.125,      0.00373231   )
        (.0625,     3.40865e-05   )
        (0.03125,   7.44076e-07  )
        (0.015625,  2.05249e-08  )
    };

    \addplot [color=green,mark=triangle] plot coordinates {

        (.5,    0.210034)
        (.25,   0.067279)
        (.125,   0.00131031)
        (.0625,   6.95576e-06)
        (0.03125,   6.12699e-08 )
        (0.015625,  2.84736e-09  )
    };

    \addplot [color=red,mark=triangle, dashed] plot coordinates {

        (.5,        0.95045)
        (.25,       0.409736)
        (.125,      0.0896597)
        (.0625,     0.0190763)
        (0.03125,   0.00447226)
        (0.015625,  0.00109614  )
    };

    \addplot [color=blue,mark=triangle, dashed] plot coordinates {

        (.5,        0.67082)
        (.25,       0.149727 )
        (.125,      0.00929605)
        (.0625,     0.000318077)
        (0.03125,   1.62167e-05 )
        (0.015625,  9.59342e-07 )
    };

    \addplot [color=black,mark=triangle, dashed] plot coordinates {

        (.5,        0.389482)
        (.25,       0.129543)
        (.125,      0.00302886)
        (.0625,     3.33622e-05)
        (0.03125,   7.14709e-07)
        (0.015625,  2.06605e-08 )
    };

    \addplot [color=green,mark=triangle, dashed] plot coordinates {

        (.5,        0.186975 )
        (.25,       0.0478791   )
        (.125,      0.00124043  )
        (.0625,     5.99822e-06)
        (0.03125,       6.66757e-08)
        (0.015625,      1.83655e-09 )
    };

    \addplot [black] plot coordinates {

        (0.03125,    0.0032)
        (0.015625,    .0008)
    } node[pos=0.5,anchor=north]{$2$};

    \addplot [black] plot coordinates {

        (0.03125,    1.1200e-05)
        (0.015625,   7e-7)
    } node[pos=0.5,anchor=north]{$4$};

    \addplot [black] plot coordinates {

        (0.03125,    3.200e-07)
        (0.015625,   1e-8)
    } node[pos=0.5,anchor=north]{$5$};

    \addplot [black] plot coordinates {

        (.0625,    2.5600e-06)
        (0.03125,  4e-8)
    } node[pos=0.5,anchor=north]{$6$};

    \legend{OG-$Q_2$\\OG-$Q_3$\\OG-$Q_4$\\OG-$Q_5$\\OB-$Q_2$\\OB-$Q_3$\\OB-$Q_4$\\OB-$Q_5$\\}
    \end{loglogaxis}
\end{tikzpicture}

%% file: tikz_scripts/three_patch_biharmonic_basic_H2.tex
\begin{tikzpicture}
	\begin{loglogaxis}[
			legend columns=4,
			legend style={at={(1,1.2)}, nodes={scale=0.8, transform shape}, column sep=.1cm},
			xlabel=$h$,
			ylabel=${\|u_{}-u^{h}\|_{H^2_*}}/{\|u_{}\|_{H^2_*}}$
		]

		\addplot [color=red,mark=triangle] plot coordinates {

				(.5,    0.888913)
				(.25,   0.61016)
				(.125,   0.3112)
				(.0625,   0.135786)
				(0.03125,   0.0648747)
				(0.015625,  0.0318627 )
			};

		\addplot [color=blue,mark=triangle] plot coordinates {

				(.5,        0.710357)
				(.25,       0.436158)
				(.125,      0.0842743)
				(.0625,     0.0154776)
				(0.03125,   0.00309599 )
				(0.015625,  0.000720296 )
			};

		\addplot [color=black,mark=triangle] plot coordinates {

				(.5,        0.595272)
				(.25,       0.252901)
				(.125,      0.0299896)
				(.0625,     0.00182205)
				(0.03125,   0.000191738)
				(0.015625,  2.32465e-05 )
			};

		\addplot [color=green,mark=triangle] plot coordinates {

				(.5,        0.36944)
				(.25,       0.190279)
				(.125,      0.0108137)
				(.0625,     0.000348933)
				(0.03125,   1.69576e-05)
				(0.015625,  4.43769e-06 )
			};

		\addplot [color=red,mark=triangle, dashed] plot coordinates {

				(.5,        0.893624 )
				(.25,       0.54519)
				(.125,      0.277021)
				(.0625,     0.132993)
				(0.03125,   0.064461)
				(0.015625,  0.031804 )
			};

		\addplot [color=blue,mark=triangle, dashed] plot coordinates {

				(.5,        0.706655)
				(.25,       0.31246)
				(.125,      0.0653824)
				(.0625,     0.0125155)
				(0.03125,   0.00289372)
				(0.015625,  0.000708549 )
			};

		\addplot [color=black,mark=triangle, dashed] plot coordinates {

				(.5,    0.515317)
				(.25,   0.242107)
				(.125,      0.0233544)
				(.0625,     0.00177646)
				(0.03125,   0.000188749)
				(0.015625,  2.3109e-05 )
			};

		\addplot [color=green,mark=triangle, dashed] plot coordinates {

				(.5,        0.356252)
				(.25,       0.12105)
				(.125,      0.00957498)
				(.0625,     0.000298767)
				(0.03125,   1.9763e-05)
				(0.015625,  4.49951e-06 )
			};

		\addplot [black] plot coordinates {

				(.03125,    0.0400)
				(0.015625,    .02)
			} node[pos=0.5,anchor=north]{$1$};

		\addplot [black] plot coordinates {

				(.03125,    0.0020)
				(0.015625,    .0005)
			} node[pos=0.5,anchor=north]{$2$};

		\addplot [black] plot coordinates {

				(.03125,    1.2000e-04)
				(0.015625,  1.5e-5)
			} node[pos=0.5,anchor=north]{$3$};

		\addplot [black] plot coordinates {

				(.0625,    1.6000e-04)
				(0.03125,   1e-5)
			} node[pos=0.5,anchor=north]{$4$};

		\legend{OG-$Q_2$\\OG-$Q_3$\\OG-$Q_4$\\OG-$Q_5$\\OB-$Q_2$\\OB-$Q_3$\\OB-$Q_4$\\OB-$Q_5$\\}
	\end{loglogaxis}
\end{tikzpicture}

%% file: tikz_scripts/three_patch_modify_projection_basic_H2.tex
\begin{tikzpicture}
    \begin{loglogaxis}[
        legend columns=4,
    	legend style={at={(1,1.2)}, nodes={scale=0.8, transform shape}, column sep=.1cm},
        xlabel=$h$,
        ylabel=${\|u_{}-u^{h}\|_{H^2_*}}/{\|u_{}\|_{H^2_*}}$ 
    ]

    \addplot [color=red,mark=triangle] plot coordinates {

        (.25,   0.514976)
        (.125,   0.252957)
        (.0625,   0.127637)
        (0.03125,   0.0634559)
        (0.015625,  0.0316453 )
    };

    \addplot [color=blue,mark=triangle] plot coordinates {

        (.25,       0.265909)
        (.125,      0.0626049)
        (.0625,     0.0123811)
        (0.03125,   0.00288704 )
        (0.015625,  0.000707917 )
    };

    \addplot [color=black,mark=triangle] plot coordinates {

        (.25,       0.240451)
        (.125,      0.0231358)
        (.0625,     0.00176459)
        (0.03125,   0.000186888)
        (0.015625,  2.23978e-05 )
    };

    \addplot [color=green,mark=triangle] plot coordinates {

        (.25,       0.120401)
        (.125,      0.00955611)
        (.0625,     0.000295519)
        (0.03125,   1.43898e-05)
        (0.015625,  8.38406e-07 )
    };

    \addplot [color=red,mark=triangle, dashed] plot coordinates {

        (.5,        0.850424 )
        (.25,       0.510382)
        (.125,      0.253449)
        (.0625,     0.127572)
        (0.03125,   0.0633865)
        (0.015625,  0.0316417)
    };

    \addplot [color=blue,mark=triangle, dashed] plot coordinates {

        (.5,        0.654709)
        (.25,       0.265234)
        (.125,      0.0623662)
        (.0625,     0.0123803)
        (0.03125,   0.00288687)
        (0.015625,  0.000707856)
    };

    \addplot [color=black,mark=triangle, dashed] plot coordinates {

        (.5,    0.415372)
        (.25,   0.237414)
        (.125,      0.0231152)
        (.0625,     0.00175536)
        (0.03125,   0.000186804)
        (0.015625,  2.23953e-05)
    };

    \addplot [color=green,mark=triangle, dashed] plot coordinates {

        (.5,        0.343907)
        (.25,       0.120176)
        (.125,      0.00955125)
        (.0625,     0.000295413)
        (0.03125,   1.43845e-05)
        (0.015625,  8.38406e-07)
    };

    \addplot [black] plot coordinates {

        (.03125,    0.0400)
        (0.015625,    .02)
    } node[pos=0.5,anchor=north]{$1$};

    \addplot [black] plot coordinates {

        (.03125,    0.0020)
        (0.015625,    .0005)
    } node[pos=0.5,anchor=north]{$2$};

    \addplot [black] plot coordinates {

        (.03125,    1.2000e-04)
        (0.015625,  1.5e-5)
    } node[pos=0.5,anchor=north]{$3$};

    \addplot [black] plot coordinates {

        (.03125,    9.6000e-06)
        (0.015625,  6e-7)
    } node[pos=0.5,anchor=north]{$4$};

    \legend{MG-$Q_2$\\MG-$Q_3$\\MG-$Q_4$\\MG-$Q_5$\\MB-$Q_2$\\MB-$Q_3$\\MB-$Q_4$\\MB-$Q_5$\\}
    \end{loglogaxis}
\end{tikzpicture}

%% file: tikz_scripts/three_patch_projection_basic_H2.tex
\begin{tikzpicture}
    \begin{loglogaxis}[
        legend columns=4,
    	legend style={at={(1,1.2)}, nodes={scale=0.8, transform shape}, column sep=.1cm},
        xlabel=$h$,
        ylabel=${\|u_{}-u^{h}\|_{H^2_*}}/{\|u_{}\|_{H^2_*}}$ 
    ]

    \addplot [color=red,mark=triangle] plot coordinates {

        (.5,    0.888234)
        (.25,   0.61016)
        (.125,   0.3112)
        (.0625,   0.135786)
        (0.03125,   0.0648747)
        (0.015625,  0.0318627 )
    };

    \addplot [color=blue,mark=triangle] plot coordinates {

        (.5,        0.710357)
        (.25,       0.436158)
        (.125,      0.0842743)
        (.0625,     0.0154776)
        (0.03125,   0.00309599 )
        (0.015625,  0.000720296 )
    };

    \addplot [color=black,mark=triangle] plot coordinates {

        (.5,        0.595272)
        (.25,       0.252901)
        (.125,      0.0299896)
        (.0625,     0.00182205)
        (0.03125,   0.000191186)
        (0.015625,  2.25474e-05 )
    };

    \addplot [color=green,mark=triangle] plot coordinates {

        (.5,        0.36876)
        (.25,       0.190279)
        (.125,      0.0108137)
        (.0625,     0.000348933)
        (0.03125,   1.49534e-05)
        (0.015625,  2.55292e-06 )
    };

    \addplot [color=red,mark=triangle, dashed] plot coordinates {

        (.5,        0.893624 )
        (.25,       0.54519)
        (.125,      0.277021)
        (.0625,     0.132993)
        (0.03125,   0.064461)
        (0.015625,  0.031804 )
    };

    \addplot [color=blue,mark=triangle, dashed] plot coordinates {

        (.5,        0.706655)
        (.25,       0.31246)
        (.125,      0.0653824)
        (.0625,     0.0125155)
        (0.03125,   0.00289372)
        (0.015625,  0.000708549 )
    };

    \addplot [color=black,mark=triangle, dashed] plot coordinates {

        (.5,    0.515317)
        (.25,   0.242107)
        (.125,      0.0233544)
        (.0625,     0.00177646)
        (0.03125,   0.000188749)
        (0.015625,  2.24312e-05)
    };

    \addplot [color=green,mark=triangle, dashed] plot coordinates {

        (.5,        0.356252)
        (.25,       0.12105)
        (.125,      0.00957498)
        (.0625,     0.000298767)
        (0.03125,   1.43877e-05)
        (0.015625,  8.38406e-07)
    };

    \addplot [black] plot coordinates {

        (.03125,    0.0400)
        (0.015625,    .02)
    } node[pos=0.5,anchor=north]{$1$};

    \addplot [black] plot coordinates {

        (.03125,    0.0020)
        (0.015625,    .0005)
    } node[pos=0.5,anchor=north]{$2$};

    \addplot [black] plot coordinates {

        (.03125,    1.2000e-04)
        (0.015625,  1.5e-5)
    } node[pos=0.5,anchor=north]{$3$};

    \addplot [black] plot coordinates {

        (.03125,    9.6000e-06)
        (0.015625,  6e-7)
    } node[pos=0.5,anchor=north]{$4$};

    \legend{OG-$Q_2$\\OG-$Q_3$\\OG-$Q_4$\\OG-$Q_5$\\OB-$Q_2$\\OB-$Q_3$\\OB-$Q_4$\\OB-$Q_5$\\}
    \end{loglogaxis}
\end{tikzpicture}

%% file: tikz_scripts/five_patch_modify_biharmonic_basic.tex
\begin{tikzpicture}
    \begin{loglogaxis}[
        legend columns=4,
    	legend style={at={(1,1.2)}, nodes={scale=0.8, transform shape}, column sep=.1cm},
        xlabel=$h$,
        ylabel=${\|u-u^{h}\|_{L^2}}/{\|u\|_{L^2}}$ 
    ]
    
    
    

    \addplot [color=red,mark=triangle] plot coordinates {

        (.125,   0.1516)
        (.0625,   0.02663)
        (0.03125,  0.006385)
        (0.015625,  0.00156659)
    };

    \addplot [color=blue,mark=triangle] plot coordinates {

        (.125,      0.0219934 )
        (.0625,     0.00112746)
        (0.03125,   4.90144e-05 )
        (0.015625,  2.61451e-06 )
    };

    \addplot [color=black,mark=triangle] plot coordinates {

        (.125,      0.0200278   )
        (.0625,     0.000194223   )
        (0.03125,   3.64379e-06  )
        (0.015625,  9.44966e-08  )
    };

    \addplot [color=green,mark=triangle] plot coordinates {

        (.125,   0.0021798)
        (.0625,   4.79902e-05)
        (0.03125,   3.80295e-07 )
        (0.015625,  4.67128e-09  )
    };

    \addplot [color=red,mark=triangle, dashed] plot coordinates {

        (.25,       0.55309)
        (.125,      0.115654)
        (.0625,     0.0252177)
        (0.03125,   0.00615597)
        (0.015625,  0.00151696 )
    };

    \addplot [color=blue,mark=triangle, dashed] plot coordinates {

        (.25,       0.154026)
        (.125,      0.0209316)
        (.0625,     0.00299256)
        (0.03125,   0.000729679 )
        (0.015625,  0.000169969 )
    };

    \addplot [color=black,mark=triangle, dashed] plot coordinates {

        (.25,       0.135409)
        (.125,      0.0218843)
        (.0625,     0.00372122)
        (0.03125,   0.000814624)
        (0.015625,  0.000185056)
    };

    \addplot [color=green,mark=triangle, dashed] plot coordinates {

        (.25,       0.0417999)
        (.125,      0.0173795 )
        (.0625,     0.00369684)
        (0.03125,       0.000814441)
        (0.015625,      0.000187864 )
    };

    \addplot [black] plot coordinates {

        (0.03125,    4.0000e-03)
        (0.015625,    .001)
    } node[pos=0.5,anchor=north]{$2$};

    \addplot [black] plot coordinates {

        (0.03125,    2.5600e-05)
        (0.015625,   1.6e-6)
    } node[pos=0.5,anchor=north]{$4$};

    \addplot [black] plot coordinates {

        (0.03125,    2.3680e-06)
        (0.015625,   7.4e-8)
    } node[pos=0.5,anchor=north]{$5$};

    \addplot [black] plot coordinates {

        (0.03125,    1.9200e-07)
        (0.015625,   3e-9)
    } node[pos=0.5,anchor=north]{$6$};

    \legend{MG-$Q_2$\\MG-$Q_3$\\MG-$Q_4$\\MG-$Q_5$\\MB-$Q_2$\\MB-$Q_3$\\MB-$Q_4$\\MB-$Q_5$\\}
    \end{loglogaxis}
\end{tikzpicture}

%% file: tikz_scripts/five_patch_modify_biharmonic_basic_H2.tex
\begin{tikzpicture}
    \begin{loglogaxis}[
        legend columns=4,
    	legend style={at={(1,1.2)}, nodes={scale=0.8, transform shape}, column sep=.1cm},
        xlabel=$h$,
        ylabel=${\|u_{}-u^{h}\|_{H^2_*}}/{\|u_{}\|_{H^2_*}}$ 
    ]

    \addplot [color=red,mark=triangle] plot coordinates {

        (.125,   0.429158)
        (.0625,   0.190387)
        (0.03125,   0.0937169)
        (0.015625,  0.0465883)
    };

    \addplot [color=blue,mark=triangle] plot coordinates {

        (.125,      0.125872)
        (.0625,     0.0260387)
        (0.03125,   0.00570847)
        (0.015625,  0.00136665)
    };

    \addplot [color=black,mark=triangle] plot coordinates {

        (.125,      0.0794173)
        (.0625,     0.00467507)
        (0.03125,   0.000462874)
        (0.015625,  5.49137e-05)
    };

    \addplot [color=green,mark=triangle] plot coordinates {

        (.125,      0.0141393)
        (.0625,     0.00108092)
        (0.03125,   4.67021e-05)
        (0.015625,  2.58318e-06)
    };

    \addplot [color=red,mark=triangle, dashed] plot coordinates {

        (.25,       0.73789)
        (.125,      0.387214)
        (.0625,     0.195203)
        (0.03125,   0.0964845)
        (0.015625,  0.0476812)
    };

    \addplot [color=blue,mark=triangle, dashed] plot coordinates {

        (.25,       0.423619)
        (.125,      0.179079)
        (.0625,     0.0729287)
        (0.03125,   0.0286441)
        (0.015625,  0.0119787)
    };

    \addplot [color=black,mark=triangle, dashed] plot coordinates {

        (.25,       0.363333)
        (.125,      0.176097)
        (.0625,     0.0724466)
        (0.03125,   0.0299489)
        (0.015625,  0.0130237)
    };

    \addplot [color=green,mark=triangle, dashed] plot coordinates {

        (.25,       0.249916)
        (.125,      0.184512)
        (.0625,     0.0758842)
        (0.03125,   0.031917)
        (0.015625,  0.0142136)
    };

    \addplot [black] plot coordinates {

        (.03125,    0.08)
        (0.015625,    .04)
    } node[pos=0.5,anchor=north]{$1$};

    \addplot [black] plot coordinates {

        (.03125,    0.0040)
        (0.015625,    .001)
    } node[pos=0.5,anchor=north]{$2$};

    \addplot [black] plot coordinates {

        (.03125,    3.2000e-04)
        (0.015625,  4e-5)
    } node[pos=0.5,anchor=north]{$3$};

    \addplot [black] plot coordinates {

        (.03125,    2.4000e-05)
        (0.015625,  1.5e-6)
    } node[pos=0.5,anchor=north]{$4$};

    \legend{MG-$Q_2$\\MG-$Q_3$\\MG-$Q_4$\\MG-$Q_5$\\MB-$Q_2$\\MB-$Q_3$\\MB-$Q_4$\\MB-$Q_5$\\}
    \end{loglogaxis}
\end{tikzpicture}

%% file: tikz_scripts/five_patch_biharmonic_basic.tex
\begin{tikzpicture}
    \begin{loglogaxis}[
        legend columns=4,
    	legend style={at={(1,1.2)}, nodes={scale=0.8, transform shape}, column sep=.1cm},
        xlabel=$h$,
        ylabel=${\|u-u^{h}\|_{L^2}}/{\|u\|_{L^2}}$ 
    ]
    
    
    

    \addplot [color=red,mark=triangle] plot coordinates {

        (.25,   0.735913)
        (.125,   0.65318)
        (.0625,   0.0422797)
        (0.03125,  0.0081658)
        (0.015625,  0.001707)
    };

    \addplot [color=blue,mark=triangle] plot coordinates {

        (.25,   0.747898 )
        (.125,      0.0646759 )
        (.0625,     0.00581605)
        (0.03125,   0.000184301 )
        (0.015625,  5.67145e-06 )
    };

    \addplot [color=black,mark=triangle] plot coordinates {

        (.25,   0.173579)
        (.125,      0.0684303   )
        (.0625,     0.000237906   )
        (0.03125,   6.28963e-06  )
        (0.015625,  1.31465e-07  )
    };

    \addplot [color=green,mark=triangle] plot coordinates {

        (.25,   0.110367)
        (.125,   0.00351538)
        (.0625,   0.000117387)
        (0.03125,   7.25403e-07 )
        (0.015625,  1.54201e-08  )
    };

    \addplot [color=red,mark=triangle, dashed] plot coordinates {

        (.25,       0.743193)
        (.125,      0.192947)
        (.0625,     0.0377209)
        (0.03125,   0.00745985)
        (0.015625,  0.00164437 )
    };

    \addplot [color=blue,mark=triangle, dashed] plot coordinates {

        (.25,       0.502288)
        (.125,      0.0338857)
        (.0625,     0.00151253)
        (0.03125,   5.74269e-05 )
        (0.015625,  2.76672e-06 )
    };

    \addplot [color=black,mark=triangle, dashed] plot coordinates {

        (.25,       0.170459)
        (.125,      0.0202788)
        (.0625,     0.000197197)
        (0.03125,   3.71669e-06)
        (0.015625,  9.89612e-08)
    };

    \addplot [color=green,mark=triangle, dashed] plot coordinates {

        (.25,       0.0739105   )
        (.125,      0.00226322 )
        (.0625,     4.85118e-05 )
        (0.03125,       4.42907e-07)
        (0.015625,      1.67135e-08 )
    };

    \addplot [black] plot coordinates {

        (0.03125,    4.0000e-03)
        (0.015625,    .001)
    } node[pos=0.5,anchor=north]{$2$};

    \addplot [black] plot coordinates {

        (0.03125,    2.5600e-05)
        (0.015625,   1.6e-6)
    } node[pos=0.5,anchor=north]{$4$};

    \addplot [black] plot coordinates {

        (0.03125,    2.5600e-06)
        (0.015625,   8e-8)
    } node[pos=0.5,anchor=north]{$5$};

    \addplot [black] plot coordinates {

        (.0625,    1.9200e-05)
        (0.03125,  3e-7)
    } node[pos=0.5,anchor=north]{$6$};

    \legend{OG-$Q_2$\\OG-$Q_3$\\OG-$Q_4$\\OG-$Q_5$\\OB-$Q_2$\\OB-$Q_3$\\OB-$Q_4$\\OB-$Q_5$\\}
    \end{loglogaxis}
\end{tikzpicture}

%% file: tikz_scripts/five_patch_biharmonic_basic_H2.tex
\begin{tikzpicture}
    \begin{loglogaxis}[
        legend columns=4,
    	legend style={at={(1,1.2)}, nodes={scale=0.8, transform shape}, column sep=.1cm},
        xlabel=$h$,
        ylabel=${\|u_{}-u^{h}\|_{H^2_*}}/{\|u_{}\|_{H^2_*}}$ 
    ]

    \addplot [color=red,mark=triangle] plot coordinates {

        (.25,   0.933142)
        (.125,   0.7876)
        (.0625,   0.239011)
        (0.03125,   0.104946)
        (0.015625,  0.048573 )
    };

    \addplot [color=blue,mark=triangle] plot coordinates {

        (.25,       0.879058)
        (.125,      0.26212)
        (.0625,     0.0699817)
        (0.03125,   0.00952877)
        (0.015625,  0.00161309)
    };

    \addplot [color=black,mark=triangle] plot coordinates {

        (.25,       0.455764)
        (.125,      0.257802)
        (.0625,     0.00558481)
        (0.03125,   0.000631124)
        (0.015625,  6.37968e-05)
    };

    \addplot [color=green,mark=triangle] plot coordinates {

        (.25,       0.360695)
        (.125,      0.0253871)
        (.0625,     0.00310925)
        (0.03125,   8.67978e-05)
        (0.015625,  9.15269e-06)
    };

    \addplot [color=red,mark=triangle, dashed] plot coordinates {

        (.25,       0.905119)
        (.125,      0.502756)
        (.0625,     0.229302)
        (0.03125,   0.102333)
        (0.015625,  0.0480131)
    };

    \addplot [color=blue,mark=triangle, dashed] plot coordinates {

        (.25,       0.775019)
        (.125,      0.186196)
        (.0625,     0.0330064)
        (0.03125,   0.00621938)
        (0.015625,  0.00140077 )
    };

    \addplot [color=black,mark=triangle, dashed] plot coordinates {

        (.25,   0.449573)
        (.125,      0.0848983)
        (.0625,     0.00500784)
        (0.03125,   0.000486555)
        (0.015625,  5.75042e-05)
    };

    \addplot [color=green,mark=triangle, dashed] plot coordinates {

        (.25,       0.256437)
        (.125,      0.0164084)
        (.0625,     0.00113918)
        (0.03125,   7.17024e-05)
        (0.015625,  1.85343e-05)
    };

    \addplot [black] plot coordinates {

        (.03125,    0.08)
        (0.015625,    .04)
    } node[pos=0.5,anchor=north]{$1$};

    \addplot [black] plot coordinates {

        (.03125,    0.0040)
        (0.015625,    .001)
    } node[pos=0.5,anchor=north]{$2$};

    \addplot [black] plot coordinates {

        (.03125,    3.2000e-04)
        (0.015625,  4e-5)
    } node[pos=0.5,anchor=north]{$3$};

    \addplot [black] plot coordinates {

        (.0625,    9.6000e-04)
        (0.03125,   6e-5)
    } node[pos=0.5,anchor=north]{$4$};

    \legend{OG-$Q_2$\\OG-$Q_3$\\OG-$Q_4$\\OG-$Q_5$\\OB-$Q_2$\\OB-$Q_3$\\OB-$Q_4$\\OB-$Q_5$\\}
    \end{loglogaxis}
\end{tikzpicture}

%% file: tikz_scripts/five_patch_modify_projection_basic_H2.tex
\begin{tikzpicture}
    \begin{loglogaxis}[
        legend columns=4,
    	legend style={at={(1,1.2)}, nodes={scale=0.8, transform shape}, column sep=.1cm},
        xlabel=$h$,
        ylabel=${\|u_{}-u^{h}\|_{H^2_*}}/{\|u_{}\|_{H^2_*}}$ 
    ]

    \addplot [color=red,mark=triangle] plot coordinates {

        (.125,   0.374903)
        (.0625,   0.186565)
        (0.03125,  0.0932418)
        (0.015625,  0.0465267)
    };

    \addplot [color=blue,mark=triangle] plot coordinates {

        (.125,      0.124056)
        (.0625,     0.0258204)
        (0.03125,   0.00569548)
        (0.015625,  0.00136586)
    };

    \addplot [color=black,mark=triangle] plot coordinates {

        (.125,      0.0796233)
        (.0625,     0.00467479)
        (0.03125,   0.000462717)
        (0.015625,  5.49055e-05)
    };

    \addplot [color=green,mark=triangle] plot coordinates {

        (.125,      0.0141655)
        (.0625,     0.00108096)
        (0.03125,   4.66999e-05)
        (0.015625,  2.58817e-06)
    };

    \addplot [color=red,mark=triangle, dashed] plot coordinates {

        (.25,       0.734819)
        (.125,      0.375986)
        (.0625,     0.187824)
        (0.03125,   0.0938761)
        (0.015625,  0.046732)
    };

    \addplot [color=blue,mark=triangle, dashed] plot coordinates {

        (.25,       0.414403)
        (.125,      0.124012)
        (.0625,     0.0261069)
        (0.03125,   0.00572515)
        (0.015625,  0.00136939 )
    };

    \addplot [color=black,mark=triangle, dashed] plot coordinates {

        (.25,   0.354993)
        (.125,      0.0798194)
        (.0625,     0.00468961)
        (0.03125,   0.000464032)
        (0.015625,  5.49953e-05)
    };

    \addplot [color=green,mark=triangle, dashed] plot coordinates {

        (.25,       0.0958095)
        (.125,      0.0141841)
        (.0625,     0.00108168)
        (0.03125,   4.67261e-05)
        (0.015625,  2.58421e-06)
    };

    \addplot [black] plot coordinates {

        (.03125,    0.08)
        (0.015625,    .04)
    } node[pos=0.5,anchor=north]{$1$};

    \addplot [black] plot coordinates {

        (.03125,    0.0040)
        (0.015625,    .001)
    } node[pos=0.5,anchor=north]{$2$};

    \addplot [black] plot coordinates {

        (.03125,    3.2000e-04)
        (0.015625,  4e-5)
    } node[pos=0.5,anchor=north]{$3$};

    \addplot [black] plot coordinates {

        (.03125,    3.2000e-05)
        (0.015625,   2e-6)
    } node[pos=0.5,anchor=north]{$4$};

    \legend{MG-$Q_2$\\MG-$Q_3$\\MG-$Q_4$\\MG-$Q_5$\\MB-$Q_2$\\MB-$Q_3$\\MB-$Q_4$\\MB-$Q_5$\\}
    \end{loglogaxis}
\end{tikzpicture}

%% file: tikz_scripts/five_patch_projection_basic_H2.tex
\begin{tikzpicture}
    \begin{loglogaxis}[
        legend columns=4,
    	legend style={at={(1,1.2)}, nodes={scale=0.8, transform shape}, column sep=.1cm},
        xlabel=$h$,
        ylabel=${\|u_{}-u^{h}\|_{H^2_*}}/{\|u_{}\|_{H^2_*}}$ 
    ]

    
    
    

    \addplot [color=red,mark=triangle] plot coordinates {

        (.25,   0.92952)
        (.125,   0.786507)
        (.0625,   0.23874)
        (0.03125,   0.104899)
        (0.015625,  0.048567)
    };

    \addplot [color=blue,mark=triangle] plot coordinates {

        (.25,       0.877915)
        (.125,      0.261873)
        (.0625,     0.0698219)
        (0.03125,   0.0095109)
        (0.015625,  0.00161265)
    };

    \addplot [color=black,mark=triangle] plot coordinates {

        (.25,       0.455764)
        (.125,      0.257802)
        (.0625,     0.00558481)
        (0.03125,   0.000628332)
        (0.015625,  6.2662e-05)
    };

    \addplot [color=green,mark=triangle] plot coordinates {

        (.25,       0.359737)
        (.125,      0.0253871)
        (.0625,     0.00310925)
        (0.03125,   8.40285e-05)
        (0.015625,  7.11817e-06)
    };

    \addplot [color=red,mark=triangle, dashed] plot coordinates {

        (.25,       0.905119)
        (.125,      0.502756)
        (.0625,     0.229302)
        (0.03125,   0.102333)
        (0.015625,  0.0480131)
    };

    \addplot [color=blue,mark=triangle, dashed] plot coordinates {

        (.25,       0.775019)
        (.125,      0.186196)
        (.0625,     0.0330064)
        (0.03125,   0.00621938)
        (0.015625,  0.00140077 )
    };

    \addplot [color=black,mark=triangle, dashed] plot coordinates {

        (.25,   0.449573)
        (.125,      0.0848983)
        (.0625,     0.00500784)
        (0.03125,   0.000486555)
        (0.015625,  5.75042e-05)
    };

    \addplot [color=green,mark=triangle, dashed] plot coordinates {

        (.25,       0.256437)
        (.125,      0.0164084)
        (.0625,     0.00113918)
        (0.03125,   4.79673e-05)
        (0.015625,  2.60337e-06)
    };

    \addplot [black] plot coordinates {

        (.03125,    0.08)
        (0.015625,    .04)
    } node[pos=0.5,anchor=north]{$1$};

    \addplot [black] plot coordinates {

        (.03125,    0.0040)
        (0.015625,    .001)
    } node[pos=0.5,anchor=north]{$2$};

    \addplot [black] plot coordinates {

        (.03125,    3.2000e-04)
        (0.015625,  4e-5)
    } node[pos=0.5,anchor=north]{$3$};

    \addplot [black] plot coordinates {

        (.03125,    3.2000e-05)
        (0.015625,   2e-6)
    } node[pos=0.5,anchor=north]{$4$};

    \legend{OG-$Q_2$\\OG-$Q_3$\\OG-$Q_4$\\OG-$Q_5$\\OB-$Q_2$\\OB-$Q_3$\\OB-$Q_4$\\OB-$Q_5$\\}
    \end{loglogaxis}
\end{tikzpicture}

%% file: tikz_scripts/spectrum.tex
\begin{tikzpicture}
	\begin{axis}[no markers,
			xmin=0, xmax=1,
            ymin=1, ymax=3,
            xlabel=$n/N$,
            ylabel=$\omega_n^h/\omega_n$,
            width=0.6\textwidth,
            height=0.4\textwidth,
            legend style={at={(.4,.9)}, nodes={scale=1, transform shape}},]
	\addplot table [] {spectrum_c0.txt};
    \addplot table [] {spectrum_c1.txt};
    \legend{OB-$Q_5$ $C^0$\\OB-$Q_5$ $C^1$\\}
    \end{axis}
\end{tikzpicture}